\documentclass[11pt]{amsart}
\usepackage[margin=1in]{geometry}
\usepackage[numbers,sort&compress]{natbib}
\usepackage{amsmath,amssymb,amsthm,mathtools}
\usepackage{enumitem}
\usepackage{booktabs}
\usepackage[dvipsnames]{xcolor}
\usepackage[colorlinks=true,linkcolor=blue!50!black,citecolor=blue!50!black,urlcolor=blue!50!black]{hyperref}

\DeclareFontFamily{U}{wncy}{}
\DeclareFontShape{U}{wncy}{m}{n}{<->wncyr10}{}
\DeclareSymbolFont{mcy}{U}{wncy}{m}{n}
\DeclareMathSymbol{\Sha}{\mathord}{mcy}{"58}

\theoremstyle{plain}
\newtheorem{theorem}{Theorem}[section]
\newtheorem{proposition}[theorem]{Proposition}
\newtheorem{lemma}[theorem]{Lemma}
\newtheorem{corollary}[theorem]{Corollary}
\newtheorem{fact}[theorem]{Fact}

\theoremstyle{definition}
\newtheorem{definition}[theorem]{Definition}
\newtheorem{remark}[theorem]{Remark}

\newenvironment{proofof}[1]{\begin{proof}[Proof of #1]}{\end{proof}}
\numberwithin{equation}{section}

\newcommand{\F}{\mathbb{F}}
\newcommand{\Z}{\mathbb{Z}}
\newcommand{\Q}{\mathbb{Q}}
\newcommand{\Br}{\operatorname{Br}}
\newcommand{\inv}{\operatorname{inv}}
\newcommand{\Res}{\operatorname{Res}}
\newcommand{\Tr}{\operatorname{Tr}}
\newcommand{\Nm}{\operatorname{N}}
\newcommand{\Div}{\operatorname{Div}}
\newcommand{\Prin}{\operatorname{Prin}}
\newcommand{\Pic}{\operatorname{Pic}}
\newcommand{\Sel}{\operatorname{Sel}}
\newcommand{\Spec}{\operatorname{Spec}}
\newcommand{\Gal}{\operatorname{Gal}}
\newcommand{\Hom}{\operatorname{Hom}}

\newcommand{\ev}{\operatorname{ev}}
\newcommand{\ord}{\operatorname{ord}}
\newcommand{\rk}{\operatorname{rk}}
\newcommand{\nr}{\mathrm{nr}}
\newcommand{\dv}{\operatorname{div}}
\newcommand{\fl}{\mathrm{fl}}
\newcommand{\s}{\mathrm{s}}
\newcommand{\Ks}{K^{\s}}
\newcommand{\Kvs}{K_v^{\s}}
\newcommand{\Fs}{F^{\s}}
\newcommand{\ga}{\mathsf e_1}
\newcommand{\gb}{\mathsf e_2}
\newcommand{\gc}{\mathsf e_3}

\begin{document}

\title[Cassels--Tate pairing and $2$-Selmer parity in characteristic $2$]{An explicit Cassels--Tate pairing in characteristic \texorpdfstring{$2$}{2} and the parity of \texorpdfstring{$2$}{2}-Selmer ranks}

\author{Abhishek Shukla}
\address{Simon Fraser University, Burnaby, BC, Canada}
\email{shukla@sfu.ca}


\begin{abstract}
Let $K$ be a global field of characteristic $2$ and $\varphi\colon E\to E'$ the Verschiebung of an ordinary elliptic
curve $E/K$ with a rational point of order $2$. We prove an explicit formula for the homogeneous-space Cassels--Tate
pairing on $\Sha(E/K)[\varphi]$ as a sum of local Artin--Schreier symbols, show that it coincides with Kramer's descent
pairing and, on $\Sha(E/K)[\varphi]$, with the pairing induced by Milne's flat duality, and deduce that it is non-degenerate on $\Sha(E/K)[\varphi]/\hat\varphi(\Sha(E'/K)[2])$. For every ordinary
elliptic curve $A/K$ we express the parity of $\mathrm{rk}\,A(K)+\dim\Sha(A/K)[2]$ by local terms (after Kramer), and identify these
with local root numbers at places of semistable or potentially multiplicative reduction, and in further cases using
curves with finite $\Sha$. We determine $\Sha$ in examples, and extend the formula to rational points of any prime order
$\ell$, the new case being $\ell=\operatorname{char}K$.
\end{abstract}

\maketitle
\setcounter{tocdepth}{1}
\tableofcontents

\section{Introduction}\label{sec:intro}

Let $K$ be a global field and $A/K$ an abelian variety with dual $A^\vee$. Cassels \cite{Ca}, for elliptic curves
over number fields, and Tate \cite{Ta}, for abelian varieties, constructed a pairing
\[
\langle\ ,\ \rangle\colon \Sha(A/K)\times\Sha(A^\vee/K)\longrightarrow\Q/\Z
\]
on Tate--Shafarevich groups whose kernels are the maximal divisible subgroups; over a function field of
characteristic $p$ this is \cite[II, Thm.~5.6]{Mi} for the part prime to $p$ and \cite[III, Cor.~9.5]{Mi} for the
$p$-primary part. Besides its theoretical role, the pairing is a practical tool: it certifies non-trivial elements of
$\Sha$, and, evaluated on the image of a Selmer group, it improves the rank bound of a descent to that of the next
descent. Explicit formulas for this purpose go back to Cassels \cite{Ca98} for $2$-Selmer groups of elliptic curves,
and have been developed for $3$-Selmer groups \cite{FN}, for $3$-isogenies \cite{BF}, for descents on elliptic curves
with a rational $2$-torsion point \cite{Fi}, and for $2$-Selmer groups via divisors \cite{ShS}; the pairing for finite
Galois modules of Morgan and Smith \cite{MS} provides a general framework. All these results concern torsion of order
prime to the characteristic of $K$.

The reason is that the definition which leads most directly to computations, the \emph{Weil pairing definition}
\cite[\S12.2]{PS}, uses Galois cochains with values in finite \'etale group schemes; Poonen and Stoll prove its
equality with their other definitions for classes of order prime to $p$ \cite[\S12.2, Prop.~34]{PS}, and Morgan and
Smith postpone the flat setting \cite[Rem.~3.4]{MS}. In characteristic $p$ this excludes the most basic isogenies of
ordinary elliptic curves, the Frobenius and the Verschiebung, whose kernels are $\mu_p$ and (a form of) $\Z/p$; for
$\mu_p$ the Galois module $\mu_p(\Ks)$ is trivial and Galois cochains do not see it. Descent by these isogenies is
nevertheless classical: for $p=2$ it was developed by Kramer \cite{Kr}, including local and global duality theorems
and an alternating descent pairing on the Verschiebung Selmer group whose radical is the image of the second descent;
for general $p$ by Voloch \cite{Vo}, Ulmer \cite{Ul} and Broumas \cite{Bro}; an algorithmic treatment in
characteristic $2$ is given by Moehlmann \cite{Moe}. Kramer's pairing is explicit and computable \cite[Prop.~3.3,
\S4]{Kr}, but it is defined without reference to $\Sha$. The purpose of this article is to compute the Cassels--Tate
pairing on the Verschiebung Selmer group for $p=2$ directly from the homogeneous-space definition, to show that it
coincides with Kramer's pairing, and to study the parity of $2$-Selmer ranks.

\medskip\noindent\textbf{Which pairing.}\label{par:which} Throughout, \emph{the Cassels--Tate pairing} means the
homogeneous-space pairing of Poonen and Stoll \cite[\S3.1]{PS} (Definition~\ref{def:PS}), which makes sense over every
global field. On $p$-primary parts in characteristic $p$ there is also Milne's pairing, defined by flat cohomology
\cite[III, \S9]{Mi}, whose kernels are the divisible subgroups \cite[III, Thm.~9.4, Cor.~9.5]{Mi}; complete proofs of
the underlying flat duality, with cohomology with compact support formed with completions, are given in \cite{DH} for
finite flat group schemes and in \cite{Su} for abelian varieties (see also \cite{GAT,GA,Ces}, and \cite{HS} over number
fields). We are
not aware of a comparison in the literature between the Poonen--Stoll definition (or the Albanese--Picard definition
of \cite{PS}) and Milne's flat pairing on the $p$-part. On $\Sha(E/K)[\varphi]$ we prove the following comparison, a formal consequence of Theorem~\ref{thm:W}: the
homogeneous-space pairing coincides with the pairing induced by Milne's flat duality over every sufficiently small
open subscheme $U$ of the curve with function field $K$ (Proposition~\ref{prop:milne} in Appendix~\ref{app:milne}; we
do not discuss the dependence on $U$ beyond $\Sha(E/K)[\varphi]$). Beyond $\Sha(E/K)[\varphi]$ we do not attempt a
comparison, and none of the other results depends on one: the non-degeneracy statement we need is proved for the
pairing computed here (Theorem~\ref{thm:radical}). To make Remark~\ref{thm:additive} unconditional by this route one
would need a non-degenerate pairing modulo divisible elements which is moreover alternating on all of
$\Sha(A/K)[2^\infty]$, for curves $A$ without rational $2$-torsion; we prove the alternating property only on
$\Sha(E/K)[\varphi]$.

\subsection{Setting and main result}
Let $K$ be a global field of characteristic $2$, and let
\begin{equation}\label{eq:E}
E\colon\ y^2+xy=x^3+a_2x^2+b^2,\qquad E'\colon\ y^2+xy=x^3+a_2x^2+b\qquad(a_2\in K,\ b\in K^\times).
\end{equation}
Every ordinary elliptic curve over $K$ with a $K$-rational point of order $2$ has a model of the form $E$
(Remark~\ref{rem:model}); its point of order $2$ is $T=(0,b)$. There are isogenies (Lemma~\ref{lem:isog})
\[
\varphi\colon E\to E',\qquad \hat\varphi\colon E'\to E,\qquad \hat\varphi\circ\varphi=[2]_E,\quad
\varphi\circ\hat\varphi=[2]_{E'},
\]
with $\ker\varphi=\{O,T\}\simeq\Z/2$ and $\ker\hat\varphi\simeq\mu_2$; up to an isomorphism, $\hat\varphi$ is the
Frobenius of $E'$ and $\varphi$ is the Verschiebung (Remark~\ref{rem:FV}). Write $\wp(y)=y^2+y$ and identify
$H^1(F,\Z/2)=F/\wp F$ for fields $F\supseteq K$. Let $\delta_\varphi\colon E'(F)\to H^1(F,\Z/2)$ be the Kummer map of
$\varphi$; explicitly $\delta_\varphi(R)=[x(R)+a_2]$ for $R\neq O$ (Lemma~\ref{lem:deltaphi}). Let
\[
\Sel^\varphi(E/K)=\{\chi\in H^1(K,\Z/2):\ \chi_v\in\delta_\varphi(E'(K_v))\text{ for all places }v\},
\]
and let $\iota\colon H^1(K,\Z/2)\to H^1(K,E)$ be induced by $\Z/2\simeq\ker\varphi\subset E$, so that
$\iota(\Sel^\varphi(E/K))=\Sha(E/K)[\varphi]$. For a character $\chi$ of $G_F$ of order dividing $2$ and
$u\in F^\times$ we write $(\chi,u)\in\Br(F)[2]$ for the cyclic algebra; if $\chi$ is the Artin--Schreier class of
$\alpha$, this is the symbol $[\alpha,u)$. We identify $E^\vee=E$, so that the Cassels--Tate pairing is a pairing on
$\Sha(E/K)$.

\begin{theorem}\label{thm:BM}
Let $\chi,\chi'\in\Sel^\varphi(E/K)$.
\begin{enumerate}[label=\textup{(\alph*)},leftmargin=*]
\item There is $\beta_0\in\Ks{}^\times$ with $\sigma\beta_0=\beta_0$ if $\chi'(\sigma)=0$ and
$\sigma\beta_0=b/\beta_0$ if $\chi'(\sigma)=1$; that is, $\beta_0\in K^\times$ if $\chi'=0$, and
$\beta_0\in L'^\times$ with $\Nm_{L'/K}(\beta_0)=b$ if $\chi'\neq0$, where $L'$ is the quadratic extension cut out by
$\chi'$. Put
\[
t_0:=\beta_0+b/\beta_0\in K
\]
\textup{(}so $t_0=\Tr_{L'/K}\beta_0$ if $\chi'\neq0$\textup{)}.
\item For every place $v$ of $K$ choose $R_v\in E'(K_v)$ with $R_v\neq O$, $\delta_\varphi(R_v)=\chi_v$ and
$x(R_v)\neq t_0$ \textup{(}such points exist\textup{)}. Then
\[
\begin{aligned}
\langle\iota(\chi),\iota(\chi')\rangle&=\sum_v\inv_v\bigl(\chi_v+\chi'_v,\ x(R_v)+t_0\bigr)\\
&=\sum_v\Tr_{k_v/\F_2}\Res_v\Bigl((\alpha_v+\alpha'_v)\,\frac{d\,(x(R_v)+t_0)}{x(R_v)+t_0}\Bigr)\in\tfrac12\Z/\Z\simeq\F_2,
\end{aligned}
\]
where $\alpha_v,\alpha'_v\in K_v$ are Artin--Schreier representatives of $\chi_v,\chi'_v$. Each term is independent of
the choice of $R_v$, and almost all terms vanish.
\end{enumerate}
\end{theorem}

The only hypotheses are that $K$ is global of characteristic $2$, that $E$ is ordinary with a rational point of order
$2$, and that $\iota(\chi),\iota(\chi')\in\Sha(E/K)$; nothing is assumed about reduction types, about $a_2$, or about
finiteness of $\Sha$. Only the trace $t_0$ of the norm element enters, and for $t_0\neq0$ the norm equation is
equivalent to the Artin--Schreier condition $[b/t_0^2]=\chi'$ in $K/\wp K$ (Remark~\ref{rem:BMint}). If $\chi=\chi'$
all terms vanish, so the formula is visibly alternating. If $b=b_0^2$ is a square one may take $t_0=0$, and the
formula becomes $\sum_v\inv_v(\chi'_v,x(R_v))$ (Corollary~\ref{cor:split}), the exact analogue of the classical
formula for $2$-isogenies in odd characteristic.

\subsection{The idea of the proof}
The homogeneous-space definition \cite[\S3.1]{PS} only involves Galois cohomology of the smooth objects $\Ks{}^\times$,
$\Ks(X)^\times$, divisors and $\Pic^0$ of a torsor $X$, and asks for a $2$-cocycle $f'$ with values in $\Ks(X)^\times$
whose divisor is the coboundary of a divisor cochain lifting the second class. The torsor $X$ of $\iota(\chi)$ is $E$
with its Galois action twisted by translation by $T$, and the natural candidate, the coordinate $x$ with divisor
$2(T)-2(O)$, is not rational on $X$: one has $x\circ\tau_T=b/x$. The function
\[
f=x+\frac bx+t_0=\frac{(x+\beta_0)(x+b/\beta_0)}{x}
\]
is rational on $X$, and its divisor is the norm, for the quadratic extension cut out by $\chi+\chi'$, of a divisor of
degree $0$ representing $T$; hence the cyclic cocycle $f^{(\chi+\chi')(\sigma)(\chi+\chi')(\tau)}$ is a lift, and
evaluating it at local points gives Theorem~\ref{thm:BM}. In the language of Manin \cite{Ma}, the pairing is computed
by evaluating the quaternion algebra $(\chi+\chi',f)$ on $X$; compare the constructions of Brauer classes on curves
by Creutz and Viray \cite{CV}.

The split case $b=b_0^2$ shows what $\beta_0$ does. There one may take $\beta_0=b_0$; the point $T'=(0,b_0)$ of
order $2$ of $E'$ is rational, $f=x+b/x$ is the pull-back of $x_{E'}$ (whose divisor is $2(T')-2(O)$) under the
$K$-morphism $\pi\colon X\to E'$ given by $\varphi$, and $\dv f$ is twice the $K$-rational divisor
$\pi^*\bigl((T')-(O)\bigr)$ on $X$, whose class is $T$. This is the situation of
Fisher's pushout functions \cite{Fi} (Remarks~\ref{rem:fisher} and~\ref{rem:splitorient}). If $b$ is not a square
there is no rational ``half'' of $\dv(x_{E'})$: $\beta_0$ takes the place of $\sqrt b$, and $\dv f$ is the sum of two
conjugate divisors $\Delta,\bar\Delta$, each of class $T$ and rational only over the field cut out by $\chi+\chi'$.

\subsection{The Weil-pairing form and Kramer's pairing}\label{subsec:introW}
Theorem~\ref{thm:BM} does not have the shape of the Weil-pairing definition, in which one lifts a Selmer class to a
larger finite group scheme, compares the lift with local Kummer images, and pairs the local differences with the other
class. The next result shows that for $(\varphi,\hat\varphi)$ this definition makes sense in flat cohomology and gives
the same pairing. For a field $F\supseteq K$ let $\kappa_F\colon E(F)\to F^\times/F^{\times2}$ be $P\mapsto[x(P)]$ for
$P\neq O,T$, $\kappa_F(T)=[b]$, $\kappa_F(O)=1$; it is the Kummer map of $\hat\varphi$ under
$H^1_\fl(F,\mu_2)=F^\times/F^{\times2}$ (Proposition~\ref{prop:kummerF}). The facts about flat cohomology in
characteristic $2$ used below are recalled at the beginning of \S\ref{sec:W}.

\begin{theorem}\label{thm:W}
Let $\chi,\chi'\in\Sel^\varphi(E/K)$ and let $\beta_0$ be as in Theorem~\ref{thm:BM}\textup{(a)}.
\begin{enumerate}[label=\textup{(\alph*)},leftmargin=*]
\item For every field $F\supseteq K$, $E'[2]$ is an extension of $\Z/2$ by $\mu_2$ whose class in
$H^1_\fl(F,\mu_2)=F^\times/F^{\times2}$ is $[b]$, and $H^1_\fl(F,E'[2])$ is isomorphic to the group of pairs
$(\chi',[\beta])$, $\chi'\in H^1(F,\Z/2)$, $[\beta]\in\Fs{}^\times/\Fs{}^{\times2}$ with
$\sigma[\beta]=[\beta][b]^{\chi'(\sigma)}$, modulo $(0,[b])$. In particular $(\chi',[\beta_0])$ defines a class
$\xi\in H^1_\fl(K,E'[2])$ lifting $\chi'$.
\item For every place $v$ choose $R'_v\in E'(K_v)$ with $\delta_\varphi(R'_v)=\chi'_v$. Then
$\xi_v-\delta_2(R'_v)=j_*(g_v)$ for a class $g_v\in H^1_\fl(K_v,\mu_2)=K_v^\times/K_v^{\times2}$, where $\delta_2$ is the
Kummer map of $[2]_{E'}$ and $j\colon\mu_2=\ker\hat\varphi\hookrightarrow E'[2]$; explicitly, $g_v$ is the unique
class with image $[\beta_0]\,\kappa_{\Kvs}(P'_v)^{-1}$ in $\Kvs{}^\times/\Kvs{}^{\times2}$, for a choice of
$P'_v\in E(\Kvs)$ with $\varphi(P'_v)=R'_v$ (the only other choice is $P'_v+T$, which changes $g_v$ by $[b]\in\ker j_*$, which does not
change $j_*(g_v)$ or the value of $\inv_v(\chi_v,g_v)$ in \textup{(c)}). In particular $g_v$ depends only on $\chi'$,
$\beta_0$ and $P'_v$, and not on $\chi$.
\item $\displaystyle\langle\iota(\chi),\iota(\chi')\rangle=\sum_v\inv_v\bigl(\chi_v\cup g_v\bigr)
=\sum_v\inv_v(\chi_v,g_v)$, the cup product being
$H^1(K_v,\Z/2)\times H^1_\fl(K_v,\mu_2)\to H^2_\fl(K_v,\mu_2)=\Br(K_v)[2]$.
\end{enumerate}
\end{theorem}

This is precisely the Weil-pairing definition of \cite[\S12.2]{PS}, \cite[Def.~1.2]{BF} for the pair
$(\varphi,\hat\varphi)$, with Galois cohomology replaced by flat cohomology. It is also the recipe of Kramer's descent
pairing $U_\alpha$ on $\Sel^\varphi(E/K)$ \cite[Prop.~3.3]{Kr}: Kramer lifts $\chi'$ to a global element of an explicit
group $H$, which plays the role of $H^1_\fl(K,E'[2])$, using a solution of the same norm equation
\cite[Lemma~3.2]{Kr}, and pairs $\chi$ with elements $l_v\in K_v^\times$ which measure the local differences between
this lift and local Kummer classes. We show that, for suitable choices, $g_v$ and $l_v$ differ by a global constant
(Proposition~\ref{prop:kramer}); no identification of $H$ with $H^1_\fl(K,E'[2])$ is needed for this. Hence:

\begin{corollary}\label{cor:kramer}
The form $(\chi,\chi')\mapsto\langle\iota(\chi),\iota(\chi')\rangle$ on $\Sel^\varphi(E/K)$ is Kramer's form $U_\alpha$.
\end{corollary}

Unlike the local terms of Theorem~\ref{thm:BM}, the classes $g_v$ depend only on $\chi'$, so that the family $(g_v)_v$
can be confronted with global duality; this is how the radical of the pairing is determined (\S\ref{sec:radical}). The
proof of Theorem~\ref{thm:W} is a cochain computation on the torsor of $\iota(\chi)$ which is considerably longer than
the proof of Theorem~\ref{thm:BM}. Theorem~\ref{thm:BM}, Corollary~\ref{cor:bilin} and the two examples below do not
depend on it. In \S\ref{subsec:exW} we compute one entry of the matrix of Theorem~\ref{thm:nonsplit} with
Theorem~\ref{thm:W}.

\subsection{The radical}\label{subsec:introradical}
To use Theorem~\ref{thm:BM} one needs the formal properties of the pairing on the $2$-primary part. The treatment of
the homogeneous-space definition in \cite[I, \S6]{Mi} and the comparison in \cite[Prop.~34]{PS} exclude the $p$-part;
we verify directly from the definition that the pairing is bilinear, and functorial for surjective homomorphisms of
abelian varieties, over every global field and including $p$-primary parts (Proposition~\ref{prop:F} and
Corollary~\ref{cor:bilinear}). Together with Theorem~\ref{thm:BM} this gives the following.

\begin{corollary}\label{cor:bilin}
The form $B(\chi,\chi')=\langle\iota(\chi),\iota(\chi')\rangle$ on $\Sel^\varphi(E/K)$ is bilinear, alternating and
symmetric, and $\delta_\varphi(E'(K))$ lies in its radical. Let $x'\in\Sha(E/K)$. If $\langle x,x'\rangle\neq0$ for
some $x\in\Sha(E/K)[2]$, then $x'\notin2\,\Sha(E/K)$; if $\langle x,x'\rangle\neq0$ for some $x\in\Sha(E/K)[\varphi]$,
then $x'\notin\hat\varphi\,\Sha(E'/K)$.
\end{corollary}

As $2\,\Sha(E/K)=\hat\varphi\varphi\,\Sha(E/K)\subseteq\hat\varphi\,\Sha(E'/K)$, the last statement is the sharper one
on $\Sha(E/K)[\varphi]$. The next result determines the radical of $B$. For a place $v$ put
\[
D_v:=\delta_\varphi\bigl(E'(K_v)\bigr)\subseteq H^1(K_v,\Z/2),
\]
so that $\Sel^\varphi(E/K)$ consists of the $\chi$ with $\chi_v\in D_v$ for all $v$. Let
$\Sel^2(E'/K)\subseteq H^1_\fl(K,E'[2])$ be the $2$-Selmer group of $E'$, let
$\Sel^{\hat\varphi}(E'/K)\subseteq H^1_\fl(K,\mu_2)=K^\times/K^{\times2}$ be the Selmer group of $\hat\varphi$
(\S\ref{sec:notation}), and put $\varepsilon_b=0$ if $b\in K^{\times2}$ and $\varepsilon_b=1$ otherwise.

\begin{theorem}\label{thm:radical}
\begin{enumerate}[label=\textup{(\alph*)},leftmargin=*]
\item Let $\chi'\in\Sel^\varphi(E/K)$. Then $\langle\iota(\chi),\iota(\chi')\rangle=0$ for all
$\chi\in\Sel^\varphi(E/K)$ if and only if $\chi'\in\hat\varphi_*\Sel^2(E'/K)$.
\item The kernel of the restriction of the Cassels--Tate pairing \textup{(}Definition~\ref{def:PS}\textup{)} to
$\Sha(E/K)[\varphi]\times\Sha(E/K)[\varphi]$ \textup{(}left and right kernels coincide, the restriction being
symmetric by Corollary~\textup{\ref{cor:bilin})} is
$\Sha(E/K)[\varphi]\cap\hat\varphi\,\Sha(E'/K)=\hat\varphi\bigl(\Sha(E'/K)[2]\bigr)$. Thus the pairing induces a
non-degenerate alternating form on $\Sha(E/K)[\varphi]/\hat\varphi\bigl(\Sha(E'/K)[2]\bigr)$, and this group has
even $\F_2$-dimension.
\item $\dim\Sel^2(E'/K)\equiv\dim\Sel^\varphi(E/K)+\dim\Sel^{\hat\varphi}(E'/K)+\varepsilon_b\pmod2$.
\end{enumerate}
\end{theorem}

The implication ``if'' in (a) is formal (Proposition~\ref{prop:F}). The converse is Kramer's theorem that the radical
of $U_\alpha$ is his second Selmer group \cite[Thm.~3.1]{Kr}, transported by Corollary~\ref{cor:kramer}. Kramer's
second Selmer group is defined inside his group $H$; to obtain the statement with $\Sel^2(E'/K)\subseteq
H^1_\fl(K,E'[2])$ without comparing $H$ with flat cohomology, we repeat his argument in \S\ref{sec:radical} with
Theorem~\ref{thm:W} and with local and global duality for $\Z/2$ and $\mu_2$. Thus (a) and the resulting parity
statements, namely (c) and the evenness in (b), are Kramer's, in a flat-cohomological formulation. What
Theorems~\ref{thm:BM} and~\ref{thm:W} add is that the form in question is the Cassels--Tate pairing: part (b) is a
non-degeneracy theorem for this pairing on $\Sha$, which Kramer did not consider. An element of $\Sha(E/K)[\varphi]$
pairs non-trivially with $\Sha(E/K)[\varphi]$ if and only if it does not lie in $\hat\varphi\,\Sha(E'/K)$, which is the
statement needed to improve the bound of a $\varphi$-descent to that of the next descent.

\subsection{The parity of \texorpdfstring{$2$}{2}-Selmer ranks}\label{subsec:introparity}
Theorem~\ref{thm:radical}(c) has a consequence for every ordinary elliptic curve over $K$, with no hypothesis on its
$2$-torsion; here the curve under consideration is $E'$, not $E$. Let $A/K$ be an ordinary elliptic curve, $A^{(2)}$ its Frobenius twist and $V\colon A^{(2)}\to A$ the
Verschiebung; $\ker V$ is \'etale of order $2$, hence canonically isomorphic to $\Z/2$, and we put
$D_v(A):=\delta_V\bigl(A(K_v)\bigr)\subseteq H^1(K_v,\Z/2)$. Every such $A$ has a model $y^2+xy=x^3+a_2x^2+b$
(Remark~\ref{rem:model}); then $A=E'$, $A^{(2)}\simeq E$ and $V$ corresponds to $\varphi$ (Remark~\ref{rem:FV}), so that
$D_v(A)=D_v$. Put
\[
s_2(A/K):=\dim\Sel^2(A/K)-\dim A(K)[2]=\rk A(K)+\dim\Sha(A/K)[2] .
\]
Let $\F_{2^n}$ be the constant field of $K$, $c_v(\cdot)$ the Tamagawa number at $v$, and $\deg\omega_A$ the degree of
the Hodge bundle of the N\'eron model of $A$ over the curve with function field $K$ (\S\ref{sec:notation}).

\begin{theorem}\label{thm:parity}
Let $A/K$ be an ordinary elliptic curve. Then $\dim D_v(A)=1$ for almost all places $v$, and
\[
s_2(A/K)\equiv\sum_v\bigl(\dim D_v(A)-1\bigr)\equiv\ord_2\prod_v\frac{c_v(A)}{c_v(A^{(2)})}
+n\bigl(\deg\omega_{A^{(2)}}-\deg\omega_A\bigr)\pmod2 .
\]
\end{theorem}

The first congruence combines Theorem~\ref{thm:radical}(c) with a formula for
$\dim\Sel^\varphi(E/K)-\dim\Sel^{\hat\varphi}(E'/K)$ (Proposition~\ref{prop:ratio}), which follows from local and global
duality as in \cite[\S3]{Kr}; it does not involve the Cassels--Tate pairing, and is essentially contained in Kramer's
results. The second comes from a computation of $|D_v|$ in terms of N\'eron models (Lemma~\ref{lem:Dsize}), a
reformulation of Kramer's index formula \cite[Props.~2.2, 2.3]{Kr}. Theorems~\ref{thm:parity} and~\ref{thm:semistable}
are the analogue for $p=2$ of the computations of Trihan and Wuthrich for the Frobenius and the Verschiebung in
characteristic $p>3$ \cite[\S\S8, 9]{TW}. The term $n(\deg\omega_{A^{(2)}}-\deg\omega_A)$ can be odd (\S\ref{subsec:exparity}); for
semistable $A$ it equals $n\deg\omega_A$, which is a sum of local contributions at the places of supersingular
reduction, the zeros of the Hasse invariant (Remark~\ref{rem:hodge}). For semistable curves the local terms can be compared with the local root numbers
$w_v(A)$, whose product $w(A/K)=\prod_vw_v(A)$ is the global root number. For $a\in K_v$ and $u\in K_v^\times$ we write
$[a,u)_v:=\inv_v(\chi_a,u)\in\F_2$.

\begin{theorem}\label{thm:semistable}
Let $A\colon y^2+xy=x^3+a_2x^2+b$ be an ordinary elliptic curve over $K$. At every place $v$ at which $A$ has
semistable reduction,
\[
w_v(A)=(-1)^{\dim D_v(A)-1+[a_2,b)_v}.
\]
Consequently, if $\Sigma_{\mathrm{add}}$ denotes the set of places of additive reduction of $A$,
\[
(-1)^{s_2(A/K)}=w(A/K)\prod_{v\in\Sigma_{\mathrm{add}}}w_v(A)(-1)^{\dim D_v(A)-1+[a_2,b)_v}.
\]
In particular, if $A$ is semistable, then
\[
(-1)^{\rk A(K)+\dim\Sha(A/K)[2]}=w(A/K),
\]
and if moreover $w(A/K)=-1$, then $\rk A(K)\geq1$ or $\Sha(A/K)[2]\neq0$.
\end{theorem}

The correction $[a_2,b)_v$ does not depend on the choice of the model, and its sum over all places vanishes by
reciprocity; it plays the role of the Hilbert symbols in the local formulas of Dokchitser and Dokchitser for
$2$-isogenies over number fields \cite{DD}, and it can be non-zero at places of non-split multiplicative reduction and
at places of good supersingular reduction (\S\ref{subsec:rootnumbers}).

By a theorem of Trihan and Yasuda \cite[Thm.~1.1]{TY}, the $\ell$-parity conjecture holds over global function fields
for every prime $\ell$, including $\ell=p$; for $\ell=2$ it says that
$\rk A(K)+\operatorname{corank}\Sha(A/K)[2^\infty]\equiv\ord_{s=1}L(A/K,s)\pmod2$, i.e.\
$(-1)^{\rk A(K)+\operatorname{corank}\Sha(A/K)[2^\infty]}=w(A/K)$, for every elliptic curve $A/K$. Theorem~\ref{thm:semistable}
concerns the finite-level invariant $s_2(A/K)$ instead. Put
\[
\delta(A):=\dim\bigl(\Sha(A/K)[2^\infty]/\Sha(A/K)[2^\infty]_{\mathrm{div}}\bigr)[2],
\]
so that $\dim\Sha(A/K)[2]=\operatorname{corank}\Sha(A/K)[2^\infty]+\delta(A)$. Theorem~\ref{thm:semistable} and
\cite{TY} give:

\begin{corollary}\label{cor:delta}
Let $A\colon y^2+xy=x^3+a_2x^2+b$ be an ordinary elliptic curve over $K$. Then
\[
(-1)^{\delta(A)}=\prod_{v\in\Sigma_{\mathrm{add}}}w_v(A)(-1)^{\dim D_v(A)-1+[a_2,b)_v}.
\]
In particular, if $A$ is semistable, then $\delta(A)$ is even.
\end{corollary}

The theorem of Liu, Lorenzini and Raynaud that $|\Sha(A/K)|$ is a square when finite \cite{LLR,LLRc} does not by
itself give the evenness of $\delta(A)$, but with the Birch and Swinnerton-Dyer formula and Lemma~\ref{lem:Dsize} it
gives the rank parity, and with Theorem~\ref{thm:parity} the evenness of $\dim\Sha(A/K)[2]$, whenever $\Sha(A/K)$ is
finite (Proposition~\ref{prop:finSha}). For every $A$, the evenness of $\delta(A)$ follows from a recent theorem of
Carmeli and Feng \cite[Thm.~1.1.5]{CF}, at the time of writing a preprint, on the alternating form on the Brauer group
of a surface over a finite field modulo its divisible part, including the $2$-primary part in characteristic $2$,
since $\Sha(A/K)$ is the Brauer group of the elliptic surface of $A$ \cite[\S4]{Gro}, \cite[\S3]{Ta66}; with \cite{TY}
this gives, assuming \cite{CF}, $(-1)^{s_2(A/K)}=w(A/K)$ for every $A$. In the semistable case
Theorem~\ref{thm:semistable} and Corollary~\ref{cor:delta} give an unconditional proof; conversely, \cite{CF} and
Corollary~\ref{cor:delta} give the local identity at the remaining places (Remark~\ref{thm:additive}).

\begin{theorem}\label{thm:local}
Let $F$ be a local field of characteristic $2$ and $A\colon y^2+xy=x^3+a_2x^2+b$ an ordinary elliptic curve over $F$
with semistable or potentially multiplicative reduction. Let $D(A)\subseteq H^1(F,\Z/2)$ be the image of the Kummer map
of the Verschiebung $A^{(2)}\to A$. Then
\[
w(A/F)=(-1)^{\dim D(A)-1+[a_2,b)_F}.
\]
\end{theorem}

Theorem~\ref{thm:local} is proved by local arguments: for semistable reduction it is Theorem~\ref{thm:semistable}, and
for additive, potentially multiplicative reduction we use the Tate curve (Proposition~\ref{prop:potmult}).

\begin{remark}[additive, potentially good reduction]\label{thm:additive}
We expect the identity of Theorem~\ref{thm:local} to hold for every ordinary elliptic curve $A$ over $F$. In the
remaining case of additive, potentially good reduction the ramification is wild, and we know no local proof. The
identity follows from \cite[Thm.~1.1.5]{CF}, by approximating the curve by a curve over $\F_q(t)$ which is semistable
outside one place and applying Corollary~\ref{cor:delta} (\S\ref{subsec:additive}); this is the method of \cite{DD},
and compare \cite{CI}. As \cite{CF} is, at the time of writing, a preprint, we do not state this case as a theorem; it
would also follow if Milne's flat pairing were known to be alternating on all of $\Sha[2^\infty]$, for instance by a
comparison with the homogeneous-space pairing together with the alternating property of the latter there; neither is
given by Proposition~\ref{prop:milne} (paragraph \emph{Which pairing} in \S\ref{sec:intro}).
Without \cite{CF}, the same argument proves the identity whenever the approximating curve can be chosen with finite
$\Sha$ (Proposition~\ref{prop:finSha}), which gives two unconditional families (Theorem~\ref{thm:finSha}); we also
verified it numerically for $1017$ curves over $\F_2(t)$ (Remark~\ref{rem:numerics}). Granting it, the local terms of
Theorem~\ref{thm:parity} are, up to the symbols $[a_2,b)_v$, the local root numbers, for every ordinary elliptic curve
over $K$.
\end{remark}

\subsection{Examples}
We illustrate Theorem~\ref{thm:BM} by two examples over rational function fields; in both, the only non-zero local
term is at the place $\infty$ and is a residue of a rational differential. In \S\ref{subsec:exrevisited} we show, using
the Birch and Swinnerton-Dyer formula, that the subgroups of $\Sha(E/K)$ exhibited by the pairing are all of
$\Sha(E/K)$ in the second example and in the first example for $k=2$.

\begin{theorem}\label{thm:example}
Let $k\geq1$, $q=2^{4k}$, $\zeta_3\in\F_4$ with $\zeta_3^2+\zeta_3+1=0$, $\Lambda=S^4+\zeta_3S^2+\zeta_3^2S$,
$Z=\{a\in\F_q:\Lambda(a)=0\}$, $K=\F_q(S)$, and $E\colon y^2+xy=x^3+\Lambda^4$. Write $\Tr=\Tr_{\F_q/\F_2}$ and
$Z^\perp=\{\lambda\in\F_q:\Tr(\lambda a)=0\ \forall a\in Z\}$, and for $\lambda\in\F_q$ let $\chi_\lambda$ be the
Artin--Schreier class of $\lambda S$.
\begin{enumerate}[label=\textup{(\alph*)},leftmargin=*]
\item $Z\subseteq Z^\perp$, and for $\lambda\in Z^\perp$ one has $\chi_\lambda\in\Sel^\varphi(E/K)$.
\item For $\lambda,\mu\in Z^\perp$: $\langle\iota(\chi_\mu),\iota(\chi_\lambda)\rangle=\Tr(\lambda\mu)$.
\item For $\lambda\in Z^\perp$: $\iota(\chi_\lambda)=0$ if and only if $\lambda\in Z$.
\item $\lambda\mapsto\iota(\chi_\lambda)$ induces an injection $Z^\perp/Z\hookrightarrow\Sha(E/K)/\hat\varphi\,\Sha(E'/K)$,
hence also an injection $Z^\perp/Z\hookrightarrow\Sha(E/K)/2\,\Sha(E/K)$.
\end{enumerate}
Consequently $\lambda\mapsto\iota(\chi_\lambda)$ induces an injection $Z^\perp/Z\hookrightarrow\Sha(E/K)[\varphi]$; its
image has order $2^{4k-4}$ and carries the non-degenerate alternating form induced by $\Tr(\lambda\mu)$. For
$k=1$ one has $Z^\perp=Z$, so the statement has content only for $k\geq2$.
\end{theorem}

In this example $b=\Lambda^2$ is a square and the group scheme $E'[2]$ is $\Z/2\times\mu_2$. In the second example $b$
is not a square, $E'[2]$ is a non-split extension of $\Z/2$ by $\mu_2$, and the norm elements $\beta_0$ are not
rational.

\begin{theorem}\label{thm:nonsplit}
Let $K=\F_8(S)$, $\varrho\in\F_8$ with $\varrho^3+\varrho+1=0$, $b=S^4+S^3+\varrho$ and $E\colon y^2+xy=x^3+b^2$. For
$c\in K$ let $\chi_c$ be the Artin--Schreier class of $c$.
\begin{enumerate}[label=\textup{(\alph*)},leftmargin=*]
\item $\chi_1,\chi_{\varrho S},\chi_{\varrho^2S}\in\Sel^\varphi(E/K)$.
\item On the basis $\chi_1,\chi_{\varrho S},\chi_{\varrho^2S}$, the matrix of
$\langle\iota(\cdot),\iota(\cdot)\rangle$ \textup{(}values in $\frac12\Z/\Z\simeq\F_2$\textup{)} is
$\left(\begin{smallmatrix}0&0&1\\0&0&1\\1&1&0\end{smallmatrix}\right)$.
\item The classes $\iota(\chi_1)$ and $\iota(\chi_{\varrho^2S})$ span a subgroup of $\Sha(E/K)[\varphi]$ of order
$4$ on which the pairing is non-degenerate, and this subgroup injects into $\Sha(E/K)/\hat\varphi\,\Sha(E'/K)$.
\end{enumerate}
\end{theorem}

\subsection{What is new, and relation to earlier work}
Much of \S\S\ref{sec:radical}--\ref{sec:parity} is due to Kramer \cite{Kr}, as indicated at each statement. New, to our
knowledge, are Theorem~\ref{thm:BM}; the identification of Kramer's pairing with the Cassels--Tate pairing
(Theorem~\ref{thm:W}, Corollary~\ref{cor:kramer}), hence Theorem~\ref{thm:radical}(b) as a statement about $\Sha$; the comparison with flat duality on
$\Sha(E/K)[\varphi]$ (Proposition~\ref{prop:milne}); the
supersingular case of Lemma~\ref{lem:good}, Lemma~\ref{lem:Dsize} in N\'eron-model form, and the root-number identities
of Theorems~\ref{thm:semistable}, \ref{thm:local} and~\ref{thm:finSha}; the growth of $\Sha[\varphi]$ in constant field extensions
(Corollary~\ref{cor:tower}); the determination of $\Sha$ in the examples; and the extension of Theorem~\ref{thm:BM} to
points of prime order $\ell=p$ (\S\ref{sec:ell}; for $\ell\neq p$ see Remark~\ref{rem:ellcases}(3)). The parity
results are the analogue for $p=2$ of \cite{TW}; for the $2^\infty$-Selmer corank see \cite{TY}, and for Brauer
groups of surfaces \cite{Fe,CF}.

\subsection{Organisation}
Theorem~\ref{thm:BM}, Corollary~\ref{cor:bilin} and the examples of \S\ref{sec:example} use only
\S\S\ref{sec:prelim}--\ref{sec:example}. Theorem~\ref{thm:W} uses \S\ref{sec:W} and Appendix~\ref{app:flat};
Theorem~\ref{thm:radical} uses in addition \S\ref{sec:radical} and, for (b) and (c), the alternating property from
Theorem~\ref{thm:BM}; the parity results of \S\ref{sec:parity} use the pairing only through Theorem~\ref{thm:radical}(c), and
\S\ref{subsec:exrevisited} also uses Corollary~\ref{cor:bilin}, Theorem~\ref{thm:radical}(b) and the matrices of
\S\ref{sec:example}. Section~\ref{sec:ell} uses only
\S\ref{sec:pairing}. Appendix~\ref{app:milne} uses \S\ref{sec:W} and the flat duality of \cite{Mi,DH}; nothing else
depends on it.

\section{Notation and conventions}\label{sec:notation}
A global field of characteristic $2$ is a finite extension of $\F_2(t)$; all its places $v$ are non-archimedean, $K_v$
is the completion, $\mathcal O_v$, $k_v$ (finite) the valuation ring and residue field, and $v$ also denotes the
normalised valuation. We fix a separable closure $\Ks$ and, for each $v$, a separable closure $\Kvs\supseteq\Ks$; this
identifies $G_v=\Gal(\Kvs/K_v)$ with a decomposition group in $G_K=\Gal(\Ks/K)$. For a field $F$, $H^i(F,-)$ is
continuous Galois cohomology of $G_F$ and $H^i_\fl(F,-)$ is fppf cohomology. We identify
$H^1(F,\Z/2)=\Hom_{\mathrm{cont}}(G_F,\Z/2)$ and view characters as functions $G_F\to\{0,1\}\subset\Z$; the sum of two
characters is taken modulo $2$. Every class that we pair in characteristic $2$ is $2$-torsion and every such pairing
takes values in $\frac12\Z/\Z\simeq\F_2$; consequently the signs of connecting homomorphisms and cup products play no
role there and are suppressed. Groups of functions and constants are written multiplicatively. Cochains are
inhomogeneous, with $(d\kappa)(\sigma,\tau)=\sigma\kappa(\tau)\cdot\kappa(\sigma)\cdot\kappa(\sigma\tau)^{-1}$ in
degree $1$ and similarly in other degrees. Dimensions are over $\F_2$.

Let $\F_{2^n}$ be the constant field of $K$ and $C$ the smooth projective curve over $\F_{2^n}$ with function field
$K$. For a place $v$ we put $q_v=|k_v|$ and $f_v=[k_v:\F_2]=n\deg(v)$, and $\pi_v$ is a uniformiser. We write
$H^1_\nr(K_v,\Z/2)$ for the subgroup of unramified classes (of order $2$) and $U_v\subseteq K_v^\times/K_v^{\times2}$ for
the image of $\mathcal O_v^\times$ (of index $2$). For an elliptic curve $A/K$, $c_v(A)$ is the Tamagawa number at $v$;
$\omega_A$ is the Hodge bundle, the pull-back along the zero section of the sheaf of relative differentials of the
N\'eron model of $A$ over $C$, and $\deg\omega_A$ is its degree; $w_v(A)\in\{\pm1\}$ is the local root number, and
$w(A/K)=\prod_vw_v(A)$. For an isogeny $\psi\colon A\to B$ with Kummer maps
$\delta_\psi\colon B(F)\to H^1_\fl(F,\ker\psi)$, the Selmer group $\Sel^\psi(A/K)\subseteq H^1_\fl(K,\ker\psi)$ consists
of the classes whose restriction to $K_v$ lies in $\delta_\psi(B(K_v))$ for every $v$, and $\Sel^2=\Sel^{[2]}$.

\section{Preliminaries}\label{sec:prelim}

\subsection{The curves and the isogenies}

\begin{remark}\label{rem:model}
Let $F$ be a field of characteristic $2$ and $E_0/F$ an elliptic curve with $j(E_0)\neq0$ (i.e.\ ordinary). Then $E_0$
has a model $y^2+xy=x^3+a_2x^2+a_6$ with $a_6\in F^\times$ \cite[Appendix~A, Prop.~A.1.1]{Si}. As $-(x,y)=(x,x+y)$, the
only point of order $2$ over $\bar F$ is $(0,\sqrt{a_6})$; it is $F$-rational if and only if $a_6=b^2$ with $b\in F$.
Thus the curves $E$ of \eqref{eq:E} are exactly the ordinary elliptic curves with an $F$-rational point of order $2$.
The discriminant of $E$ is $b^2$ and that of $E'$ is $b$.
\end{remark}

From now on $F$ denotes a field of characteristic $2$ containing $a_2$ and $b$ (in practice $F=K$, $K_v$, or a
separable closure). We write $x(P)$ for the $x$-coordinate of $P\neq O$.

\begin{lemma}\label{lem:elem}
On $E$ and over any $F$:
\begin{enumerate}[label=\textup{(E\arabic*)},leftmargin=*]
\item $T=(0,b)$ is the only point with $x=0$; $2T=O$; $\dv(x)=2(T)-2(O)$.
\item $x(P+T)=b/x(P)$ for $P\notin\{O,T\}$.
\item If $P_1,P_2,P_3\neq O$ are the three intersection points \textup{(}with multiplicity\textup{)} of $E$ with a
line $y=mx+\nu$, then $x(P_1)x(P_2)x(P_3)=(\nu+b)^2$.
\item \textup{(}doubling\textup{)} For $P\notin\{O,T\}$: $2P=T$ if and only if $x(P)^2=b$, and otherwise
$x(2P)=\bigl((x(P)^2+b)/x(P)\bigr)^2$. If $P=(x_1,y_1)$ and $2P\ne O$, then
$2P=\bigl(x_3,(m+1)x_3+x_1^2\bigr)$ with $m=(y_1+x_1^2)/x_1$, $x_3=m^2+m+a_2$.
\item \textup{(}addition\textup{)} If $P_1\neq\pm P_2$ and $P_1,P_2,P_1+P_2\notin\{O,T\}$, then
$x(P_1+P_2)=x(P_1)x(P_2)\,w^2$ with $w=(\nu+b)/(x(P_1)x(P_2))\in F(P_1,P_2)^\times$, where $y=mx+\nu$ is the line
through $P_1,P_2$.
\item The map $\kappa_F\colon E(F)\to F^\times/F^{\times2}$ given by $\kappa_F(P)=[x(P)]$ for $P\neq O,T$,
$\kappa_F(T)=[b]$ and $\kappa_F(O)=1$ is a homomorphism.
\end{enumerate}
\end{lemma}

\begin{proof}
Substituting $y=mx+\nu$ into the equation of $E$ gives
\begin{equation}\label{eq:cubic}
x^3+(a_2+m^2+m)x^2+\nu x+(\nu+b)^2=0,
\end{equation}
whose roots are the $x(P_i)$; this is (E3) (Vieta). (E1): $x=0$ forces $y^2=b^2$, i.e.\ $y=b$; since
$-T=(0,0+b)=T$, $2T=O$; $x$ has a double pole at $O$ and its only zero is $T$, hence $\dv(x)=2(T)-2(O)$. (E2): the line
through $T$ and $P$ (non-vertical, as $x(P)\neq0$) has $\nu=b$, so \eqref{eq:cubic} is $x\,(x^2+(a_2+m^2+m)x+b)$, and the
two non-zero roots $x(P)$, $x(-(P+T))=x(P+T)$ have product $b$. (E4): the tangent line at $P=(x_1,y_1)$ has slope
$m=(y_1+x_1^2)/x_1$ (differentiate: $x\,dy=(x^2+y)\,dx$ on $E$) and intercept $\nu=y_1+mx_1=x_1^2$. The third
intersection point is $-2P$, so by (E3) $x_1^2\,x(2P)=(x_1^2+b)^2$ whenever $2P\neq O$; hence $x(2P)=0$, i.e.\ $2P=T$, if
and only if $x_1^2=b$, and otherwise $x(2P)=((x_1^2+b)/x_1)^2$. The sum of the roots of \eqref{eq:cubic} gives
$x_3=x(2P)=m^2+m+a_2$, and $-2P=(x_3,mx_3+\nu)$, so $2P=(x_3,mx_3+\nu+x_3)$. (E5): apply (E3) to $P_1,P_2,P_3=-(P_1+P_2)$;
$\nu+b\neq0$ because $x(P_3)\neq0$. (E6): we must show $\kappa(P_1+P_2)=\kappa(P_1)\kappa(P_2)$. If $P_1$ or $P_2$ is
$O$ this is clear. If $P_2=T$ it follows from (E2) (and $\kappa(T)^2=1$). If $P_2=-P_1$ it holds because $x(-P)=x(P)$.
If $P_2=P_1$ it is (E4) (both sides are trivial). If $P_1+P_2=T$ with $P_1\notin\{O,T\}$ then $P_2=-P_1+T$,
and it follows from (E2) applied to $-P_1$, since $x(-P_1)=x(P_1)$. In the remaining case it is (E5).
\end{proof}

\begin{lemma}\label{lem:isog}
Define
\[
\varphi(x,y)=\Bigl(\frac{x^2+b}{x},\ \frac{(x^2+b)\,y+b\,(x^2+x+b)}{x^2}\Bigr),\qquad
\hat\varphi(x,y)=\bigl(x^2,\ y^2+a_2x^2\bigr).
\]
\begin{enumerate}[label=\textup{(\roman*)},leftmargin=*]
\item $\varphi$ is a separable isogeny $E\to E'$ of degree $2$ with kernel $\{O,T\}\simeq\Z/2$.
\item $\hat\varphi$ is a purely inseparable isogeny $E'\to E$ of degree $2$. Its kernel is
$\Spec\mathcal O_{E',O}/(z^2)=\Spec F[z]/(z^2)$, $z=x/y$, and $z\mapsto\zeta-1$ is an isomorphism of group schemes
$\ker\hat\varphi\simeq\mu_2=\Spec F[\zeta]/(\zeta^2-1)$.
\item $\hat\varphi\circ\varphi=[2]_E$ and $\varphi\circ\hat\varphi=[2]_{E'}$.
\end{enumerate}
\end{lemma}

\begin{proof}
Write $f_E(x,y)=y^2+xy+x^3+a_2x^2+b^2$ and $f_{E'}(x,y)=y^2+xy+x^3+a_2x^2+b$. A direct expansion in characteristic $2$
gives
\[
x^4\,f_{E'}(\varphi(x,y))=(x^2+b)^2f_E(x,y),\qquad f_E(\hat\varphi(x,y))=f_{E'}(x,y)^2 .
\]
Hence $\varphi$ and $\hat\varphi$ are morphisms of curves; they map $O$ to $O$ (and $\varphi(T)=O$), so they are
homomorphisms \cite[III.4.8]{Si}. These are the isogenies of \cite[\S1]{Kr}, where $\hat\varphi$ is the Frobenius
$\pi$ and $\varphi$ its dual $\psi$, in the model \eqref{eq:E}.

(i) $[F(x):F(x(\varphi))]=2$ since $x$ is a root of $X^2+x(\varphi)X+b$, which is separable; as $F(E)/F(x)$ and
$\varphi^*F(E')/F(x(\varphi))$ both have degree $2$, $\deg\varphi=2$ and $\varphi$ is separable. So $\ker\varphi$ is
\'etale of order $2$; it contains $T$.

(ii) $\hat\varphi^*F(E)=F(x^2,y^2)$ has index $2$ in $F(x,y)$ and $F(x,y)/F(x^2,y^2)$ is purely inseparable. With
$z_E=x/y$ on $E$, $\hat\varphi^*z_E=x^2/(y^2+a_2x^2)=z^2/(1+a_2z^2)$, so the scheme-theoretic fibre
$\hat\varphi^{-1}(O)$ is $\Spec\mathcal O_{E',O}/(z^2)=\Spec F[z]/(z^2)$. The group law of $E'$ in the parameter $z$
is given by the formal group law $z_1+z_2-a_1z_1z_2-a_2(z_1^2z_2+z_1z_2^2)+(\text{terms of degree}\geq4)$
\cite[IV.1]{Si}, with $a_1=1$; every monomial of degree $\geq3$ lies in $(z_1^2,z_2^2)$, so on $\ker\hat\varphi$ the
comultiplication is $z\mapsto z\otimes1+1\otimes z+z\otimes z$. With $\zeta=1+z$ this is $\zeta\mapsto\zeta\otimes\zeta$,
i.e.\ $\ker\hat\varphi\simeq\mu_2$.

(iii) Let $P=(x_1,y_1)$ with $P,2P\neq O$. Then $x(\hat\varphi\varphi P)=((x_1^2+b)/x_1)^2=x(2P)$ by (E4), and
comparing $y(\hat\varphi\varphi P)=Y^2+a_2X^2$, $(X,Y)=\varphi(P)$, with the formula for $y(2P)$ in (E4), using the
equation of $E$, gives $y(\hat\varphi\varphi P)=y(2P)$ (cf.\ \cite[(4)]{Kr}). Hence $\hat\varphi\circ\varphi=[2]_E$ on a
dense open set, hence everywhere. Then $\varphi\circ\hat\varphi\circ\varphi=\varphi\circ[2]_E=[2]_{E'}\circ\varphi$, and
$\varphi$ is surjective, so $\varphi\circ\hat\varphi=[2]_{E'}$.
\end{proof}

In particular $\hat\varphi$ is the dual isogeny of $\varphi$, and $\ker\varphi\simeq\Z/2$, $\ker\hat\varphi\simeq\mu_2$ are
Cartier dual. The unique perfect pairing $\Z/2\times\mu_2\to\mu_2$ is $(n,\zeta)\mapsto\zeta^n$; it is the Weil pairing
$e_\varphi$.

\begin{remark}[Frobenius and Verschiebung]\label{rem:FV}
Let $E'^{(2)}\colon y^2+xy=x^3+a_2^2x^2+b^2$ be the Frobenius twist of $E'$ and
$\mathrm{Fr}\colon E'\to E'^{(2)}$, $(x,y)\mapsto(x^2,y^2)$. The map $(x,y)\mapsto(x,y+a_2x)$ is an isomorphism
$E'^{(2)}\to E$ (because $a_2^2+a_2=\wp(a_2)$), and $\hat\varphi$ is its composite with $\mathrm{Fr}$. Thus, up to
this isomorphism, $\hat\varphi$ is the Frobenius isogeny of $E'$ and $\varphi$ is the Verschiebung, and
$\Sel^\varphi(E/K)$ is the Selmer group of the Verschiebung descent of \cite{Kr,Vo,Ul,Moe}. We keep the model
\eqref{eq:E} because it makes the point $T$, the function $x$ with $\dv(x)=2(T)-2(O)$ and the relation
$x\circ\tau_T=b/x$ explicit.
\end{remark}

\begin{lemma}[Kummer map of $\varphi$]\label{lem:deltaphi}
Let $\delta_\varphi\colon E'(F)\to H^1(F,\ker\varphi)=H^1(F,\Z/2)=F/\wp F$ be the connecting map of
$0\to\{O,T\}\to E(\Fs)\xrightarrow{\varphi}E'(\Fs)\to0$. Then $\delta_\varphi(O)=0$ and
$\delta_\varphi(R)=[x(R)+a_2]$ for $R\neq O$; if $x(R)\neq0$ this equals $[b/x(R)^2]$.
\end{lemma}

\begin{proof}
For $R\in E'(F)$ write $X=x(R)$, $Y=y(R)$. Dividing the equation of $E'$ by $X^2$ shows $(Y/X)^2+Y/X=X+a_2+b/X^2$ when
$X\neq0$, so $[X+a_2]=[b/X^2]$. Let $X\neq0$ and $Q\in E(\Fs)$ with $\varphi(Q)=R$; then $Q\notin\{O,T\}$ and
$x(Q)$ is a root of $u^2+Xu+b$, the other root being $b/x(Q)=x(Q+T)$ by (E2). A point of $\varphi^{-1}(R)=\{Q,Q+T\}$ is
determined by its $x$-coordinate: if $x(Q+T)=x(Q)$ then $x(Q)^2=b$, contradicting $X=x(Q)+b/x(Q)\neq0$. Hence
$\sigma Q=Q+T$ if and only if $\sigma$ exchanges the roots of $u^2+Xu+b$, i.e.\ of $z^2+z+b/X^2$ ($z=u/X$), so
$\delta_\varphi(R)=[b/X^2]$. If $X=0$ then $R=(0,\sqrt b)$ with $\sqrt b\in F$, and $\varphi^{-1}(R)$ consists of the
two points $(\sqrt b,y)$ with $y^2+\sqrt b\,y=b\sqrt b+a_2b+b^2$, i.e.\ $(y/\sqrt b)^2+y/\sqrt b=\sqrt b+a_2+b$; so
$\delta_\varphi(R)=[\sqrt b+b+a_2]=[a_2]$, since $\sqrt b+b=\wp(\sqrt b)$.
\end{proof}

\subsection{Cyclic algebras and local symbols}\label{subsec:symbols}
Let $F$ be a field of characteristic $2$. The Artin--Schreier map identifies $F/\wp F$ with $H^1(F,\Z/2)$: $c$
corresponds to $\chi_c(\sigma)=\sigma\alpha-\alpha$, $\wp(\alpha)=c$. For $\chi\in H^1(F,\Z/2)$ and $u\in F^\times$, the
cyclic algebra $(\chi,u)\in\Br(F)[2]$ is $\delta\chi\cup u$, where $\delta\colon H^1(F,\Q/\Z)\to H^2(F,\Z)$
\cite[XIV, \S1]{Se}, \cite[\S4.7]{GS}. Lifting $\chi$ to $\{0,\frac12\}\subset\Q$ one gets
$\delta\chi(\sigma,\tau)=\chi(\sigma)\chi(\tau)$ (with $\chi$ now $\{0,1\}$-valued), so $(\chi,u)$ is represented by
the $2$-cocycle
\begin{equation}\label{eq:cycliccocycle}
(\chi,u)(\sigma,\tau)=u^{\chi(\sigma)\chi(\tau)} .
\end{equation}
It is bi-additive in $(\chi,u)$ and compatible with restriction to extensions of $F$. We use the following standard facts.

\begin{fact}[norm criterion]\label{fact:norm}
If $\chi\neq0$ cuts out $L/F$, then $(\chi,u)=0$ if and only if $u\in\Nm_{L/F}(L^\times)$
\cite[XIV, \S1, Prop.~2, Cor.~1]{Se}, \cite[Cor.~4.7.5]{GS}.
\end{fact}

\begin{fact}[Albert--Brauer--Hasse--Noether]\label{fact:ABHN}
For a global field $K$ the sequence
\[
0\to\Br K\to\bigoplus_v\Br K_v\xrightarrow{\sum\inv_v}\Q/\Z\to0
\]
is exact
\cite[Thm.~8.1.17]{NSW}; for function fields see also \cite[I, Thm.~A.7]{Mi} \textup{(}the group $H^1(G_k,\mathrm{Jac})$
occurring there vanishes for finite $k$ by Lang's theorem\textup{)}. In particular, for $\chi\in H^1(K,\Z/2)$ and
$u\in K^\times$: $\sum_v\inv_v(\chi_v,u)=0$; and $(\chi,u)=0$ as soon as $(\chi_v,u)=0$ for all $v$.
\end{fact}

\begin{fact}[Schmid's formula]\label{fact:schmid}
Let $F=k((\pi))$ with $k$ finite of characteristic $2$. For $c\in F$ and $u\in F^\times$,
$\inv(\chi_c,u)=\Tr_{k/\F_2}\Res(c\,du/u)\in\F_2=\frac12\Z/\Z$ \cite{Sch}, \cite[XIV, \S5, Prop.~15 and its
Corollary]{Se} \textup{(}there $[c,u):=(\chi_c,u)$\textup{)}. In particular: $\chi_c=0$ if $c\in\pi k[[\pi]]$;
and if $c\in k[[\pi]]$ and $v(u)$ is even, then $(\chi_c,u)=0$, since $\Res(c\,du/u)=v(u)\,c(0)$.
\end{fact}

\section{The homogeneous-space pairing}\label{sec:pairing}

\subsection{The definition}\label{subsec:PS}
Let $A$ be an abelian variety over a global field $K$ (of any characteristic). For a torsor $X$ under $A$ write
$X_\s=X\times_K\Ks$, $\Ks(X)$ for its function field, and $\Div(X_\s)\supseteq\Div^0(X_\s)\supseteq\Prin(X_\s)$.
Translations act trivially on $\Pic^0$, so $\Pic^0(X_\s)=A^\vee(\Ks)$ canonically and $G_K$-equivariantly.

\begin{lemma}\label{lem:H3}
$H^3(G_K,\Ks{}^\times)=0$.
\end{lemma}

\begin{proof}
This is recalled in \cite[\S3.1]{PS}; it amounts to $\operatorname{scd}(G_K)\leq2$. If $p=\operatorname{char}K>0$,
the part prime to $p$ vanishes by class field theory \cite[I, Cor.~4.21]{Mi} (Chapter~I of \cite{Mi} treats global
fields of both characteristics), and the $p$-primary part vanishes because
$\operatorname{scd}_p(G_K)\leq\operatorname{cd}_p(G_K)+1\leq2$, as $\operatorname{cd}_p(G_K)\leq1$ in characteristic $p$
\cite[I, \S3.2, II, \S2.2]{SeGC}.
\end{proof}

\begin{definition}[{\cite[\S3.1]{PS}; cf.\ \cite[I, Rem.~6.11]{Mi}}]\label{def:PS}
Let $a\in\Sha(A/K)$, $a'\in\Sha(A^\vee/K)$, and let $X$ be a torsor representing $a$. Represent $a'$ by a cocycle with
values in $\Pic^0(X_\s)$, lift it to a continuous cochain $a''\in C^1(G_K,\Div^0(X_\s))$, and put
$\eta=da''\in Z^2(G_K,\Prin(X_\s))$. Choose $f'\in Z^2(G_K,\Ks(X)^\times)$ with $\dv f'=\eta$ (it exists by
Lemma~\ref{lem:H3}). For each place $v$, let $o_v\in\Br K_v=H^2(G_v,\Kvs{}^\times)$ be the class whose
image in $H^2(G_v,\Kvs(X)^\times)$, where $\Kvs(X)$ is the function field of $X\times_K\Kvs$, is the class of the
restriction $f'_v$. Then
$\langle a,a'\rangle:=\sum_v\inv_v(o_v)\in\Q/\Z$.
\end{definition}
In words: the exact sequences of $G_K$-modules
\[
0\to\Prin(X_\s)\to\Div^0(X_\s)\to\Pic^0(X_\s)\to0,\qquad
1\to\Ks{}^\times\to\Ks(X)^\times\xrightarrow{\ \dv\ }\Prin(X_\s)\to0
\]
give connecting maps $H^1(G_K,\Pic^0(X_\s))\to H^2(G_K,\Prin(X_\s))\to H^3(G_K,\Ks{}^\times)=0$, so the class $a'$ can
be carried to a $2$-cocycle $f'$ of functions; $f'$ is determined up to a coboundary and a cocycle of constants, i.e.\
up to an element of $\Br K$. Over $K_v$ the class $a'_v$ vanishes and $X$ has points, so $f'_v$ is cohomologous to a
cocycle of constants, whose class $o_v$ is the local obstruction. The pairing is the sum of these local obstructions;
an ambiguity by an element of $\Br K$ does not change the sum (Fact~\ref{fact:ABHN}). For example, if $a'=0$ one may
take $a''=0$ and $f'=1$, and all $o_v$ vanish.

\emph{Finiteness.} Two admissible lifts $f'$ differ by a coboundary times a $2$-cocycle of constants, i.e.\ by a
global Brauer class, whose localisations vanish for almost all $v$; so it suffices that $o_v=0$ for almost all $v$ for
one choice of data. For the lifts used in this article this is proved directly: for the lift of \S\ref{subsec:proofBM}
the local terms vanish outside the finite set of places where $a_2\notin\mathcal O_v$, $b\notin\mathcal O_v^\times$,
$\chi'$ is ramified, or $t_0\notin\mathcal O_v$ (Lemma~\ref{lem:BMgood}), and for the
lift of \S\ref{subsec:thmA} see \S\ref{subsec:A6}; this covers all pairs of classes for which the pairing is computed
in this article. For arbitrary pairs, as in Proposition~\ref{prop:F} and Corollary~\ref{cor:bilinear}, this is the standard argument
of \cite[\S3.1]{PS}, \cite[I, Rem.~6.11]{Mi}, which does not use the characteristic; we sketch it for function fields
(for number fields replace $C$ by $\operatorname{Spec}\mathcal O_K$ and exclude the archimedean places). The cochains
$a''$ and $f'$ factor through $\Gal(L/K)$ for a finite Galois extension $L/K$. Choose a dense open $U$ of the curve $C$
with function field $K$ over which $A$ extends to an abelian scheme, $X$ to a smooth proper torsor $\mathcal X$, the
normalisation $U_L$ is finite \'etale over $U$, and the functions occurring in $f'$ are units on $\mathcal X_{U_L}$
outside the closure $\mathcal Z$ of the supports of their divisors. By the Lang--Weil estimates \cite{LW} the image of
$\mathcal Z$ in $\mathcal X$ has $O(q_v^{d-1})$ points over $k_v$ ($d=\dim A$), while $\mathcal X_{k_v}$ is trivial by
Lang's theorem \cite{Lan} and has at least $(\sqrt{q_v}-1)^{2d}$ points. So for almost all $v\in U$ there is
$P_v\in X(K_v)$ whose reduction avoids $\mathcal Z$ (Hensel's lemma); by Lemma~\ref{lem:ev}(b), $o_v$ is then the class of
a cocycle with values in $\mathcal O_{L_w}^\times$ for the unramified extension $L_w/K_v$, and
$H^2(\Gal(L_w/K_v),\mathcal O_{L_w}^\times)=0$ \cite[XII, \S1, Prop.~1]{Se} gives $o_v=0$.

The class $o_v$ exists because $a'_v=0$: then $a''_v=d D+\dv(\varpi)$ for some $D\in\Div^0(X_{\Kvs})$ and some
$1$-cochain $\varpi$ with values in $\Kvs(X)^\times$, so $f'_v/d\varpi$ has trivial divisor, i.e.\ is a constant
$2$-cocycle. It is unique by Lemma~\ref{lem:ev}(a), since $X(K_v)\neq\emptyset$ ($a_v=0$).

\begin{lemma}\label{lem:ev}
Let $F$ be a field and $X$ a smooth geometrically integral $F$-variety.
\begin{enumerate}[label=\textup{(\alph*)},leftmargin=*]
\item If $X(F)$ is Zariski dense in $X$, then $\Br F=H^2(G_F,\Fs{}^\times)\to H^2(G_F,\Fs(X)^\times)$ is injective.
\item If $P\in X(F)$, $f\in Z^2(G_F,\Fs(X)^\times)$ has all its values regular and invertible at $P$,
and $c\in H^2(G_F,\Fs{}^\times)$ has the same image as $f$ in $H^2(G_F,\Fs(X)^\times)$, then $c=[\ev_P f]$.
\end{enumerate}
If $F$ is a local field and $X$ is a torsor under an abelian variety with $X(F)\neq\emptyset$, the hypothesis of (a)
holds.
\end{lemma}

\begin{proof}
(a) Let $c\in Z^2(G_F,\Fs{}^\times)$ with $c=d\kappa$, $\kappa\in C^1(G_F,\Fs(X)^\times)$ continuous. A continuous cochain
from a profinite group to a discrete module takes finitely many values, so there is $P\in X(F)$ outside the zeros and
poles of all values of $\kappa$. Evaluation at $P$ is $G_F$-equivariant on functions regular and invertible at $P$
(because $P$ is $F$-rational) and is the identity on constants, so $c=d(\ev_P\kappa)$ is a coboundary.
(b) Let $\mathcal O_P$ be the local ring of $X_\s$ at $P$. It is regular, hence factorial \cite{AB}, so
$M=\Fs(X)^\times/\mathcal O_P^\times$ is the free abelian group on the set of height-one primes of $\mathcal O_P$, on
which $G_F$ acts with open stabilisers; by Shapiro's lemma $H^1(G_F,M)$ is a direct sum of groups
$\Hom_{\mathrm{cont}}(H,\Z)=0$ for open subgroups $H$. Hence $H^2(G_F,\mathcal O_P^\times)\to H^2(G_F,\Fs(X)^\times)$ is
injective, and $f$, $c$ have the same class in $H^2(G_F,\mathcal O_P^\times)$. Apply the
$G_F$-equivariant retraction $\ev_P\colon\mathcal O_P^\times\to\Fs{}^\times$. Finally, let $F$ be a local field and
$X\simeq A$ with $A$ an abelian variety of dimension $g$, and suppose $A(F)\subseteq Z$ for a closed subset
$Z\subsetneq A$. Choose $P\in A(F)$, an open neighbourhood $U$ of $P$ and an \'etale morphism
$f\colon U\to\mathbb A^g_F$. By the implicit function theorem (Hensel's lemma), $f(U(F))$ contains a non-empty open
ball $D\subseteq F^g$. Then $D\subseteq f(Z\cap U)$, whose Zariski closure has dimension $<g$ because $f$ is
quasi-finite. But a polynomial vanishing on a non-empty open ball of $F^g$ is zero ($F$ is infinite), so $D$ is
Zariski dense in $\mathbb A^g$, a contradiction. Hence $A(F)$ is Zariski dense in $A$.
\end{proof}

\begin{remark}[independence of choices]\label{rem:PSuse}
The value $\langle a,a'\rangle$ does not depend on the choices: two lifts $f'$ of the class of $\eta$ differ by a
coboundary times a cocycle with values in $\Ks{}^\times$, which changes $\sum_v\inv_v o_v$ by the sum of the invariants
of a global Brauer class, i.e.\ by $0$ (Fact~\ref{fact:ABHN}); changing $a''$ changes $\eta$ by the coboundary of a
$\Prin$-valued cochain, hence $f'$ by a coboundary; changing the cocycles representing $a,a'$, or $X$ within its
isomorphism class, transports all data. These verifications use only Fact~\ref{fact:ABHN} and are valid in every
characteristic. Poonen and Stoll state the definition for an arbitrary global field and refer to Milne for it
\cite[\S3, \S3.1]{PS}; for what is known about its formal properties see \S\ref{subsec:introradical}.
\end{remark}

\subsection{Functoriality and bilinearity}
We now prove that the pairing is functorial and bilinear. For the part prime to the characteristic, the functoriality is stated in
\cite[I, Rem.~6.10(a)]{Mi}. The proofs below use no
Weil pairing and no cohomology of finite group schemes, only Definition~\ref{def:PS}, so they apply to the $p$-primary
parts in characteristic $p$.

\begin{proposition}\label{prop:F}
Let $\rho\colon A\to B$ be a surjective homomorphism of abelian varieties over a global field $K$ of arbitrary
characteristic. For $x\in\Sha(A/K)$ and $y\in\Sha(B^\vee/K)$,
\[
\langle x,\rho^\vee_*y\rangle_A=\langle\rho_*x,y\rangle_B .
\]
\end{proposition}

\begin{proof}
Let $X$ be a torsor representing $x$ and $Y=X\times^{A,\rho}B$ the induced torsor, which represents $\rho_*x$. The
natural map $\pi\colon X\to Y$ is a $K$-morphism which over $\Ks$ is $\rho$ up to translations; in particular it is
flat and surjective. As $X$ and $Y$ are smooth, Weil and Cartier divisors agree, so pull-back of divisors is defined;
it maps principal divisors to principal divisors and preserves algebraic equivalence to zero. Thus pull-back by $\pi$
maps $\Ks(Y)^\times\to\Ks(X)^\times$, $\Div^0(Y_\s)\to\Div^0(X_\s)$ and
$\Prin(Y_\s)\to\Prin(X_\s)$ compatibly with $\dv$ and with the $G_K$-actions, and on
$\Pic^0(Y_\s)=B^\vee(\Ks)\to\Pic^0(X_\s)=A^\vee(\Ks)$ it is $\rho^\vee$ (translations act trivially on $\Pic^0$). Hence
if $(a'',\eta,f')$ are data as in Definition~\ref{def:PS} for $(\rho_*x,y)$ on $Y$, then $(\pi^*a'',\pi^*\eta,\pi^*f')$
are such data for $(x,\rho^\vee_*y)$ on $X$. If $o^Y_v\in\Br K_v$ has image $[f'_v]$ in $H^2(G_v,\Kvs(Y)^\times)$, then,
$\pi^*$ being the identity on constants, $o^Y_v$ has image $[\pi^*f'_v]$ in $H^2(G_v,\Kvs(X)^\times)$. Since
$X(K_v)\neq\emptyset$, the uniqueness statement of Lemma~\ref{lem:ev}(a) shows that $o^X_v=o^Y_v$ for all $v$.
\end{proof}

\begin{corollary}\label{cor:bilinear}
For every abelian variety $A$ over a global field $K$, the pairing of Definition~\ref{def:PS} is bilinear. In
particular $\langle x,0\rangle=\langle0,y\rangle=0$, and $\langle nx,y\rangle=\langle x,ny\rangle$ for $n\in\Z$.
\end{corollary}

\begin{proof}
\emph{Additivity in the second variable.} Let $X$ represent $a\in\Sha(A/K)$ and let $a'_1,a'_2\in\Sha(A^\vee/K)$. If
$(a''_i,f'_i)$ are data as in Definition~\ref{def:PS} for $(a,a'_i)$ on $X$, then $(a''_1+a''_2,f'_1f'_2)$ are data for
$(a,a'_1+a'_2)$, and the corresponding local classes are $o_{1,v}+o_{2,v}$. Hence
$\langle a,a'_1+a'_2\rangle=\langle a,a'_1\rangle+\langle a,a'_2\rangle$.

\emph{Additivity in the first variable.} Let $x_1,x_2\in\Sha(A/K)$ and $y\in\Sha(A^\vee/K)$. Let
$m,p_1,p_2\colon A\times A\to A$ be the sum and the two projections; they are surjective, $m=p_1+p_2$, and
$(A\times A)^\vee=A^\vee\times A^\vee$ with $p_1^\vee(y)=(y,0)$, $p_2^\vee(y)=(0,y)$, $m^\vee(y)=(y,y)=p_1^\vee(y)+p_2^\vee(y)$.
By Proposition~\ref{prop:F} for $m$, additivity in the second variable, and Proposition~\ref{prop:F} for $p_1$ and
$p_2$,
\[
\langle x_1+x_2,y\rangle_A=\langle(x_1,x_2),m^\vee_*y\rangle_{A\times A}
=\sum_{i=1,2}\langle(x_1,x_2),(p_i^\vee)_*y\rangle_{A\times A}=\langle x_1,y\rangle_A+\langle x_2,y\rangle_A .
\]
\end{proof}

For an elliptic curve $E$ we identify $E^\vee=E$ via $P\mapsto[(P)-(O)]$; under this identification the dual of an isogeny is the dual
isogeny, and the dual of $[n]_E$ is $[n]_E$.

\section{The explicit formula}\label{sec:main}

Throughout this section $K$ is a global field of characteristic $2$, $E,E',\varphi,\hat\varphi$ are as in
\eqref{eq:E} and Lemma~\ref{lem:isog}, and $\iota\colon H^1(K,\Z/2)\to H^1(K,E)$, $\iota(\chi)=[\sigma\mapsto\chi(\sigma)T]$.
By the Kummer sequence of $\varphi$ over $K_v$, $\iota(\chi)_v=0$ if and only if $\chi_v\in\delta_\varphi(E'(K_v))$, so
$\iota(\chi)\in\Sha(E/K)$ if and only if $\chi\in\Sel^\varphi(E/K)$.

\subsection{The norm element}\label{subsec:A1}

\begin{lemma}[{cf.\ \cite[Lemma~3.2]{Kr}}]\label{lem:norm}
Let $\chi'\in\Sel^\varphi(E/K)$. Then $(\chi',b)=0$ in $\Br K$, and there is $\beta_0\in\Ks{}^\times$ such that
$\sigma\beta_0=\beta_0$ if $\chi'(\sigma)=0$ and $\sigma\beta_0=\bar\beta_0:=b/\beta_0$ if $\chi'(\sigma)=1$.
Equivalently: $\beta_0\in K^\times$ if $\chi'=0$, and $\beta_0\in L'^\times$ with $\Nm_{L'/K}\beta_0=b$ if
$\chi'\neq0$, $L'=K(\chi')$.
\end{lemma}

\begin{proof}
Let $v$ be a place. Since $\iota(\chi')_v=0$ there is $P'_v\in E(\Kvs)$ with $\sigma P'_v-P'_v=\chi'(\sigma)T$ for
$\sigma\in G_v$. If $\chi'_v\neq0$, let $L'_v$ be the quadratic extension of $K_v$ cut out by $\chi'_v$; then
$P'_v\in E(L'_v)\smallsetminus E(K_v)$, so $P'_v\notin\{O,T\}$ and, for $\sigma$ with $\chi'(\sigma)=1$,
$\sigma x(P'_v)=x(P'_v+T)=b/x(P'_v)$ by (E2). Thus $\Nm_{L'_v/K_v}x(P'_v)=b$, and $(\chi'_v,b)=0$ by
Fact~\ref{fact:norm}; if $\chi'_v=0$ this is trivial. By Fact~\ref{fact:ABHN}, $(\chi',b)=0$ in $\Br K$, so by
Fact~\ref{fact:norm} $b=\Nm_{L'/K}\beta_0$ for some $\beta_0\in L'^\times$ (if $\chi'\neq0$; if $\chi'=0$ take any
$\beta_0\in K^\times$). Regarding $\beta_0\in\Ks$, $\sigma\beta_0$ is $\beta_0$ or its conjugate $b/\beta_0$
according as $\chi'(\sigma)=0$ or $1$.
\end{proof}

This proves Theorem~\ref{thm:BM}(a); note that $t_0=\beta_0+\bar\beta_0$ is fixed by $G_K$.

\begin{remark}[an explicit $\beta_0$]\label{rem:beta0}
Let $K=\F_q(S)$, $b\in\F_q[S]$ with $b=\mathrm{lc}(b)\prod_a(S+a)^{e_a}$, all $a\in\F_q$, and
$\chi'=\chi_{\lambda S+c_0}$ with $\lambda\in\F_q^\times$, $c_0\in\F_q$. Let $L'=K(y)$, $y^2+y=\lambda S+c_0$, and suppose
that for every root $a$ there is $y_a\in\F_q$ with $y_a^2+y_a=\lambda a+c_0$ (i.e.\ $\Tr_{\F_q/\F_2}(\lambda a+c_0)=0$, which
says that $\chi'$ is trivial at the place $S=a$). Since $\Nm(y+y_a)=(y+y_a)(y+y_a+1)=\wp(y)+\wp(y_a)=\lambda(S+a)$, the
element $\beta_0=d\prod_a(y+y_a)^{e_a}$ with $d^2=\mathrm{lc}(b)\lambda^{-\deg b}$ satisfies $\Nm\beta_0=b$.
\end{remark}

\subsection{The torsor}\label{subsec:torsor}
Let $\chi\in\Sel^\varphi(E/K)$. The torsor of $\iota(\chi)$ is $X=E$ with the twisted Galois action
$\sigma^*Q=\sigma Q+\chi(\sigma)T$ on points of $X(\Ks)=E(\Ks)$; this is an action because $T$ is rational and $2T=O$,
and for $Q_0\in X(\Ks)$ the cocycle $\sigma\mapsto\sigma^*Q_0-Q_0$ is cohomologous to $\chi(\sigma)T$. On functions the
twisted action is $\sigma^*f=(\sigma f)\circ\tau_{\chi(\sigma)T}$ (so that $(\sigma^*f)(\sigma^*Q)=\sigma(f(Q))$), and
on divisors $\sigma^*(Q)=(\sigma^*Q)$. For the coordinate function $x\in K(E)$, (E2) gives
\begin{equation}\label{eq:twisth}
\sigma^*x=b^{\chi(\sigma)}\,x^{(-1)^{\chi(\sigma)}},\qquad\text{i.e.}\quad \sigma^*x=b/x\ \text{ if }\chi(\sigma)=1,\quad
\sigma^*x=x\ \text{ if }\chi(\sigma)=0 .
\end{equation}
Thus $x$ is not a rational function on $X$ unless $\chi=0$, whereas $x+b/x$ is; it is the pull-back of the
coordinate $x_{E'}$ under the $K$-morphism $\pi\colon X\to E'$ which is $\varphi$ on points. A point $P\in E(\Kvs)$ lies
in $X(K_v)$ if and only if $\sigma P=P+\chi(\sigma)T$ for all $\sigma\in G_v$, i.e.\ if and only if
$\varphi(P)\in E'(K_v)$ and $\delta_\varphi(\varphi(P))=\chi_v$; hence
\begin{equation}\label{eq:piX}
\pi(X(K_v))=\{R\in E'(K_v):\ \delta_\varphi(R)=\chi_v\}.
\end{equation}
Since translations act trivially on $\Pic^0$, the map $\sum n_i(Q_i)\mapsto\sum n_iQ_i$ identifies $\Pic^0(X_\s)$ with
$E(\Ks)$ with its usual Galois action, and $\iota(\chi')$ is represented by the cocycle
$\sigma\mapsto\chi'(\sigma)T$.

\subsection{Proof of Theorem~\ref{thm:BM}(b)}\label{subsec:proofBM}
\begin{proof}
Let $X$ be the torsor of $\iota(\chi)$ of \S\ref{subsec:torsor}: $X_\s=E_\s$ with $\sigma^*Q=\sigma Q+\chi(\sigma)T$, and
$\sigma^*x=b/x$ or $x$ according as $\chi(\sigma)=1$ or $0$ \eqref{eq:twisth}. Put $\psi:=\chi+\chi'$, viewed as a function
$G_K\to\{0,1\}$ with $\psi(\sigma)\equiv\chi(\sigma)+\chi'(\sigma)\pmod2$.

\emph{Step 1: the function.} Put
\[
f:=x+\frac bx+t_0=\frac{(x+\beta_0)(x+\bar\beta_0)}{x}\in\Ks(X)^\times .
\]
For $\sigma\in G_K$, $\sigma^*$ fixes or exchanges $x$ and $b/x$, and $\sigma$ fixes $t_0$ (it fixes or exchanges
$\beta_0$ and $\bar\beta_0$). Hence $\sigma^*f=f$, i.e.\ $f\in K(X)^\times$. (In fact $f=\pi^*(x_{E'}+t_0)$ for the
$K$-morphism $\pi\colon X\to E'$ which is $\varphi$ on points.)

\emph{Step 2: its divisor.} Choose $Q_0\in E(\Ks)$ with $x(Q_0)=\beta_0$ (the equation $y^2+\beta_0y=\beta_0^3+a_2\beta_0^2+b^2$ is
separable in $y$ because $\beta_0\neq0$); then $Q_0\notin\{O,T\}$. The points with $x$-coordinate $\beta_0$ are $\pm Q_0$, and those with
$x$-coordinate $\bar\beta_0$ are $\pm Q_0+T$ by (E2). Put
\[
\mathcal Z=(Q_0)+(-Q_0),\qquad \bar{\mathcal Z}=(Q_0+T)+(-Q_0+T),\qquad \Delta=\mathcal Z-(O)-(T),\qquad
\bar\Delta=\bar{\mathcal Z}-(O)-(T).
\]
We subtract $(O)+(T)$ rather than $2(O)$ because the set $\{O,T\}$ is stable under the twisted action
($\sigma^*(O)=(\chi(\sigma)T)$), whereas $(O)$ alone is not. By (E1), $\dv(x+\beta_0)=\mathcal Z-2(O)$, $\dv(x+\bar\beta_0)=\bar{\mathcal Z}-2(O)$ and $\dv(x)=2(T)-2(O)$, so
\begin{equation}\label{eq:divf}
\dv f=\mathcal Z+\bar{\mathcal Z}-2(O)-2(T)=\Delta+\bar\Delta .
\end{equation}
Since $x(\sigma Q_0)=\sigma\beta_0$ is $\beta_0$ or $\bar\beta_0$ according as $\chi'(\sigma)=0$ or $1$, we have
$\sigma Q_0\in\pm Q_0+\chi'(\sigma)T$, hence $\sigma^*(\pm Q_0)\in\{Q_0+\psi(\sigma)T,\,-Q_0+\psi(\sigma)T\}$ (as
$2T=O$); also $\sigma^*$ preserves $\{O,T\}$. Therefore
\begin{equation}\label{eq:sigmaDelta}
\sigma^*\Delta=\Delta\ \text{ if }\psi(\sigma)=0,\qquad \sigma^*\Delta=\bar\Delta\ \text{ if }\psi(\sigma)=1 .
\end{equation}
The divisor $\Delta$ has degree $0$ and its class in $\Pic^0(X_\s)=E(\Ks)$ is $Q_0-Q_0-O-T=T$.

\emph{Step 3: the lift.} Define $a''\in C^1(G_K,\Div^0(X_\s))$ by
\[
a''_\sigma:=\psi(\sigma)\,\Delta+\bigl(\sigma^*(O)-(O)\bigr).
\]
Since $\sigma^*(O)=(\chi(\sigma)T)$, the second summand has class $\chi(\sigma)T$, so the class of $a''_\sigma$ in
$E(\Ks)$ is $\psi(\sigma)T+\chi(\sigma)T=\chi'(\sigma)T$; thus $a''$ lifts the cocycle
$\sigma\mapsto\chi'(\sigma)T$ representing $\iota(\chi')$. The second summand is the coboundary of the $0$-cochain
$(O)$, so it does not contribute to $da''$. For the first, using \eqref{eq:sigmaDelta} and
$\psi(\sigma\tau)=\psi(\sigma)+\psi(\tau)-2\psi(\sigma)\psi(\tau)$: if $\psi(\sigma)=0$ then
$\psi(\tau)\sigma^*\Delta-\psi(\sigma\tau)\Delta+\psi(\sigma)\Delta=\psi(\tau)\Delta-\psi(\tau)\Delta=0$, and if
$\psi(\sigma)=1$ it equals $\psi(\tau)\bar\Delta-(1-\psi(\tau))\Delta+\Delta=\psi(\tau)(\Delta+\bar\Delta)$. By
\eqref{eq:divf},
\[
(da'')(\sigma,\tau)=\psi(\sigma)\psi(\tau)\,\dv f=\dv\bigl(f^{\psi(\sigma)\psi(\tau)}\bigr).
\]
As $f$ is $G_K$-invariant and $(\sigma,\tau)\mapsto\psi(\sigma)\psi(\tau)$ is a $\Z$-valued $2$-cocycle,
$f'(\sigma,\tau):=f^{\psi(\sigma)\psi(\tau)}$ lies in $Z^2(G_K,\Ks(X)^\times)$. Thus $f'$ is a lift as in
Definition~\ref{def:PS} for the pair $(\iota(\chi),\iota(\chi'))$.

\emph{Step 4: evaluation.} Fix $v$ and let $P_v\in E(\Kvs)$ with $\varphi(P_v)=R_v$. As $\delta_\varphi(R_v)=\chi_v$,
$\sigma P_v=P_v+\chi(\sigma)T$ for $\sigma\in G_v$, i.e.\ $P_v\in X(K_v)$; conversely $\varphi$ maps $X(K_v)$ onto
$\{R\in E'(K_v):\delta_\varphi(R)=\chi_v\}$, and $X(K_v)$ is infinite, so $R_v$ can be chosen as required. The function
$f$ is regular and invertible at $P_v$ because $P_v\notin\{O,T\}$ ($R_v\neq O$) and
$f(P_v)=x(P_v)+b/x(P_v)+t_0=x(R_v)+t_0\neq0$. By Lemma~\ref{lem:ev}(b) and \eqref{eq:cycliccocycle},
\[
o_v=\bigl[(\sigma,\tau)\mapsto(x(R_v)+t_0)^{\psi(\sigma)\psi(\tau)}\bigr]=\bigl(\psi_v,\ x(R_v)+t_0\bigr).
\]
Definition~\ref{def:PS} gives the first formula and Fact~\ref{fact:schmid} the second. The class $o_v$ is determined by
$f'$ (Lemma~\ref{lem:ev}(a)), so the term at $v$ does not depend on $R_v$; and $o_v=0$ for almost all $v$ by
Lemma~\ref{lem:BMgood} below, whose proof does not use the present theorem.
\end{proof}

\begin{remark}[interpretation]\label{rem:BMint}
(1) If $\chi=\chi'$ then $\psi=0$ and all terms vanish: the alternating property is evident in this form.

(2) Only the trace of $\beta_0$ enters. If $t_0\neq0$, then $u=\beta_0/t_0$ satisfies $u^2+u=b/t_0^2$ and
$\sigma u-u=\chi'(\sigma)$, so
\begin{equation}\label{eq:t0}
[\,b/t_0^2\,]=\chi'\quad\text{in }K/\wp K ;
\end{equation}
conversely any $t_0\in K^\times$ with \eqref{eq:t0} is $\beta_0+\bar\beta_0$ for $\beta_0=t_0u$, which satisfies
Lemma~\ref{lem:norm}. (The case $t_0=0$ occurs exactly when $\beta_0^2=b$.) Thus the norm equation
$\Nm_{L'/K}\beta_0=b$ can be replaced by the Artin--Schreier condition \eqref{eq:t0} on $t_0\in K^\times$.

(3) Heuristically, Theorem~\ref{thm:BM} is an evaluation of a Brauer class in the sense of Manin: $f'$ is the cocycle
of the quaternion algebra $(\psi,f)$ over $K(X)$, which is the pull-back by $\pi\colon X\to E'$ of
$(\theta+\chi',\,x_{E'}+t_0)$, where $\theta\in H^1(E',\Z/2)$ is the class of the \'etale double cover
$\varphi\colon E\to E'$ (so $\theta(R)=\delta_\varphi(R)$). The function $x_{E'}+t_0$ has a double pole at $O$ and
vanishes exactly at $\varphi(\pm Q_0)$, where, if $t_0\neq0$, $\theta$ specialises to $\chi'$ by \eqref{eq:t0} (if $t_0=0$
the zero is a double zero at $(0,\sqrt b)$); this is why the algebra is unramified. We do not use this interpretation.
\end{remark}

\begin{remark}[a variant with a rational divisor]\label{rem:fisher}
Suppose that $\mathfrak b$ is a $K$-rational divisor of degree $0$ on $X$ whose class in $\Pic^0(X_\s)=E(\Ks)$ is $T$,
and that $f_1\in K(X)^\times$ satisfies $\dv f_1=2\mathfrak b$. Then $a''_\sigma=\chi'(\sigma)\mathfrak b$ lifts the
cocycle $\chi'(\sigma)T$, and $(da'')(\sigma,\tau)=2\chi'(\sigma)\chi'(\tau)\mathfrak b
=\dv\bigl(f_1^{\chi'(\sigma)\chi'(\tau)}\bigr)$, so the same argument as in Steps~3 and~4 gives
\[
\langle\iota(\chi),\iota(\chi')\rangle=\sum_v\inv_v\bigl(\chi'_v,\ f_1(P_v)\bigr),\qquad P_v\in X(K_v)\smallsetminus\operatorname{supp}(\dv f_1).
\]
This is the form of the formula used by Fisher \cite[\S3]{Fi} over number fields, where $f_1$ is called a pushout
function. A divisor $\mathfrak b$ exists because $X$ is everywhere locally soluble (the obstruction to representing a
rational divisor class by a rational divisor lies in $\Br K$ and vanishes locally), but finding one requires solving a
conic over $K$; the divisor $\Delta$ used above is only rational over the quadratic extension cut out by $\psi$, and is
given directly by $\beta_0$. We do not pursue this variant.
\end{remark}

\begin{remark}[the split case]\label{rem:splitorient}
Suppose $b=b_0^2$ with $b_0\in K^\times$ and take $\beta_0=b_0$, so $t_0=0$ and $f=(x+b_0)^2/x$. Then $x(Q_0)^2=b$, so
$2Q_0=T$ by (E4) and $-Q_0=Q_0+T$; hence $\mathcal Z=\bar{\mathcal Z}$ and $\Delta=\bar\Delta$. Moreover
$\varphi(Q_0)$ has $x$-coordinate $b_0+b/b_0=0$, i.e.\ $\varphi(Q_0)=T':=(0,b_0)\in E'(K)$, and
$\pi^*(T')=(Q_0)+(Q_0+T)$, $\pi^*(O)=(O)+(T)$; so $\Delta=\pi^*\bigl((T')-(O)\bigr)$ is $K$-rational, $f=\pi^*x_{E'}$,
and $\dv f=2\Delta$. This is Remark~\ref{rem:fisher} with $\mathfrak b=\Delta$ and $f_1=f$, and it gives the second
formula of Corollary~\ref{cor:split} directly. In general $\Delta$ is rational only over the field cut out by $\psi$,
and $\beta_0$ replaces the missing square root of $b$.
\end{remark}

\subsection{Consequences}

We can now prove Corollary~\ref{cor:bilin}.

\begin{proof}
The map $\iota$ is a homomorphism and the pairing is bilinear (Corollary~\ref{cor:bilinear}), so $B$ is bilinear. By
Theorem~\ref{thm:BM}, $B(\chi,\chi)=0$, because $\psi=\chi+\chi=0$; so $B$ is alternating, hence symmetric, its values being killed by $2$. If
$\chi'=\delta_\varphi(R)$ with $R\in E'(K)$ then $\iota(\chi')=0$, so $B(\chi,\chi')=B(\chi',\chi)=0$ for all $\chi$.
Finally, if $x\in\Sha(E/K)[2]$ and $x'=2y$ with $y\in\Sha(E/K)$, then
$\langle x,x'\rangle=\langle x,2y\rangle=\langle2x,y\rangle=0$; and if $x\in\Sha(E/K)[\varphi]$ and $x'=\hat\varphi_*y$
with $y\in\Sha(E'/K)$, then $\langle x,x'\rangle=\langle\varphi_*x,y\rangle_{E'}=0$ by Proposition~\ref{prop:F}
\textup{(}with $\rho=\varphi$, $\rho^\vee=\hat\varphi$\textup{)}.
\end{proof}

\begin{corollary}[split case]\label{cor:split}
Suppose $b=b_0^2$ with $b_0\in K^\times$. For $\chi,\chi'\in\Sel^\varphi(E/K)$ choose $R_v\in E'(K_v)$ with
$R_v\neq O$, $\delta_\varphi(R_v)=\chi_v$ and $x(R_v)\neq0$. Then
\[
\langle\iota(\chi),\iota(\chi')\rangle=\sum_v\inv_v\bigl(\chi_v+\chi'_v,\ x(R_v)\bigr)
=\sum_v\inv_v\bigl(\chi'_v,\ x(R_v)\bigr)
=\sum_v\Tr_{k_v/\F_2}\Res_v\Bigl(\alpha'_v\,\frac{dx(R_v)}{x(R_v)}\Bigr).
\]
\end{corollary}

\begin{proof}
Take $\beta_0=b_0$ in Lemma~\ref{lem:norm} (it is $G_K$-fixed and $b/b_0=b_0$); then $t_0=0$ and the first equality
is Theorem~\ref{thm:BM}. For the second, let $X=x(R_v)$ and $c=b_0/X\in K_v^\times$. By Lemma~\ref{lem:deltaphi},
$\chi_v=[b/X^2]=[c^2]=[c]$. If $u$ is a root of $u^2+u=c$, then $c=u(u+1)$ is a norm from $K_v(u)$, so $(\chi_v,c)=0$
by Fact~\ref{fact:norm} (trivially if $\chi_v=0$), and therefore $(\chi_v,X)=(\chi_v,b_0)$. Since
$\sum_v\inv_v(\chi_v,b_0)=0$ by Fact~\ref{fact:ABHN}, the terms $(\chi_v,x(R_v))$ do not contribute to the sum. The
last expression is Fact~\ref{fact:schmid}.
\end{proof}

By symmetry the pairing is also $\sum_v\inv_v(\chi_v,x(R'_v))$ with $\delta_\varphi(R'_v)=\chi'_v$; this is the shape of
the classical formula for $2$-isogenies in odd characteristic, where $x(R'_v)$ represents the Kummer class of a local
point.

For computations one needs to know where the local terms of Theorem~\ref{thm:BM} vanish.

\begin{lemma}\label{lem:BMgood}
Let $v$ be a place with $a_2\in\mathcal O_v$, $b\in\mathcal O_v^\times$, $t_0\in\mathcal O_v$ and $\chi'$ unramified at
$v$. Then $\bigl(\delta_\varphi(R)+\chi'_v,\ x(R)+t_0\bigr)=0$ for every $R\in E'(K_v)$ with $R\neq O$ and
$x(R)\neq t_0$.
\end{lemma}

The proof uses purity for the Brauer group and the following observation.

\begin{lemma}\label{lem:dvrsymbol}
Let $\mathcal O$ be a discrete valuation ring containing $\F_2$, with fraction field $L$ and residue field $k$. Let
$\psi\in H^1(\mathcal O,\Z/2)$ have trivial image in $H^1(k,\Z/2)$, and let $f\in L^\times$. Then $(\psi,f)\in\Br\mathcal O$.
\end{lemma}

\begin{proof}
As $\mathcal O$ is local, $H^1(\mathcal O,\mathbb G_a)=0$ and the Artin--Schreier sequence gives
$H^1(\mathcal O,\Z/2)=\mathcal O/\wp(\mathcal O)$; so $\psi=\chi_c$ with $c\in\mathcal O$, and as the image of $c$ in $k$
lies in $\wp(k)$ we may, after adding $\wp(a)$ for a lift $a$ of a root, assume $c\in\mathfrak m$. Let $\pi$ be a
uniformiser. If $c=0$ there is nothing to prove. Otherwise, replacing $c$ by $c+\wp(\pi)$ if $v(c)\geq2$, we may assume
$c=\pi w$ with $w\in\mathcal O^\times$. For $c'\in\mathcal O$ and $u\in\mathcal O^\times$ the class $(\chi_{c'},u)$ is that
of the Azumaya algebra $\mathcal O\langle i,j\rangle/(i^2+i-c',\ j^2-u,\ ji-(i+1)j)$, so lies in $\Br\mathcal O$. In $L$,
$c=i(i+1)$ is a norm from $L(i)$, $i^2+i=c$, so $(\chi_c,c)=0$ (Fact~\ref{fact:norm}) and
$(\chi_c,\pi)=(\chi_c,w)\in\Br\mathcal O$. Writing $f=\pi^nu$ with $u\in\mathcal O^\times$,
$(\psi,f)=n(\chi_c,\pi)+(\chi_c,u)\in\Br\mathcal O$.
\end{proof}

\begin{proofof}{Lemma~\textup{\ref{lem:BMgood}}}
Let $\mathcal E$, $\mathcal E'$ be the Weierstrass models of $E$, $E'$ over $\mathcal O_v$; they are smooth, as their
discriminants $b^2$ and $b$ are units. By the N\'eron mapping property \cite[1.2/8]{BLR} $\varphi$ extends to an isogeny $\mathcal E\to\mathcal E'$,
whose kernel is the constant group scheme generated by the section $T=(0,b)$ (which is disjoint from $O$); so
$\mathcal E\to\mathcal E'$ is a $\Z/2$-torsor and defines $\theta\in H^1(\mathcal E',\Z/2)$, whose value at a point
$P$ of $E'$ is $\delta_\varphi(P)$ over the residue field of $P$. As $\chi'$ is unramified,
$\psi:=\theta+\chi'_v\in H^1(\mathcal E',\Z/2)$. Let $f=x+t_0$ and $\alpha=(\psi,f)\in\Br K_v(E')$.

The section $s\colon\operatorname{Spec}\mathcal O_v\to\mathcal E'$ extending $R$ maps the closed point to a closed
point $\mathfrak m$ of the special fibre, and since every point of its image specialises to $\mathfrak m$, it factors
through $\operatorname{Spec}B$, $B=\mathcal O_{\mathcal E',\mathfrak m}$, a regular local ring of dimension $2$. By
\cite[Prop.~7.4]{AG}, $\Br B=\bigcap_{\mathfrak p}\Br B_{\mathfrak p}$ inside the Brauer group of the function field,
the intersection over the primes $\mathfrak p$ of height one of $B$. These are the generic point of the special fibre,
which is irreducible, and the closed points of $E'_{K_v}$ whose closure contains $\mathfrak m$. We check that
$\alpha\in\Br B_{\mathfrak p}$ for each such $\mathfrak p$. If $\mathfrak p$ is the generic point of the
special fibre, $f$ is a unit at $\mathfrak p$, because its reduction $\tilde x+\tilde t_0$ is a non-zero function on the
special fibre; as $\psi$ is defined over $B_{\mathfrak p}$, $\alpha\in\Br B_{\mathfrak p}$. Otherwise $\mathfrak p$ is a
closed point $P$ of $E'_{K_v}$. If $P$ is not a zero or pole of $f$, again $\alpha\in\Br B_{\mathfrak p}$. At $P=O$ the
function $f$ has a pole of order $2$, and $(\psi,f)=(\psi,f\pi_P^2)$ with $f\pi_P^2$ a unit. If $t_0=0$, then $\beta_0^2=b$ and $\sigma\beta_0=b/\beta_0=\beta_0$ whenever $\chi'(\sigma)=1$, so $\beta_0\in K$
and $b$ is a square in $K$; the only zero of $f=x$ is then the rational point $(0,\beta_0)$, of multiplicity $2$, and the
same argument applies. If $t_0\neq0$ and $x(P)=t_0$, then
$\psi(P)=\delta_\varphi(P)+\chi'=[b/t_0^2]+\chi'=0$ over the residue field of $P$, by Lemma~\ref{lem:deltaphi} and
\eqref{eq:t0}, and Lemma~\ref{lem:dvrsymbol} applies. Hence $\alpha\in\Br B$, and $s^*\alpha\in\Br\mathcal O_v=0$. As $R$ is neither a zero nor
a pole of $f$, the restriction of $s^*\alpha$ to $K_v$ is $(\psi(R),f(R))=(\delta_\varphi(R)+\chi'_v,x_R+t_0)$.
\end{proofof}

Lemma~\ref{lem:BMgood} makes precise the heuristic of Remark~\ref{rem:BMint}(3): the algebra $(\theta+\chi',x_{E'}+t_0)$ is
unramified on the smooth model of $E'$ at $v$.

\section{Examples}\label{sec:example}

\subsection{A split example: proof of Theorem~\ref{thm:example}}\label{subsec:split}
Let $k\geq1$, $q=2^{4k}$, so that $\F_{16}\subseteq\F_q$. Let $\zeta_3\in\F_4$ with $\zeta_3^2+\zeta_3+1=0$ and
\[
\Lambda(S)=S^4+\zeta_3S^2+\zeta_3^2S=S(S+1)(S^2+S+\zeta_3^2)\in\F_4[S].
\]
Since $\Tr_{\F_4/\F_2}(\zeta_3^2)=1$, the roots $\theta,\theta+1$ of $S^2+S+\zeta_3^2$ lie in $\F_{16}\smallsetminus\F_4$; the
polynomial $\Lambda$ is additive and $Z:=\{a\in\F_q:\Lambda(a)=0\}=\{0,1,\theta,\theta+1\}$ is an $\F_2$-subspace. Write $\Tr$ for
$\Tr_{\F_q/\F_2}$, and $Z^\perp=\{\lambda\in\F_q:\Tr(\lambda a)=0\ \forall a\in Z\}$, a subspace of dimension $4k-2$.

Let $K=\F_q(S)$ and consider \eqref{eq:E} with $a_2=0$ and $b=\Lambda^2$:
\[
E\colon y^2+xy=x^3+\Lambda^4,\qquad E'\colon y^2+xy=x^3+\Lambda^2,\qquad T=(0,\Lambda^2).
\]
For $\lambda\in\F_q$ let $\chi_\lambda\in H^1(K,\Z/2)$ be the Artin--Schreier class of $\lambda S$.

\emph{Why the example works.} Since $\Lambda$ is additive, its roots form the $\F_2$-subspace $Z$, and since $b=\Lambda^2$
is a square, Corollary~\ref{cor:split} applies with $t_0=0$. The bad places of $E$ are $\infty$ and the roots of $\Lambda$.
The class of $\lambda S$ is trivial at the place $S=a$ if and only if $\Tr(\lambda a)=0$, and it is unramified at all
other finite places, so for $\lambda\in Z^\perp$ only the Selmer condition at $\infty$ remains. At $\infty$ the
explicit points $R_\lambda$ below exist when $\Tr(\lambda)=0$, which holds on $Z^\perp$ because $1\in Z$. The pairing
then reduces to a single residue at $\infty$, which turns out to be the trace form. The field $\F_q$ is taken to
contain $\F_{16}$ so that $Z\subseteq\F_q$ and $\Tr(1)=\Tr(\zeta_3)=0$.

\begin{proof}
\emph{Traces.} For $x\in\F_{16}$, $\Tr(x)=k\Tr_{\F_{16}/\F_2}(x)$. We have $\Tr_{\F_{16}/\F_2}(1)=0$ and
$\theta^4=(\theta+\zeta_3^2)^2=\theta+\zeta_3^2+\zeta_3=\theta+1$, so $\Tr_{\F_{16}/\F_4}(\theta)=\theta+\theta^4=1$ and
$\Tr_{\F_{16}/\F_2}(\theta)=\Tr_{\F_4/\F_2}(1)=0$; likewise $\Tr_{\F_{16}/\F_2}(\theta^2)=0$ and
$\Tr_{\F_{16}/\F_2}(\theta(\theta+1))=\Tr_{\F_{16}/\F_2}(\zeta_3^2)=\Tr_{\F_4/\F_2}(2\zeta_3^2)=0$. As $aa'$ for $a,a'\in Z$ is an
$\F_2$-combination of $1,\theta,\theta^2,\theta(\theta+1)$, $Z\subseteq Z^\perp$. Moreover $\Tr(1)=4k\cdot1=0$ and
$\Tr(\zeta_3)=\Tr_{\F_4/\F_2}(2k\,\zeta_3)=0$, $\Tr(\zeta_3^2)=0$; and $\Tr(\lambda)=\Tr(\lambda\cdot1)=0$ for $\lambda\in Z^\perp$.

\emph{The local points at $\infty$.} Let $t=1/S$, $K_\infty=\F_q((t))$; recall $t\F_q[[t]]\subseteq\wp(K_\infty)$ and
$\alpha S^2\equiv\sqrt\alpha\,S\pmod{\wp K_\infty}$. For $\lambda\in\F_q$ put $c=c_\lambda=1+\lambda$ and
$X_\lambda=S^2+cS=S(S+c)$. Division gives
\begin{equation}\label{eq:division}
\frac{\Lambda}{X_\lambda}=\frac{S^3+\zeta_3S+\zeta_3^2}{S+c}=S^2+cS+c^2+\zeta_3+\frac{r_c}{S+c},\qquad r_c=c^3+\zeta_3c+\zeta_3^2,
\end{equation}
and $r_c=\Lambda(c)/c=0$ if $c\in Z\smallsetminus\{0\}$. Hence modulo $\wp(K_\infty)$:
$\Lambda^2/X_\lambda^2\equiv\Lambda/X_\lambda\equiv S+cS+c^2+\zeta_3=\lambda S+c^2+\zeta_3$ and $X_\lambda\equiv S+cS=\lambda S$, so
\[
X_\lambda+\Lambda^2/X_\lambda^2\equiv c^2+\zeta_3,\qquad \Tr(c^2+\zeta_3)=\Tr(1)+\Tr(\lambda)+\Tr(\zeta_3)=\Tr(\lambda).
\]
For $\lambda\in Z^\perp$ this is $0$, so $X_\lambda+\Lambda^2/X_\lambda^2\in\wp(K_\infty)$ and there is
$R_\lambda\in E'(K_\infty)$ with $x(R_\lambda)=X_\lambda$ (take $y=X_\lambda z$, $z^2+z=X_\lambda+\Lambda^2/X_\lambda^2$). By
Lemma~\ref{lem:deltaphi}, $\delta_\varphi(R_\lambda)=[X_\lambda]=[\lambda S]=\chi_{\lambda,\infty}$.

(a) We check $\chi_\lambda$ locally, $\lambda\in Z^\perp$. At $\infty$: $\chi_{\lambda,\infty}=\delta_\varphi(R_\lambda)$.
At the place $v_a$ given by $S=a$, $a\in Z$ (residue field $\F_q$): $\lambda S=\lambda a+\lambda(S+a)$ with
$\lambda(S+a)\in\wp(K_{v_a})$, so $\chi_{\lambda,v_a}$ is the class of the constant $\lambda a$, which is trivial because
$\Tr(\lambda a)=0$. At every other place $v$: $E$ has good reduction ($\Delta(E)=\Lambda^4$ is a unit), $\chi_{\lambda,v}$ is
unramified, so $\iota(\chi_\lambda)_v$ is the inflation of a class in $H^1(\Gal(K_v^{\mathrm{nr}}/K_v),E(K_v^{\mathrm{nr}}))$,
which vanishes \cite[I, Prop.~3.8]{Mi}. Hence $\iota(\chi_\lambda)\in\Sha(E/K)$.

(b) Let $\lambda,\mu\in Z^\perp$ and $\nu=\lambda+\mu$. Here $b=\Lambda^2$ is a square, so Corollary~\ref{cor:split} (with
$b_0=\Lambda$, $\chi=\chi_\mu$, $\chi'=\chi_\lambda$) gives
\[
\langle\iota(\chi_\mu),\iota(\chi_\lambda)\rangle=\sum_v\inv_v\bigl(\chi_{\nu,v},x(R_v)\bigr),\qquad R_v\in E'(K_v),\
\delta_\varphi(R_v)=\chi_{\mu,v},\ x(R_v)\neq0 .
\]
At $v_a$, $\chi_{\nu,v_a}=0$ because $\nu\in Z^\perp$. At a place $v\neq\infty$ of good reduction the hypotheses of
Lemma~\ref{lem:BMgood} hold with $t_0=0$ and $\chi'=\chi_\lambda$, which is unramified at $v$ as $\lambda S\in\mathcal O_v$;
so $\inv_v(\chi_{\nu,v},x(R_v))=0$. At $\infty$ take $R_\infty=R_\mu$, so that
$x(R_\mu)=S(S+c)$ with $c=1+\mu$:
\begin{align*}
\inv_\infty(\chi_\nu,X_\mu)&=\Tr\Res_\infty\Bigl(\nu S\Bigl(\frac{dS}{S}+\frac{dS}{S+c}\Bigr)\Bigr)
=\Tr\Res_\infty\Bigl(\nu\,dS+\nu\,dS+\frac{\nu c\,dS}{S+c}\Bigr)=\Tr(\nu c),
\end{align*}
using $\Res_\infty(dS)=0$ and $\Res_\infty(dS/(S+c))=-1=1$. Finally
$\Tr(\nu c)=\Tr\bigl((\lambda+\mu)(1+\mu)\bigr)=\Tr(\lambda)+\Tr(\mu)+\Tr(\lambda\mu)+\Tr(\mu^2)=\Tr(\lambda\mu)$, since
$\Tr(\mu^2)=\Tr(\mu)=0$ and $\Tr(\lambda)=0$. This proves (b).

(c) If $\lambda\in Z\smallsetminus\{1\}$, then $c=1+\lambda\in Z\smallsetminus\{0\}$, $r_c=0$, and
$\Lambda/X_\lambda=S^2+cS+c^2+\zeta_3\in\F_q[S]$ by \eqref{eq:division}. Then
$X_\lambda+\Lambda^2/X_\lambda^2=S^4+(1+c^2)S^2+cS+c^4+\zeta_3^2=z^2+z$ with
$z=S^2+cS+\gamma$, $\gamma^2+\gamma=c^4+\zeta_3^2$, and such $\gamma\in\F_q$ exists because $\Tr(c^4+\zeta_3^2)=\Tr(c)+\Tr(\zeta_3^2)=0$. So
$R_\lambda=(X_\lambda,X_\lambda z)\in E'(K)$ and $\chi_\lambda=[X_\lambda]=\delta_\varphi(R_\lambda)$ globally; hence
$\iota(\chi_\lambda)=0$. Also $\chi_1=\chi_\theta+\chi_{\theta+1}$, so $\iota(\chi_1)=0$. Conversely, if
$\iota(\chi_\lambda)=0$ then $\langle\iota(\chi_\mu),\iota(\chi_\lambda)\rangle=0$ for all $\mu\in Z^\perp$
(Corollary~\ref{cor:bilinear}), so by (b) $\Tr(\lambda\mu)=0$ for all $\mu\in Z^\perp$, i.e.\
$\lambda\in(Z^\perp)^\perp=Z$ (the trace form is non-degenerate).

(d) Suppose $\iota(\chi_\lambda)\in\hat\varphi\,\Sha(E'/K)$ with $\lambda\in Z^\perp$. By Corollary~\ref{cor:bilin},
$\langle\iota(\chi_\mu),\iota(\chi_\lambda)\rangle=0$ for all $\mu\in Z^\perp$, since $\iota(\chi_\mu)\in\Sha(E/K)[\varphi]$.
By (b), $\Tr(\lambda\mu)=0$ for all $\mu\in Z^\perp$, so $\lambda\in Z$. As $\lambda\mapsto\iota(\chi_\lambda)$ is
additive, this gives the injection $Z^\perp/Z\hookrightarrow\Sha(E/K)/\hat\varphi\,\Sha(E'/K)$; since
$2\,\Sha(E/K)\subseteq\hat\varphi\,\Sha(E'/K)$, it also gives the injection into $\Sha(E/K)/2\,\Sha(E/K)$.

Finally, by (c) the induced map $Z^\perp/Z\to\Sha(E/K)[\varphi]$ is injective; $\dim_{\F_2}Z^\perp/Z=4k-4$, and the
radical of $\Tr(\lambda\mu)$ on $Z^\perp$ is $Z^\perp\cap(Z^\perp)^\perp=Z$.
\end{proof}

\begin{remark}
For $k=1$ we get $Z^\perp=Z$ and no information. For $k=2$ the pairing exhibits a subgroup of order $16$ of
$\Sha(E/\F_{256}(S))[\varphi]$ with a non-degenerate pairing, injecting into $\Sha/2\Sha$; in fact this subgroup is all
of $\Sha(E/\F_{256}(S))$, and $E$ has rank $3$ (Proposition~\ref{prop:exsplit}). In this example
$E'[2]\simeq\Z/2\times\mu_2$ (Proposition~\ref{prop:B}) and $\beta_0=\Lambda$ is rational. A non-split example is given
in the next subsection.
\end{remark}

\subsection{A non-split example: proof of Theorem~\ref{thm:nonsplit}}\label{subsec:nonsplit}
Let $K=\F_8(S)$, let $\varrho\in\F_8$ with $\varrho^3+\varrho+1=0$, so that $\F_8^\times=\langle\varrho\rangle$ and
\[
\varrho^3=\varrho+1,\quad \varrho^4=\varrho^2+\varrho,\quad \varrho^5=\varrho^2+\varrho+1,\quad \varrho^6=\varrho^2+1,\quad
\varrho^7=1 .
\]
Write $\Tr=\Tr_{\F_8/\F_2}$; then $\Tr(1)=1$, $\Tr(\varrho)=\Tr(\varrho^2)=\Tr(\varrho^4)=0$ and
$\Tr(\varrho^3)=\Tr(\varrho^5)=\Tr(\varrho^6)=1$. Let
\[
b=S^4+S^3+\varrho,\qquad E\colon y^2+xy=x^3+b^2,\qquad E'\colon y^2+xy=x^3+b\qquad(a_2=0).
\]
Since $b$ is not a square in $K$, $E'[2]$ is a non-split extension of $\Z/2$ by $\mu_2$ (Proposition~\ref{prop:B}). The
roots of $\lambda^3+\lambda+1$ are $\varrho,\varrho^2,\varrho^4$. For $c\in K$ write $\chi_c$ for the
Artin--Schreier class of $c$; thus $\chi_{\varrho^4S}=\chi_{\varrho S}+\chi_{\varrho^2S}$.

\begin{proof}
Let $t=1/S$ and $K_\infty=\F_8((t))$. The polynomial $b$ has no root in $\F_{64}$. Indeed, a root $S_0$ would be
non-zero, and $u=1/S_0$ would satisfy $u^4+cu=c$ with $c=\varrho^{-1}$. The equation $z^3=\varrho$ has three solutions in
$\F_{64}$ and only one in $\F_8$; let $z\in\F_{64}\smallsetminus\F_8$ be a solution and $w=z/c$. Then
$w^{15}=\varrho^5c^{-15}=\varrho^{20}=c$, so $w^{16}=cw$, and with $\Tr'=\Tr_{\F_{64}/\F_2}$ we get
$\Tr'(wc)=\Tr'(w(u^4+cu))=\Tr'((w^{16}+cw)u)=0$. On the other hand $wc=z$ has minimal polynomial $x^2+\tau x+\nu$ over
$\F_8$ with $\nu=z^9=\varrho^3$, and $\varrho=z^3=(\tau^2+\nu)z+\tau\nu$ gives $\tau^2=\nu$, i.e.\ $\tau=\varrho^5$; so
$\Tr'(z)=\Tr(z+z^8)=\Tr(\varrho^5)=1$, a contradiction. Hence $b$ has no
factor of degree $1$ or $2$ over $\F_8$, so it is irreducible, and its zeros form one place $v_b$ of degree $4$,
with residue field $\F_{8^4}$ in which $S$ maps to a root $a$ of $b$; the sum of the roots of $b$ is $1$. The bad
places of $E$ are $v_b$ and $\infty$. For $s_0\in\F_8$ we write $v_{s_0}$ for the place $S=s_0$; its residue field is
$\F_8$, and for a polynomial $c\in\F_8[S]$ the class $\chi_{c,v_{s_0}}$ is the class of the constant $c(s_0)$, which is
trivial if and only if $\Tr(c(s_0))=0$.

(a) At $v_b$: $\chi_{1,v_b}$ and $\chi_{\lambda S,v_b}$ are the classes of the constants $1$ and $\lambda a$ of the residue
field, and $\Tr_{\F_{8^4}/\F_2}(1)=0$, $\Tr_{\F_{8^4}/\F_2}(\lambda a)=\Tr(\lambda\cdot1)=0$ for
$\lambda\in\{\varrho,\varrho^2\}$; so these classes are trivial at $v_b$. At the other finite places $E$ has good
reduction and the classes are unramified, so they satisfy the Selmer condition as in the proof of
Theorem~\ref{thm:example}(a). At $\infty$ we exhibit local points (Lemma~\ref{lem:deltaphi}): for $X=S^2+S+1$,
$b/X^2=1+t+O(t^2)$, so $X+b/X^2\equiv(S^2+S+1)+1\equiv0$ and $X\equiv1$ modulo $\wp(K_\infty)$; hence there is
$R_1\in E'(K_\infty)$ with $x(R_1)=S^2+S+1$ and $\delta_\varphi(R_1)=\chi_{1,\infty}$. For $\lambda^3+\lambda+1=0$ and
$X=\lambda S$, $b/X^2=\lambda^{-2}(S^2+S)+\varrho\lambda^{-2}t^2$, so
$X+b/X^2\equiv(\lambda+\lambda^{-1}+\lambda^{-2})S=\lambda^{-2}(\lambda^3+\lambda+1)S=0$; hence there is
$R_\lambda\in E'(K_\infty)$ with $x(R_\lambda)=\lambda S$ and $\delta_\varphi(R_\lambda)=\chi_{\lambda S,\infty}$.

(b) By Corollary~\ref{cor:bilin} the matrix is symmetric with zero diagonal, so it suffices to compute the three
entries above the diagonal. We use Theorem~\ref{thm:BM} with the following norm elements. For
$\chi'=\chi_{\varrho^2S}$, $L'=K(y)$ with $y^2+y=\varrho^2S$, and
\[
\beta_0=A+By,\qquad A=S^2+\varrho^4S+1,\qquad B=\varrho^3(S+1),\qquad t_0=\Tr_{L'/K}\beta_0=B=\varrho^3(S+1),
\]
satisfies $\Nm_{L'/K}\beta_0=A^2+AB+\varrho^2SB^2=b$. For $\chi'=\chi_{\varrho S}$, $L'=K(y)$ with $y^2+y=\varrho S$, and
\[
\beta_0=A+By,\qquad A=S^2+S+1,\qquad B=\varrho S+\varrho^3,\qquad t_0=B=\varrho(S+\varrho^2),
\]
satisfies $\Nm_{L'/K}\beta_0=A^2+AB+\varrho SB^2=b$. (Both identities are verified by expanding, using
$\varrho^3=\varrho+1$; for instance in the second case $A^2=S^4+S^2+1$, $AB=\varrho S^3+S^2+S+\varrho^3$ and
$\varrho SB^2=\varrho^3S^3+S$.) Such norm elements are found by a small search: by Remark~\ref{rem:BMint}(2) it
suffices to find $t_0\in K^\times$ with $b/t_0^2+\alpha'\in\wp(K)$, where $\alpha'=\varrho^2S$, resp.\ $\varrho S$;
writing $b/t_0^2+\alpha'=\wp(A/t_0)$ one gets $\beta_0=A+t_0y$. Above, $t_0=B$ is linear.

\emph{The finite places.} Let $\chi\in\{\chi_1,\chi_{\varrho S}\}$, let $\chi'\in\{\chi_{\varrho S},\chi_{\varrho^2S}\}$ with
$t_0$ as above, and $\psi=\chi+\chi'$. We claim that the term of Theorem~\ref{thm:BM} at every finite place $v$ vanishes,
for every admissible $R_v$. At $v_b$, $\chi_{v_b}=\chi'_{v_b}=0$ by (a), so $\psi_{v_b}=0$. At a finite place
$v\neq v_b$ we have $a_2=0$, $b\in\mathcal O_v^\times$, $t_0\in\mathcal O_v$ and $\chi'$ unramified, so the term
vanishes by Lemma~\ref{lem:BMgood}. This proves the claim, so in each case the pairing is the term at $\infty$.

\emph{The place $\infty$.} We use the points $R_1$, $R_\varrho$ of (a) and Fact~\ref{fact:schmid}; note $dS=t^{-2}dt$.

\smallskip\noindent
$\chi=\chi_1$, $\chi'=\chi_{\varrho^2S}$: $x(R_1)+t_0=S^2+(1+\varrho^3)S+(1+\varrho^3)=S^2+\varrho S+\varrho$, and
$\psi$ is the class of $1+\varrho^2S$. Then
\[
(1+\varrho^2S)\,\frac{d(S^2+\varrho S+\varrho)}{S^2+\varrho S+\varrho}
=\frac{\varrho\,(1+\varrho^2t^{-1})}{1+\varrho t+\varrho t^2}\,dt
\]
has residue $\varrho^3$, and $\Tr(\varrho^3)=1$: the entry is $1$.

\smallskip\noindent
$\chi=\chi_1$, $\chi'=\chi_{\varrho S}$: $x(R_1)+t_0=S^2+(1+\varrho)S+(1+\varrho^3)=S^2+\varrho^3S+\varrho$, and $\psi$ is
the class of $1+\varrho S$. Then
\[
(1+\varrho S)\,\frac{d(S^2+\varrho^3S+\varrho)}{S^2+\varrho^3S+\varrho}
=\frac{\varrho^3\,(1+\varrho t^{-1})}{1+\varrho^3t+\varrho t^2}\,dt
\]
has residue $\varrho^4$, and $\Tr(\varrho^4)=0$: the entry is $0$.

\smallskip\noindent
$\chi=\chi_{\varrho S}$, $\chi'=\chi_{\varrho^2S}$: $x(R_\varrho)+t_0=(\varrho+\varrho^3)S+\varrho^3=S+\varrho^3$, and $\psi$ is
the class of $(\varrho+\varrho^2)S=\varrho^4S$. Then
\[
\varrho^4S\,\frac{dS}{S+\varrho^3}=\frac{\varrho^4\,t^{-2}}{1+\varrho^3t}\,dt=\varrho^4\bigl(t^{-2}+\varrho^3t^{-1}+\cdots\bigr)dt
\]
has residue $\varrho^7=1$, and $\Tr(1)=1$: the entry is $1$. This proves (b).

(c) By (b), $B(\chi_1,\chi_{\varrho^2S})=\frac12$, so $\iota(\chi_1)\neq0$ and $\iota(\chi_{\varrho^2S})\neq0$
(Corollary~\ref{cor:bilinear}), and $B(\chi_1,\chi_1+\chi_{\varrho^2S})=\frac12$, so
$\iota(\chi_1)+\iota(\chi_{\varrho^2S})\neq0$. Hence the two classes span a group of order $4$, on which the pairing
has matrix $\left(\begin{smallmatrix}0&1\\1&0\end{smallmatrix}\right)$. Every non-zero element $x'$ of this group
satisfies $\langle x,x'\rangle\neq0$ for some $x$ in the group, and $x\in\Sha(E/K)[\varphi]$; so
$x'\notin\hat\varphi\,\Sha(E'/K)$ by Corollary~\ref{cor:bilin}.
\end{proof}

\begin{remark}
The radical of the matrix in Theorem~\ref{thm:nonsplit}(b) is spanned by $\chi_1+\chi_{\varrho S}$. In fact
$\chi_1,\chi_{\varrho S},\chi_{\varrho^2S}$ span $\Sel^\varphi(E/K)$, $\chi_1+\chi_{\varrho S}$ is the Kummer class of a
point of infinite order of $E'(K)$, and $\Sha(E/K)$ is the group of order $4$ of (c) (Proposition~\ref{prop:exnonsplit});
its $2$-primary part is determined by descent alone.
One entry of the matrix is computed with Theorem~\ref{thm:W} in \S\ref{subsec:exW}.
\end{remark}

\section{The flat Weil-pairing form}\label{sec:W}

In this section we prove Theorem~\ref{thm:W}. The plan is as follows. In \S\ref{subsec:kummerhat} we identify the
Kummer map of $\hat\varphi$ with $\kappa_F$ and prove a local orthogonality statement. In \S\ref{subsec:thmA} we prove,
again directly from Definition~\ref{def:PS} but with a different lift $f'$, that the Cassels--Tate pairing is
$\sum_v\inv_v(\chi_v,g_v)$ for explicit classes $g_v\in K_v^\times/K_v^{\times2}$ (Theorem~\ref{thm:A}). In
\S\ref{sec:fppf} we describe $H^1_\fl(F,E'[2])$, identify $g_v$ with the difference of a global lift and a local
Kummer class, and identify $(\chi_v,g_v)$ with a flat cup product. In \S\ref{subsec:exW} we work out an example, and in \S\ref{subsec:kramer} we compare with Kramer's pairing.
Subsections~\ref{subsec:kummerhat} and~\ref{subsec:thmA} use only Galois cohomology, apart from the interpretation of
$\kappa_F$ in Proposition~\ref{prop:kummerF}; flat cohomology enters in \S\ref{sec:fppf}.

\smallskip\noindent\emph{Flat cohomology in characteristic $2$.} The group scheme
$\mu_2=\Spec F[\zeta]/(\zeta^2-1)=\Spec F[z]/(z^2)$ ($\zeta=1+z$) is infinitesimal: $\mu_2(\Fs)=1$, so Galois
cohomology does not see it, and we use cohomology $H^i_\fl$ on the fppf site. Squaring $\mathbb G_m\to\mathbb G_m$ is
surjective for the fppf topology (adjoining a square root is a flat extension) but not for the \'etale topology. The
Kummer sequence $1\to\mu_2\to\mathbb G_m\xrightarrow{2}\mathbb G_m\to1$ and $H^1_\fl(F,\mathbb G_m)=\Pic(\Spec F)=0$
give $H^1_\fl(F,\mu_2)=F^\times/F^{\times2}$, the class of $u$ being the torsor $\Spec F[e]/(e^2-u)$, and
$H^2_\fl(F,\mu_2)=\Br(F)[2]$. A second feature of characteristic $2$ is that $\Fs$ is not closed under square roots, so
$\Fs{}^\times/\Fs{}^{\times2}\neq1$; this group, with its Galois action, is where the classes $[\beta_0]$ and
$\kappa(P)$ live (Lemma~\ref{lem:squares}, Proposition~\ref{prop:B}). Smooth group schemes have the same cohomology in
the fppf and \'etale topologies \cite[Thm.~11.7]{Gro}, so only the finite group schemes $\mu_2$ and $E'[2]$ require
flat cohomology.

\subsection{The Kummer map of \texorpdfstring{$\hat\varphi$}{phi-hat} and local isotropy}\label{subsec:kummerhat}
Recall that for a field $F\supseteq K$, $\kappa_F\colon E(F)\to F^\times/F^{\times2}$ is the homomorphism of
Lemma~\ref{lem:elem}\textup{(E6)}.

\begin{proposition}[Kummer map of $\hat\varphi$]\label{prop:kummerF}
Let $\delta_{\hat\varphi}\colon E(F)\to H^1_\fl(F,\ker\hat\varphi)$ be the connecting map of the fppf sequence
$0\to\ker\hat\varphi\to E'\xrightarrow{\hat\varphi}E\to0$. Under $\ker\hat\varphi\simeq\mu_2$
\textup{(}Lemma~\ref{lem:isog}\textup{)} and the Kummer isomorphism $H^1_\fl(F,\mu_2)=F^\times/F^{\times2}$, one has
$\delta_{\hat\varphi}=\kappa_F$; that is, $\delta_{\hat\varphi}(P)=[x(P)]$ for $P\notin\{O,T\}$ and
$\delta_{\hat\varphi}(T)=[b]$.
\end{proposition}

The proof is given in Appendix~\ref{app:flat}. The idea is the following. In characteristic $2$ an action of
$\mu_2=\Spec F[z]/(z^2)$, with comultiplication $z\mapsto z\otimes1+1\otimes z+z\otimes z$, on an affine scheme
$\Spec A$ is the same as a $\Z/2$-grading of $A$, i.e.\ a derivation $D$ of $A$ with $D^2=D$, and a $\mu_2$-torsor
$\Spec A$ is the Kummer torsor of $u$ when $A$ is generated by an element of weight one with square $u$. For the
translation action of $\ker\hat\varphi$ on $E'$ the derivation is $D(f)=df/(dx/x)$, so $D(x)=x$: the coordinate $x$
has weight one, and on the fibre $\hat\varphi^{-1}(P)$ it is a square root of $x(P)$.

\begin{lemma}[local isotropy; cf.\ {\cite[Prop.~1.6]{Kr}}]\label{lem:iso}
Let $F\supseteq K$ be a field, $\chi\in H^1(F,\Z/2)$, and suppose there is $P\in E(\Fs)$ with $\sigma P=P+\chi(\sigma)T$
for all $\sigma\in G_F$ \textup{(}i.e.\ $\chi=\delta_\varphi(\varphi P)$\textup{)}. Then $(\chi,\kappa_F(Q))=0$ in
$\Br(F)$ for every $Q\in E(F)$.
\end{lemma}

\begin{proof}
We may assume $\chi\neq0$; let $L=F(\chi)$ be the quadratic extension it defines. Then $P\in E(L)\smallsetminus E(F)$,
so $P\notin\{O,T\}$, and for $\sigma\notin\ker\chi$ we get $\Nm_{L/F}x(P)=x(P)\,x(P+T)=b$ by (E2). This proves the
claim for $Q=T$ by the norm criterion (Fact~\ref{fact:norm}); $Q=O$ is trivial. Let $Q\notin\{O,T\}$. Then
$P\neq\pm Q$ and $P+Q\notin\{O,T\}$ (otherwise $P$ would be $F$-rational), and $P+Q$ satisfies the same twisted
relation as $P$, so $\Nm_{L/F}x(P+Q)=b$. By (E5) over $L$, $x(P+Q)=x(P)x(Q)w^2$ with $w\in L^\times$. Taking norms,
$b=b\,x(Q)^2\Nm(w)^2$, so $(x(Q)\Nm(w))^2=1$ and $x(Q)=\Nm(w)^{-1}$ is a norm from $L$. Hence $(\chi,x(Q))=0$.
\end{proof}

\begin{lemma}\label{lem:squares}
Let $F$ be a field of characteristic $2$. Then $F^\times\cap(\Fs{}^{\times})^2=F^{\times2}$, the map
$F^\times/F^{\times2}\to\Fs{}^\times/\Fs{}^{\times2}$ is injective, and $H^0(G_F,\Fs{}^\times/\Fs{}^{\times2})
=F^\times/F^{\times2}$.
\end{lemma}

\begin{proof}
If $z\in\Fs$ and $z^2\in F$ then $z$ is a root of $X^2-z^2=(X-z)^2$, so $z$ is purely inseparable and separable over
$F$, i.e.\ $z\in F$. For the last claim use the exact sequence
$1\to\Fs{}^\times\xrightarrow{2}\Fs{}^\times\to\Fs{}^\times/\Fs{}^{\times2}\to1$ (squaring is injective) and Hilbert's
Theorem~90, $H^1(G_F,\Fs{}^\times)=0$.
\end{proof}

\subsection{The formula with local norm classes}\label{subsec:thmA}
We abbreviate $\kappa=\kappa_F$ when the field is clear.

\begin{theorem}\label{thm:A}
Let $\chi,\chi'\in\Sel^\varphi(E/K)$ and let $\beta_0$ be as in Lemma~\ref{lem:norm}.
\begin{enumerate}[label=\textup{(\alph*)},leftmargin=*]
\item For each place $v$ let $P'_v\in E(\Kvs)$ satisfy $\sigma P'_v=P'_v+\chi'(\sigma)T$ for all $\sigma\in G_v$
\textup{(}such points exist; equivalently $R'_v:=\varphi(P'_v)\in E'(K_v)$ and $\delta_\varphi(R'_v)=\chi'_v$\textup{)}.
Then $[\beta_0]\,\kappa(P'_v)^{-1}\in\Kvs{}^\times/\Kvs{}^{\times2}$ is the image of a unique
$g_v\in K_v^\times/K_v^{\times2}$. It depends only on $\chi'$, $\beta_0$ and $P'_v$, and not on $\chi$.
\item One has
\[
\langle\iota(\chi),\iota(\chi')\rangle=\sum_v\inv_v(\chi_v,g_v)
=\sum_v\Tr_{k_v/\F_2}\Res_v\Bigl(\alpha_v\,\frac{dg_v}{g_v}\Bigr),
\]
where $\alpha_v\in K_v$ is any element with $\chi_v=\chi_{\alpha_v}$. Almost all terms vanish. The class $(\chi_v,g_v)$ does not
depend on the choice of $P'_v$, and the sum does not depend on the choice of $\beta_0$.
\end{enumerate}
\end{theorem}

\begin{remark}\label{rem:hyp}
If $\chi'_v=0$ one may take $P'_v\in E(K_v)$, and then $g_v=[\beta_0]\kappa(P'_v)^{-1}$ with $\beta_0\in K_v$; by
Lemma~\ref{lem:iso} (applied over $K_v$ to $\chi_v$ and a point $P$ of the torsor of $\iota(\chi)_v$, which exists
since $\iota(\chi)_v=0$) the term is then simply $\inv_v(\chi_v,\beta_0)$.
\end{remark}

The proof occupies \S\S\ref{subsec:A1b}--\ref{subsec:A6}. We put $\Gamma=(\Z/2)^2$, write its elements as
$1,\ga,\gb,\gc$ with $(\chi,\chi')$-values $(0,0),(1,0),(0,1),(1,1)$, and for a place $v$ let
$\Pi_v=(\chi,\chi')|_{G_v}\colon G_v\to\Gamma$.

\subsubsection{The classes \texorpdfstring{$g_v$}{g\_v}; proof of (a)}\label{subsec:A1b}
Let $v$ be a place and $P'_v$ as in (a); such a point exists since $\iota(\chi')_v=0$. If $\chi'_v\neq0$ then, as in
the proof of Lemma~\ref{lem:norm}, $P'_v\in E(L'_v)\smallsetminus E(K_v)$ for the quadratic extension $L'_v$ cut out by
$\chi'_v$, and $\sigma x(P'_v)=b/x(P'_v)$ for $\sigma$ with $\chi'(\sigma)=1$.
Suppose first $P'_v\notin\{O,T\}$, and put $H'=x(P'_v)$ and
\[
\gamma_v:=\beta_0/H'\in\Kvs{}^\times .
\]
For $\sigma\in G_v$ we have $\sigma H'=H'$ or $b/H'$, and $\sigma\beta_0=\beta_0$ or $b/\beta_0$, according as
$\chi'(\sigma)=0$ or $1$; hence $\sigma\gamma_v=\gamma_v$ or $\gamma_v^{-1}$ accordingly. If $\chi'_v=0$ then
$\gamma_v\in K_v^\times$. If $\chi'_v\neq0$ then $\gamma_v\in L'_v$ has norm $1$, so by Hilbert's Theorem~90
\begin{equation}\label{eq:H90}
\gamma_v=w_0/\bar w_0=w_0^2/\Nm(w_0),\qquad w_0\in L'^\times_v,\ \bar w_0=\sigma_1w_0\ (\chi'(\sigma_1)=1).
\end{equation}
In both cases $[\gamma_v]\in\Kvs{}^\times/\Kvs{}^{\times2}$ comes from $K_v^\times/K_v^{\times2}$: from $\gamma_v$ itself,
resp.\ from $\Nm(w_0)^{-1}$. If $P'_v\in\{O,T\}$ then $\chi'_v=0$ and $[\beta_0]\kappa(P'_v)^{-1}$ is $[\beta_0]$ or
$[\beta_0/b]$, with $\beta_0\in K_v$. Uniqueness of $g_v$ is Lemma~\ref{lem:squares}.

\subsubsection{The global lift}\label{subsec:A2}
Let $X$ be the torsor of $\iota(\chi)$ and $h:=x$, as in \S\ref{subsec:torsor}; recall \eqref{eq:twisth}. We lift the
cocycle $\sigma\mapsto\chi'(\sigma)T$ representing $\iota(\chi')$ to the divisor cochain
\[
a''_\sigma=\chi'(\sigma)D,\qquad D=(T)-(O).
\]
Unlike in \S\ref{subsec:proofBM} we use this naive divisor, which is $K$-rational on $E$ and satisfies
$\sigma^*D=\pm D$ on $X$. The price is that $\eta$ is the divisor of a cochain of functions $\mathcal F$ only up to a
constant $3$-cocycle \eqref{eq:dPhi}, which is cancelled by a constant correction $\varepsilon$ built from $\beta_0$.
Thus $\beta_0$ enters here through constants, whereas in \S\ref{subsec:proofBM} it enters through the divisor
$\Delta$.
As $\sigma^*D=(-1)^{\chi(\sigma)}D$ and $\chi'(\sigma\tau)=\chi'(\sigma)+\chi'(\tau)-2\chi'(\sigma)\chi'(\tau)$,
\[
\eta(\sigma,\tau):=\sigma^*a''_\tau-a''_{\sigma\tau}+a''_\sigma=2\,m(\sigma,\tau)\,D=\dv\bigl(h^{m(\sigma,\tau)}\bigr),
\qquad m(\sigma,\tau):=\chi'(\tau)\bigl(\chi'(\sigma)-\chi(\sigma)\bigr),
\]
using (E1). Put $\mathcal F(\sigma,\tau):=h^{m(\sigma,\tau)}$. Since $d\eta=0$, $d\mathcal F$ has trivial divisor, hence is a
constant $3$-cocycle; by \eqref{eq:twisth} its value is the power of $b$ contributed by $\sigma^*\mathcal F(\tau,\rho)$:
\begin{equation}\label{eq:dPhi}
(d\mathcal F)(\sigma,\tau,\rho)=b^{\chi(\sigma)m(\tau,\rho)}=b^{\chi(\sigma)\chi'(\rho)(\chi'(\tau)-\chi(\tau))} .
\end{equation}
Define the constant $2$-cochain
\[
\varepsilon(\sigma,\tau):=(\sigma\tau\bar\beta_0)^{\chi(\sigma)\chi(\tau)}\,(\sigma\beta_0)^{-\chi(\sigma)\chi'(\tau)} .
\]

\begin{lemma}\label{lem:global}
$d\varepsilon=d\mathcal F$. Consequently $f':=\mathcal F/\varepsilon\in Z^2(G_K,\Ks(X)^\times)$ and $\dv f'=\eta$, so $f'$ is a lift
as in Definition~\ref{def:PS}.
\end{lemma}

\begin{proof}
For a $\Z$-valued $p$-cochain $\alpha$ (trivial action) and a cochain $\vartheta$ with values in a multiplicative module put
$(\alpha\cup\vartheta)(\sigma_1,\dots,\sigma_{p+q})=\bigl(\sigma_1\cdots\sigma_p\,\vartheta(\sigma_{p+1},\dots)\bigr)^{\alpha(\sigma_1,\dots,\sigma_p)}$;
then $d(\alpha\cup\vartheta)=(d\alpha\cup\vartheta)\cdot(\alpha\cup d\vartheta)^{(-1)^p}$ \cite[V.3]{Br}. Let
$\tilde\chi\colon G_K\to\{0,1\}\subset\Z$ be the lift of $\chi$ and $Y(\sigma,\tau)=\chi(\sigma)\chi(\tau)$. Then
$d\tilde\chi=2Y$ and $dY=0$ (direct check from $\chi(\sigma\tau)=\chi(\sigma)+\chi(\tau)-2\chi(\sigma)\chi(\tau)$).
Let $\vartheta(\tau)=\beta_0^{\chi'(\tau)}$, and regard $\bar\beta_0$ as a $0$-cochain. From Lemma~\ref{lem:norm}:
$(d\vartheta)(\tau,\rho)=b^{\chi'(\tau)\chi'(\rho)}=:\omega'(\tau,\rho)$ and $(d\bar\beta_0)(\rho)=(\beta_0/\bar\beta_0)^{\chi'(\rho)}$.
By definition $\varepsilon=(Y\cup\bar\beta_0)\cdot(\tilde\chi\cup\vartheta)^{-1}$, so
\[
d\varepsilon=(Y\cup d\bar\beta_0)\cdot(d\tilde\chi\cup\vartheta)^{-1}\cdot(\tilde\chi\cup\omega') .
\]
Here $(Y\cup d\bar\beta_0)(\sigma,\tau,\rho)\,(2Y\cup\vartheta)(\sigma,\tau,\rho)^{-1}
=\bigl(\sigma\tau(\beta_0/\bar\beta_0)\cdot\sigma\tau(\beta_0)^{-2}\bigr)^{Y(\sigma,\tau)\chi'(\rho)}
=b^{-\chi(\sigma)\chi(\tau)\chi'(\rho)}$ because $\beta_0\bar\beta_0=b\in K$, and
$(\tilde\chi\cup\omega')(\sigma,\tau,\rho)=b^{\chi(\sigma)\chi'(\tau)\chi'(\rho)}$. Comparing with \eqref{eq:dPhi}
gives $d\varepsilon=d\mathcal F$. The remaining assertions follow since $\varepsilon$ is constant.
\end{proof}

\subsubsection{The local classes}\label{subsec:A3}
Fix $v$. As $\iota(\chi)_v=0$, $X(K_v)\neq\emptyset$: there is $P_v\in E(\Kvs)$ with $\sigma P_v=P_v+\chi(\sigma)T$ for
$\sigma\in G_v$. Replacing $P_v$ by $P_v+Q$ with $Q\in E(K_v)$ (an infinite group) we may assume, given $P'_v$ as in
(a) with $P'_v\notin\{O,T\}$,
\begin{equation}\label{eq:generic}
P_v\notin\{O,T\},\qquad P_v\neq\pm P'_v,\qquad P_v+P'_v\notin\{O,T\}.
\end{equation}
The values of $f'$ are $h^{0,\pm1}$ times constants, hence regular and invertible at $P_v$, and evaluation at $P_v$ is
$G_v$-equivariant for the twisted action ($(\sigma^*f)(P_v)=(\sigma^*f)(\sigma^*P_v)=\sigma(f(P_v))$). By
Lemma~\ref{lem:ev}(b) applied to the curve $X_{K_v}$ and the rational point $P_v$, with $H:=x(P_v)$,
\begin{equation}\label{eq:localclass}
o_v=[c(H,\beta_0)],\qquad c(H,\beta_0)(\sigma,\tau):=H^{m(\sigma,\tau)}\,(\sigma\beta_0)^{\chi(\sigma)\chi'(\tau)}\,
(\sigma\tau\bar\beta_0)^{-\chi(\sigma)\chi(\tau)}\qquad(\sigma,\tau\in G_v).
\end{equation}
Note $\sigma H=b/H$ if $\chi(\sigma)=1$ and $\sigma H=H$ otherwise.

\subsubsection{Two local cochain identities}\label{subsec:A4}
Keep $v$, $P_v$, $P'_v\notin\{O,T\}$ as in \eqref{eq:generic}; put $H'=x(P'_v)$, $J=x(P_v+P'_v)$, and let
$w\in\Kvs{}^\times$ be the element of (E5) with $J=HH'w^2$. Since $\sigma(P_v+P'_v)=P_v+P'_v+(\chi+\chi')(\sigma)T$,
$\sigma J=J$ or $b/J$ according as $(\chi+\chi')(\sigma)=0$ or $1$. Comparing $\sigma J$ with
$\sigma H\cdot\sigma H'\cdot(\sigma w)^2$ and using that squaring is injective in characteristic $2$, the action of
$G_v$ on $H,H',w$ depends only on $\Pi_v(\sigma)$ and is given by Table~\ref{tab:action}.

\begin{table}[h]\small
\centering
\begin{tabular}{c|ccccc}
\toprule
$\Pi_v(\sigma)$ & $H$ & $H'$ & $w$ & $w_0$ & $\bar w_0$\\
\midrule
$\ga=(1,0)$ & $b/H$ & $H'$ & $1/(H'w)$ & $w_0$ & $\bar w_0$\\
$\gb=(0,1)$ & $H$ & $b/H'$ & $1/(Hw)$ & $\bar w_0$ & $w_0$\\
$\gc=(1,1)$ & $b/H$ & $b/H'$ & $HH'w/b$ & $\bar w_0$ & $w_0$\\
\bottomrule
\end{tabular}
\caption{The action of $\sigma\in G_v$ through $\Pi_v(\sigma)\in\Gamma$ ($b$ is fixed). The last two columns refer to
\eqref{eq:H90}.}\label{tab:action}
\end{table}

Let $M$ be the free abelian group on symbols $b,H,H',w,w_0,\bar w_0$, with the $\Gamma$-action of
Table~\ref{tab:action} (one checks $\ga^2=\gb^2=1$ and $\ga\gb=\gb\ga=\gc$ on generators, e.g.\
$\ga(\gb(w))=\ga(1/(Hw))=HH'w/b=\gc(w)$). The evaluation map $M\to\Kvs{}^\times$ is $G_v$-equivariant when $G_v$ acts on $M$
through $\Pi_v$ (when $\chi'_v=0$ the symbols $w_0,\bar w_0$ are not used). Hence it suffices to verify identities of
cochains on $\Gamma$ with values in $M$, which are then pulled back along $\Pi_v$ and evaluated. All cochains below are
normalised (value $1$ if one argument is $1$), so only the nine pairs $(\sigma,\tau)\in\{\ga,\gb,\gc\}^2$ matter.

The two identities have the following meaning. By \eqref{eq:localclass} the local class $o_v$ is built from
$H=x(P_v)$ and the global element $\beta_0$. Identity (B) says that if $\beta_0$ is replaced by the $x$-coordinate $H'$
of the local point $P'_v$, which satisfies the same twisted relation locally, the resulting cochain is a coboundary:
a lift coming from a local point contributes nothing. Identity (C) says that the discrepancy $\nu=\beta_0/H'$
contributes exactly the cyclic algebra $(\chi_v,g_v)$. Table~\ref{tab:B} is the verification: each entry is the common
value of the two sides.

\smallskip\noindent\textbf{Identity (B).} Let $c(H,H')$ be defined by \eqref{eq:localclass} with $\beta_0$ replaced by
the symbol $H'$ (so $\bar\beta_0$ is replaced by $b/H'$), and let $\kappa_B(\ga)=1/(Hw)$, $\kappa_B(\gb)=1/w$,
$\kappa_B(1)=\kappa_B(\gc)=1$. Then $c(H,H')=d\kappa_B$; both sides are given by Table~\ref{tab:B}.

\begin{table}[h]\small
\centering
\begin{tabular}{c|ccc}
\toprule
$\sigma\backslash\tau$ & $\ga$ & $\gb$ & $\gc$\\
\midrule
$\ga$ & $H'/b$ & $H'/H$ & $1/H$\\
$\gb$ & $1$ & $H$ & $H$\\
$\gc$ & $1/H'$ & $b/H'$ & $1$\\
\bottomrule
\end{tabular}\qquad
\begin{tabular}{c|ccc}
\toprule
$\sigma\backslash\tau$ & $\ga$ & $\gb$ & $\gc$\\
\midrule
$\ga$ & $\nu$ & $\nu$ & $1$\\
$\gb$ & $1$ & $1$ & $1$\\
$\gc$ & $\nu^{-1}$ & $\nu^{-1}$ & $1$\\
\bottomrule
\end{tabular}
\caption{Left: $c(H,H')=d\kappa_B$. Right: $c_\nu$, where $\nu=w_0/\bar w_0$; it equals $(\chi,\Nm w_0)\,d\kappa_C$.}
\label{tab:B}
\end{table}

For instance, $m(\ga,\gb)=-1$, so $c(H,H')(\ga,\gb)=H^{-1}\cdot\ga(H')=H'/H$, while
\[
d\kappa_B(\ga,\gb)=\ga(\kappa_B(\gb))\,\kappa_B(\ga)/\kappa_B(\gc)=\ga(1/w)/(Hw)=H'w/(Hw)=H'/H;
\]
and, as $\gc\ga=\gb$, $c(H,H')(\gc,\ga)=\bigl(\gb(b/H')\bigr)^{-1}=1/H'$, while $\gc(Hw)=(b/H)(HH'w/b)=H'w$ gives
$d\kappa_B(\gc,\ga)=\gc(\kappa_B(\ga))\,\kappa_B(\gc)/\kappa_B(\gb)=w/(H'w)=1/H'$. The other seven entries are equally direct.
Note that (B) uses the relation $J=HH'w^2$, i.e.\ the group law of $E$.

\smallskip\noindent\textbf{Identity (C).} Put $\nu:=\gamma_v=\beta_0/H'$. Writing $\beta_0=H'\nu$ and
$\bar\beta_0=(b/H')\nu^{-1}$ in \eqref{eq:localclass} gives
\[
c(H,\beta_0)=c(H,H')\cdot c_\nu,\qquad c_\nu(\sigma,\tau):=(\sigma\nu)^{\chi(\sigma)\chi'(\tau)}\,
(\sigma\tau\nu)^{\chi(\sigma)\chi(\tau)} .
\]
If $\chi'_v=0$, then $\nu\in K_v^\times$ is $G_v$-fixed and $c_\nu(\sigma,\tau)=\nu^{\chi(\sigma)\chi(\tau)}$ is the cocycle
\eqref{eq:cycliccocycle} of $(\chi_v,\nu)$. If $\chi'_v\neq0$, write $\nu=w_0/\bar w_0$ as in \eqref{eq:H90} and let
$\kappa_C(\ga)=1/\bar w_0$, $\kappa_C(\gc)=1/w_0$, $\kappa_C(1)=\kappa_C(\gb)=1$; then
\[
c_\nu=(\chi,\Nm w_0)\cdot d\kappa_C,\qquad \Nm w_0=w_0\bar w_0 ,
\]
where $(\chi,\Nm w_0)(\sigma,\tau)=(\Nm w_0)^{\chi(\sigma)\chi(\tau)}$ is nontrivial exactly for
$\sigma,\tau\in\{\ga,\gc\}$. E.g.\ at $(\ga,\ga)$: $c_\nu(\ga,\ga)=\nu$ and
$\Nm w_0\cdot\ga(1/\bar w_0)/\bar w_0=w_0/\bar w_0$; at $(\gc,\gc)$: $c_\nu(\gc,\gc)=\gc(\nu)\nu=1$ and
$\Nm w_0\cdot\gc(1/w_0)/w_0=1$; the remaining entries are checked in the same way.

\subsubsection{Conclusion of the proof of (b)}\label{subsec:A5}
Let $P'_v\notin\{O,T\}$. By (B) and (C), $o_v=[c(H,H')c_\nu]=[c_\nu]$ equals $(\chi_v,\nu)$ if $\chi'_v=0$ and
$(\chi_v,\Nm w_0)$ if $\chi'_v\neq0$. By \eqref{eq:H90} and the definition of $g_v$, $g_v=[\nu]$, resp.\
$g_v=[\Nm w_0]^{-1}$; since $(\chi_v,-)$ is $2$-torsion, $o_v=(\chi_v,g_v)$ in both cases. If $P'_v\in\{O,T\}$ (so
$\chi'_v=0$), then $g_v$ differs from the value for any other choice $P''_v\in E(K_v)\smallsetminus\{O,T\}$ by
$\kappa(P''_v-P'_v)$ with $P''_v-P'_v\in E(K_v)$, which is orthogonal to $\chi_v$ by Lemma~\ref{lem:iso} (applied with
$P=P_v$). Hence in all cases
\[
o_v=(\chi_v,g_v),
\]
and Definition~\ref{def:PS} gives $\langle\iota(\chi),\iota(\chi')\rangle=\sum_v\inv_v(\chi_v,g_v)$; the residue form is
Fact~\ref{fact:schmid}. Since $o_v$ is determined by $f'$, i.e.\ by $\beta_0$, the term $(\chi_v,g_v)$ does not depend on
$P'_v$. (Directly: replacing $P'_v$ by $P'_v+Q$, $Q\in E(K_v)$, multiplies $g_v$ by $\kappa(Q)$, and
$(\chi_v,\kappa(Q))=0$ by Lemma~\ref{lem:iso}.) Replacing $\beta_0$ by another solution $\beta_0\rho$, with
$\Nm\rho=1$, i.e.\ $\rho=r/\bar r$ ($r\in L'^\times$) by Hilbert's Theorem~90, multiplies every $g_v$ by the global element
$\Nm(r)^{-1}\in K^\times$ (modulo squares), which changes the sum by $\sum_v\inv_v(\chi_v,\Nm r)=0$
(Fact~\ref{fact:ABHN}); if $\chi'=0$, $\beta_0$ changes by an element of $K^\times$ and the same argument applies.

\subsubsection{Finiteness}\label{subsec:A6}
Let $\Sigma$ be a finite set of places containing those where $\chi$ or $\chi'$ ramifies, where $a_2$ or $b$ is not
integral or $b$ is not a unit, and those below which $\beta_0$ or $\bar\beta_0$ is not a unit. Let $v\notin\Sigma$; then
$\chi_v$ is unramified, so $\inv_v(\chi_v,g)=0$ whenever $v(g)$ is even (Fact~\ref{fact:schmid}). If $\chi'_v\neq0$,
$L'_v/K_v$ is unramified and $g_v=[\Nm w_0]^{-1}$ has even valuation (norms from an unramified quadratic extension have
even valuation). If $\chi'_v=0$, the term is $\inv_v(\chi_v,\beta_0)$ by Remark~\ref{rem:hyp}, and
$\beta_0\in\mathcal O_v^\times$. Hence the term at $v$ vanishes. This completes the proof of
Theorem~\ref{thm:A}. \qed

\begin{lemma}[a formula for $g_v$]\label{lem:trace}
In Theorem~\ref{thm:A}(a) suppose $P'_v\notin\{O,T\}$ and put $\gamma_v=\beta_0/x(P'_v)$.
\begin{enumerate}[label=\textup{(\alph*)},leftmargin=*]
\item If $\chi'_v=0$ then $g_v=[\gamma_v]$. If $\chi'_v\neq0$ and $\gamma_v\neq1$ then
$g_v=[\Tr_{L'_v/K_v}(\gamma_v)]$, where $L'_v=K_v(\chi'_v)$.
\item Suppose $\chi'_v\neq0$, $\chi'_v=\chi_{\alpha'}$ with $\alpha'\in K_v$, so $L'_v=K_v(y)$ with $y^2+y=\alpha'$, and
write $\beta_0=A+By$ with $A,B\in K_v$. Let $X=x(R'_v)\neq0$ and let $h\in K_v$ satisfy $h^2+h=b/X^2+\alpha'$ \textup{(}it
exists since $[b/X^2]=\delta_\varphi(R'_v)=[\alpha']$, Lemma~\ref{lem:deltaphi}\textup{)}. Choose $P'_v$ with
$x(P'_v)=X(y+h)$. If $\gamma_v\neq1$, then $g_v=[\,b\,X\,(A+Bh)\,]$.
\end{enumerate}
\end{lemma}

\begin{proof}
(a) The first case is the definition. In the second, put $w_0=1+\gamma_v$ and let $\sigma_1$ generate
$\Gal(L'_v/K_v)$. Since $\sigma_1\gamma_v=\gamma_v^{-1}$, $\sigma_1w_0=1+\gamma_v^{-1}$, so $w_0/\sigma_1w_0=\gamma_v$
and $\Nm(w_0)=(1+\gamma_v)(1+\gamma_v^{-1})=\gamma_v+\gamma_v^{-1}=\Tr(\gamma_v)$ (characteristic $2$), which is non-zero
as $\gamma_v\neq1$. As in \eqref{eq:H90}, $g_v=[\Nm w_0]^{-1}=[\Tr\gamma_v]$.
(b) The roots of $U^2+XU+b$ are the $x$-coordinates of the two points of $\varphi^{-1}(R'_v)=\{P'_v,P'_v+T\}$ (proof of
Lemma~\ref{lem:deltaphi}); they are $X(y+h)$ and $X(y+h+1)$, so the choice is possible, and
$\sigma_1x(P'_v)=x(P'_v+T)=X(y+1+h)$. Hence $\Tr(\gamma_v)=\Tr\bigl(\beta_0\,\sigma_1x(P'_v)\bigr)/\Nm(x(P'_v))$ with
$\Nm(x(P'_v))=b$. Using $\Tr(y)=1$ and $y^2=y+\alpha'$,
$\Tr\bigl((A+By)\,X(y+1+h)\bigr)=X\bigl(A+B+B(1+h)\bigr)=X(A+Bh)$. So $g_v=[X(A+Bh)/b]=[bX(A+Bh)]$.
\end{proof}

\begin{remark}[comparison with Theorem~\ref{thm:BM}]\label{rem:BMcompare}
The lifts $f'$ of \S\ref{subsec:A2} and of \S\ref{subsec:proofBM} belong to different divisor cochains $a''$, and the two
families of local terms, $(\chi_v,g_v)$ and $(\psi_v,x(R_v)+t_0)$, need not a priori agree place by place; their sums
agree because both compute the pairing of Definition~\ref{def:PS} (Remark~\ref{rem:PSuse}). They do agree at every
place $v$ with $\chi'_v=0$: there $\beta_0\in K_v^\times$, and if $\chi_v\neq0$, $M=K_v(\chi_v)$ and $P_v\in E(M)$ is a
point of $X(K_v)$ with $\varphi(P_v)=R_v$, then by (E2)
\[
\Nm_{M/K_v}\bigl(x(P_v)+\beta_0\bigr)=\bigl(x(P_v)+\beta_0\bigr)\bigl(b/x(P_v)+\beta_0\bigr)=\beta_0\,\bigl(x(R_v)+t_0\bigr),
\]
so $(\chi_v,x(R_v)+t_0)=(\chi_v,\beta_0)$ by Fact~\ref{fact:norm}, which is the term of Theorem~\ref{thm:A} by
Remark~\ref{rem:hyp}. We do not know whether the local terms agree in general when $\chi'_v\neq0$; in the example of
\S\ref{subsec:exW} they do. The formula of Theorem~\ref{thm:A} is
visibly additive in $\chi$, since the classes $g_v$ depend only on $\chi'$, whereas that of Theorem~\ref{thm:BM} is
visibly alternating.
\end{remark}

\subsection{The flat-cohomological form}\label{sec:fppf}

In this subsection $F$ is any field of characteristic $2$ with $a_2\in F$, $b\in F^\times$. By
Lemma~\ref{lem:isog}(iii), $E'[2]=\ker(\varphi\circ\hat\varphi)=\hat\varphi^{-1}(\ker\varphi)$; as $\hat\varphi$ is finite
flat and surjective, it maps $E'[2]$ fppf-surjectively onto $\ker\varphi=\{O,T\}\simeq\Z/2$, with kernel
$\ker\hat\varphi\simeq\mu_2$. Write $\kappa=\kappa_{\Fs}\colon E(\Fs)\to\Fs{}^\times/\Fs{}^{\times2}$.

\begin{lemma}\label{lem:kappaonto}
$\kappa$ is surjective with kernel $\hat\varphi(E'(\Fs))$, and $\hat\varphi$ is injective on $E'(\Fs)$.
\end{lemma}

\begin{proof}
For $u\in\Fs{}^\times$, $y^2+uy=u^3+a_2u^2+b^2$ is separable in $y$, so there is a point with $x=u$. If $\kappa(P)=1$ with
$P\notin\{O,T\}$, write $x(P)=u^2$, $u\in\Fs$, and let $y_0\in\Fs$ be a root of the separable polynomial
$y_0^2+uy_0+u^3+a_2u^2+b$; then $(u,y_0)\in E'(\Fs)$ and $\hat\varphi(u,y_0)$ has $x$-coordinate $u^2$, so
$P=\pm\hat\varphi(u,y_0)\in\hat\varphi(E'(\Fs))$. If $\kappa(T)=[b]=1$, $b=c^2$ with $c\in\Fs$ and $T=\hat\varphi(0,c)$.
Conversely $\kappa\circ\hat\varphi=1$ because $x(\hat\varphi(R))=x(R)^2$ and $\hat\varphi(0,c)=T$ with $c^2=b$.
Injectivity: $\ker\hat\varphi(\Fs)=\mu_2(\Fs)=1$.
\end{proof}

Kramer \cite[\S1]{Kr} describes $H^1_\fl(F,E'[2])$ by an explicit group of pairs; the following model is better
adapted to our cochain computations.

\begin{proposition}\label{prop:B}
\begin{enumerate}[label=\textup{(\roman*)},leftmargin=*]
\item The sequence $0\to\mu_2\xrightarrow{j}E'[2]\xrightarrow{\hat\varphi}\Z/2\to0$ of fppf group schemes is exact,
and the fibre $\hat\varphi^{-1}(T)\subset E'[2]$ over the generator is a $\mu_2$-torsor of class $[b]$. In particular
the extension splits if and only if $b\in F^{\times2}$.
\item Let $\mathcal D(F)$ be the group of pairs $(\chi',[\beta])$ with $\chi'\in H^1(F,\Z/2)$,
$[\beta]\in\Fs{}^\times/\Fs{}^{\times2}$, and $\sigma[\beta]=[\beta][b]^{\chi'(\sigma)}$ for all $\sigma\in G_F$. There is an
isomorphism
\[
\Phi\colon H^1_\fl(F,E'[2])\xrightarrow{\ \sim\ }\mathcal D(F)/\langle(0,[b])\rangle
\]
such that: (a) $\mathrm{pr}_1\circ\Phi=\hat\varphi_*\colon H^1_\fl(F,E'[2])\to H^1(F,\Z/2)$; (b) for
$u\in F^\times/F^{\times2}=H^1_\fl(F,\mu_2)$, $\Phi(j_*u)=(0,[u])$; (c) for $R\in E'(F)$,
$\Phi(\delta_2R)=(\delta_\varphi R,\kappa(Q))$ for any $Q\in E(\Fs)$ with $\varphi(Q)=R$, where $\delta_2$ is the Kummer map
of $[2]_{E'}$; (d) if $\chi'\neq0$, $L'=F(\chi')$ and $b=\Nm_{L'/F}\beta_0$, or if $\chi'=0$ and $\beta_0\in F^\times$, then
$(\chi',[\beta_0])\in\mathcal D(F)$, so it defines a lift $\xi=\Phi^{-1}(\chi',[\beta_0])$ of $\chi'$.
\end{enumerate}
\end{proposition}

\begin{proof}
(i) Exactness was explained above. The fibre $\hat\varphi^{-1}(T)$ lies in $E'[2]$ (as $\varphi(T)=O$) and is the torsor
$\delta_{\hat\varphi}(T)$, of class $[b]$ by Proposition~\ref{prop:kummerF}. An extension of $\Z/2$ by $\mu_2$ splits if
and only if the fibre over the generator has an $F$-point.

(ii) As $[2]\colon E'\to E'$ is fppf surjective, $E'[2]$ is quasi-isomorphic to the complex $[E'\xrightarrow{2}E']$ (in
degrees $0,1$), and since $E'$ is smooth its fppf cohomology is Galois cohomology \cite[Thm.~11.7]{Gro}. Hence
$H^1_\fl(F,E'[2])$ is the first Galois hypercohomology of $E'(\Fs)\xrightarrow{2}E'(\Fs)$: classes of pairs $(\xi,W)$
with $\xi\in Z^1(G_F,E'(\Fs))$, $W\in E'(\Fs)$, $2\xi_\sigma=\sigma W-W$, modulo pairs $(dU,2U)$. Similarly
$\ker\varphi\simeq[E\xrightarrow{\varphi}E']$ and $\mu_2\simeq[E'\xrightarrow{\hat\varphi}E]$.

Given $(\xi,W)$ choose $Q\in E(\Fs)$ with $\varphi(Q)=W$. Then $\varphi(\hat\varphi(\xi_\sigma)-(\sigma Q-Q))=0$, so
$\hat\varphi(\xi_\sigma)-(\sigma Q-Q)=\chi'(\sigma)T$ with $\chi'(\sigma)\in\{0,1\}$; $\chi'$ is a homomorphism (a direct
computation using the cocycle relation of $\xi$), and it is the image of $(\xi,W)$ under the map of complexes
$(\hat\varphi,\mathrm{id})\colon[E'\xrightarrow2E']\to[E\xrightarrow{\varphi}E']$. Put $\Phi(\xi,W)=(\chi',\kappa(Q))$. Since
$\kappa$ is a homomorphism killing $\hat\varphi(E'(\Fs))$ (Lemma~\ref{lem:kappaonto}),
$\sigma\kappa(Q)=\kappa(Q+\hat\varphi(\xi_\sigma)-\chi'(\sigma)T)=\kappa(Q)[b]^{\chi'(\sigma)}$, so the value lies in
$\mathcal D(F)$. Replacing $Q$ by $Q+T$ multiplies $\kappa(Q)$ by $[b]$; replacing $(\xi,W)$ by $(\xi+dU,W+2U)$ is absorbed
by $Q\mapsto Q+\hat\varphi(U)$, which changes neither $\chi'$ nor $\kappa(Q)$. So $\Phi$ is a well-defined homomorphism,
and (a) holds by construction. (b): a class of $H^1_\fl(F,\mu_2)$ is represented by $(\zeta,P)$, $\zeta\in Z^1(G_F,E'(\Fs))$,
$P\in E(\Fs)$, $\hat\varphi(\zeta_\sigma)=\sigma P-P$; $j_*$ is induced by $(\mathrm{id},\varphi)$, giving
$(\zeta,\varphi(P))$, and with $Q=P$ we get $\Phi=(0,\kappa(P))$. The class of $(\zeta,P)$ restricts over $\Fs$ to that of
$(0,P)$, i.e.\ to $\delta_{\hat\varphi}(P)=\kappa(P)$ (Proposition~\ref{prop:kummerF}); since
$F^\times/F^{\times2}\to\Fs{}^\times/\Fs{}^{\times2}$ is injective (Lemma~\ref{lem:squares}), $\kappa(P)$ is the image of
$u$. (c): $\delta_2(R)$ is the class of $(0,R)$. (d): if $\chi'(\sigma)=1$ then $\sigma\beta_0=b/\beta_0$, so
$\sigma[\beta_0]=[b][\beta_0]$.

Surjectivity: given $(\chi',[\beta])\in\mathcal D(F)$, choose $Q$ with $\kappa(Q)=[\beta]$ (Lemma~\ref{lem:kappaonto}) and
put $W=\varphi(Q)$. Then $\kappa(\sigma Q-Q-\chi'(\sigma)T)=1$, so $\sigma Q-Q-\chi'(\sigma)T=\hat\varphi(\xi_\sigma)$ for a
unique $\xi_\sigma\in E'(\Fs)$; uniqueness gives the cocycle relation, and
$2\xi_\sigma=\varphi\hat\varphi(\xi_\sigma)=\sigma W-W$. Thus $\Phi(\xi,W)=(\chi',[\beta])$. Injectivity: if $\Phi(c)=0$,
then $\hat\varphi_*c=0$, so $c=j_*(u)$ for some $u$, and by (b) $[u]\in\{1,[b]\}$ in $\Fs{}^\times/\Fs{}^{\times2}$, hence
in $F^\times/F^{\times2}$; and $j_*([b])=0$ because $[b]$ is the image of the generator of $H^0(F,\Z/2)$ under the
connecting map, by (i). So $c=0$.
\end{proof}

The cup product $H^1(F,\Z/2)\times H^1_\fl(F,\mu_2)\to H^2_\fl(F,\mu_2)$ is induced by the Weil pairing
$\Z/2\times\mu_2\to\mu_2$. The Kummer sequence and $H^1_\fl(F,\mathbb G_m)=0$ identify $H^2_\fl(F,\mu_2)$ with
$H^2_\fl(F,\mathbb G_m)[2]=\Br(F)[2]$. Cup products in fppf cohomology are recalled in Appendix~\ref{app:flat}.

\begin{lemma}\label{lem:cup}
For $\chi\in H^1(F,\Z/2)$ and $u\in F^\times$, the image of $\chi\cup\delta_{\mathrm{Kum}}(u)\in H^2_\fl(F,\mu_2)$ in
$\Br(F)[2]$ is $(\chi,u)$ \textup{(}up to sign, which is irrelevant\textup{)}.
\end{lemma}

The proof, a chain-level computation with the complexes $[\Z\xrightarrow2\Z]$ and $[\mathbb G_m\xrightarrow2\mathbb G_m]$,
is given in Appendix~\ref{app:flat}.

\begin{remark}
The proof in Appendix~\ref{app:flat} does not use a general compatibility of cup products with connecting maps on the fppf site; it only uses
the comparison with Galois cohomology for complexes of smooth groups and an explicit chain-level pairing. Only
Corollary~\ref{cor:Aprime} depends on Lemma~\ref{lem:cup}; Theorem~\ref{thm:A} is stated with cyclic algebras and does
not.
\end{remark}

\begin{corollary}[Weil-pairing form]\label{cor:Aprime}
In the situation of Theorem~\ref{thm:A}, let $\xi\in H^1_\fl(K,E'[2])$ be the lift $\Phi^{-1}(\chi',[\beta_0])$ of
$\chi'$ and $R'_v=\varphi(P'_v)\in E'(K_v)$. Then $\xi_v-\delta_2(R'_v)=j_*(g_v)$ in $H^1_\fl(K_v,E'[2])$, and
\[
\langle\iota(\chi),\iota(\chi')\rangle=\sum_v\inv_v\bigl(\chi_v\cup g_v\bigr),
\]
the cup product being $H^1(K_v,\Z/2)\times H^1_\fl(K_v,\mu_2)\to H^2_\fl(K_v,\mu_2)\subset\Br K_v$. The class
$g_v\in H^1_\fl(K_v,\mu_2)$ is determined by $\xi_v-\delta_2(R'_v)$ up to $j_*$-kernel $\langle[b]\rangle$, and changing
$R'_v$ changes it by $\kappa(E(K_v))$; both are orthogonal to $\chi_v$ (Lemma~\ref{lem:iso}).
\end{corollary}

\begin{proof}
The construction of $\Phi$ is compatible with restriction from $K$ to $K_v$ (use $E(\Ks)\subseteq E(\Kvs)$). By
Proposition~\ref{prop:B}(c),(d), $\Phi(\xi_v-\delta_2R'_v)=(\chi'_v,[\beta_0])-(\chi'_v,\kappa(P'_v))=(0,[\beta_0]\kappa(P'_v)^{-1})=(0,g_v)$,
which is $\Phi(j_*g_v)$ by (b). Now apply Theorem~\ref{thm:A} and Lemma~\ref{lem:cup}.
\end{proof}

\begin{remark}
This is the Weil-pairing definition of \cite[\S12.2]{PS}, \cite[Def.~1.2]{BF} for the pair $(\varphi,\hat\varphi)$: lift
the class $\chi'\in H^1(K,\ker\varphi)$ through $\hat\varphi_*\colon H^1_\fl(K,E'[2])\to H^1(K,\ker\varphi)$, compare the lift
with local Kummer images of $[2]$, and pair the differences, which lie in $H^1_\fl(K_v,\ker\hat\varphi)$, with
$\chi_v$ via $e_\varphi$. The global lift exists if and only if $(\chi',b)=0$: by Proposition~\ref{prop:B}(ii)(a) and
surjectivity of $\Phi$, $\chi'$ lifts if and only if there is $[\beta]$ with $\sigma[\beta]/[\beta]=[b]^{\chi'(\sigma)}$,
i.e.\ if and only if the cocycle $\sigma\mapsto[b]^{\chi'(\sigma)}$ is a coboundary in
$Z^1(G_K,\Ks{}^\times/\Ks{}^{\times2})$. The connecting map $\partial$ of
$1\to\Ks{}^\times\xrightarrow2\Ks{}^\times\to\Ks{}^\times/\Ks{}^{\times2}\to1$ is an isomorphism
$H^1(G_K,\Ks{}^\times/\Ks{}^{\times2})\simeq\Br(K)[2]$ (Hilbert's Theorem~90), and lifting $[b]^{\chi'(\sigma)}$ to
$b^{\chi'(\sigma)}$ one finds $d(b^{\chi'})=(b^{\chi'(\sigma)\chi'(\tau)})^2$, so $\partial$ maps the class to $(\chi',b)$
by \eqref{eq:cycliccocycle}. The class $(\chi',b)$ is locally trivial, hence trivial by Fact~\ref{fact:ABHN} --- this is
the proof of Lemma~\ref{lem:norm}.
\end{remark}

\begin{proofof}{Theorem~\ref{thm:W}}
Part (a) is Proposition~\ref{prop:B}. Parts (b) and (c) are
Corollary~\ref{cor:Aprime}, together with Theorem~\ref{thm:A}(a) for the explicit description of $g_v$ and
Lemma~\ref{lem:cup} for the equality $\chi_v\cup g_v=(\chi_v,g_v)$.
\end{proofof}

\subsection{Comparison with Kramer's pairing}\label{subsec:kramer}
We recall Kramer's definition \cite[\S1, \S3]{Kr} for the curve $A=E'$, i.e.\ for $a_1=1$ and $a_6=b$ in his notation;
his maps $\alpha$ and $\beta$ are $\delta_\varphi$ and $\kappa$, and his first Selmer group $M_\alpha^{(1)}$ is
$\Sel^\varphi(E/K)$. For a field $F\supseteq K$ let $\mathbb E_F=F[r]$ with $r^2=b$, and for $\zeta=f+gr\in\mathbb
E_F^\times$ ($f,g\in F$) put
\[
\theta_{\mathrm{Kr}}(\zeta)=\bigl(\zeta^2,\ \wp(bg^2/\zeta^2)\bigr)\in F^\times\times bF^2,\qquad\zeta^2=f^2+bg^2 .
\]
This is a homomorphism, and $H_{\mathrm{Kr}}(F)$ denotes its cokernel (Kramer's $\theta$ and $H$; the subscript
distinguishes them from our notation); the class of $(l,m)$ is written $((l,m))$. For
$R\in E'(F)$ with $x(R)\neq0$ put $\gamma_{\mathrm{Kr}}(R)=((x(R),\,b/x(R)^2))$ \cite[Prop.~1.4]{Kr}. Let $\chi'\in\Sel^\varphi(E/K)$.
Choose $c\in K$ with $\chi'=[bc^2]$ \cite[Lemma~3.2]{Kr}, for every place $v$ a point $R'_v\in E'(K_v)$ with
$\delta_\varphi(R'_v)=\chi'_v$, and $l_v\in K_v^\times$ with $\gamma_{\mathrm{Kr}}(R'_v)=((l_v,\,bc^2))$ in $H_{\mathrm{Kr}}(K_v)$; such $l_v$ exist
and are determined up to $\{1,b\}K_v^{\times2}$ \cite[Prop.~1.3(c)]{Kr}. Then
\[
U_\alpha(\chi,\chi')=\sum_v\inv_v(\chi_v,l_v)\qquad(\chi\in\Sel^\varphi(E/K)),
\]
and $U_\alpha$ is a well-defined alternating bilinear form on $\Sel^\varphi(E/K)$ \cite[Prop.~3.3]{Kr}.

\begin{proposition}\label{prop:kramer}
Let $\chi'\in\Sel^\varphi(E/K)$, $\chi'\neq0$, and let $c\in K^\times$ with $\chi'=[bc^2]$. Put $t_0=1/c$, let
$y\in\Ks$ with $y^2+y=b/t_0^2$, and $\beta_0=t_0y$. For every place $v$ let $R'_v\in E'(K_v)$ with
$\delta_\varphi(R'_v)=\chi'_v$ and $X_v:=x(R'_v)\neq0$, and let $h_v\in K_v^\times$ with
$h_v^2+h_v=b/X_v^2+b/t_0^2$.
\begin{enumerate}[label=\textup{(\roman*)},leftmargin=*]
\item $\beta_0$ satisfies the conditions of Lemma~\ref{lem:norm}, and $\gamma_{\mathrm{Kr}}(R'_v)=((X_vh_v,\,bc^2))$; thus one may
take $l_v=X_vh_v$ in Kramer's definition.
\item There is $P'_v\in E(\Kvs)$ with $\varphi(P'_v)=R'_v$ and $x(P'_v)=X_v(y+h_v)$, and for this choice the class of
Theorem~\ref{thm:A}\textup{(a)} is $g_v=[\,b\,t_0\,X_vh_v\,]=[b\,t_0]\cdot[l_v]$.
\end{enumerate}
Consequently $U_\alpha(\chi,\chi')=\sum_v\inv_v(\chi_v,g_v)=\langle\iota(\chi),\iota(\chi')\rangle$ for all
$\chi,\chi'\in\Sel^\varphi(E/K)$.
\end{proposition}

\begin{proof}
Such $R'_v$ exist because the fibres of $\delta_\varphi$ on $E'(K_v)$ are cosets of the infinite group
$\varphi(E(K_v))$, and $h_v$ exists because $[b/X_v^2]=\delta_\varphi(R'_v)=\chi'_v=[b/t_0^2]$ in $K_v/\wp K_v$
(Lemma~\ref{lem:deltaphi}); replacing $h_v$ by $h_v+1$ if necessary, $h_v\neq0$. (If $X_v=t_0$, this means
$h_v=1$; the arguments below apply, with $d=0$ and $\zeta=1$.)

(i) As $\Nm_{L'/K}(t_0y)=t_0^2\,y(y+1)=b$ and $y\notin K$, the element $\beta_0$ is as in Lemma~\ref{lem:norm}, with
$\beta_0+b/\beta_0=t_0$ (Remark~\ref{rem:BMint}(2)). Put $d=1/X_v+1/t_0$, so that $h_v^2+h_v=bd^2$, and
$\zeta=h_v+dr\in\mathbb E_{K_v}$. Then $\zeta^2=h_v^2+bd^2=h_v\neq0$ and $bd^2/\zeta^2=h_v+1$, so
$\theta_{\mathrm{Kr}}(\zeta)=(h_v,\wp(h_v))$. Hence, in $H_{\mathrm{Kr}}(K_v)$,
\[
\gamma_{\mathrm{Kr}}(R'_v)=((X_v,\,b/X_v^2))=((X_vh_v,\,b/X_v^2+\wp(h_v)))=((X_vh_v,\,b/t_0^2)).
\]

(ii) We have $\wp(y+h_v)=b/X_v^2$, so $X_v(y+h_v)$ is a root of $U^2+X_vU+b$, i.e.\ the $x$-coordinate of a point
$P'_v$ of $\varphi^{-1}(R'_v)$ (proof of Lemma~\ref{lem:deltaphi}); it satisfies the condition of
Theorem~\ref{thm:A}(a). If $\chi'_v\neq0$, we apply Lemma~\ref{lem:trace}(b) with $\alpha'=b/t_0^2$, $A=0$, $B=t_0$
and $h=h_v$; here $\gamma_v=t_0y/(X_v(y+h_v))\neq1$, because $1,y$ are linearly independent over $K_v$ and
$X_vh_v\neq0$. This gives $g_v=[bX_vt_0h_v]$. If $\chi'_v=0$, then $y\in K_v$, $P'_v\in E(K_v)$, and
$g_v=[\beta_0]\kappa(P'_v)^{-1}=[t_0y\,X_v(y+h_v)]$. From $y^2+y=bc^2$ and $h_v^2+h_v=bd^2$ we get
\[
y\,h_v\,(y+h_v)=y^2h_v+yh_v^2=(y+bc^2)h_v+y(h_v+bd^2)=b\,(c^2h_v+d^2y)=b\,(ch_v+dy)^2,
\]
the last step because $(ch_v+dy)^2=c^2(h_v+bd^2)+d^2(y+bc^2)=c^2h_v+d^2y$. The left-hand side is non-zero, so
$[y(y+h_v)]=[bh_v]$ and again $g_v=[b\,t_0X_vh_v]$.

For the last assertion, we may assume $\chi'\neq0$: for $\chi'=0$ the left-hand side vanishes by
Corollary~\ref{cor:bilin}, and so does $U_\alpha(\chi,0)$, by bilinearity \cite[Prop.~3.3]{Kr} or directly from the
definition with $c=0$, $R'_v=O$; and both sums are independent of the choices made (Theorem~\ref{thm:A}(b), \cite[Prop.~3.3]{Kr}). By (i)
and (ii),
\[
\sum_v\inv_v(\chi_v,g_v)=\sum_v\inv_v(\chi_v,l_v)+\sum_v\inv_v(\chi_v,bt_0)=U_\alpha(\chi,\chi'),
\]
since $bt_0\in K^\times$ (Fact~\ref{fact:ABHN}). The first sum is $\langle\iota(\chi),\iota(\chi')\rangle$ by
Theorem~\ref{thm:A}(b).
\end{proof}

This proves Corollary~\ref{cor:kramer}. By \cite[Thm.~3.1]{Kr} the radical of $U_\alpha$ is Kramer's second Selmer
group $M_\alpha^{(2)}$, the image in $K/\wp K$ of the group of elements of $H_{\mathrm{Kr}}(K)$ which lie in $\gamma_{\mathrm{Kr}}(E'(K_v))$ for
all $v$; comparing with Theorem~\ref{thm:radical}(a) we obtain $M_\alpha^{(2)}=\hat\varphi_*\Sel^2(E'/K)$, without
comparing $H_{\mathrm{Kr}}(K)$ with $H^1_\fl(K,E'[2])$.

\subsection{An example}\label{subsec:exW}
We compute the entry $\langle\iota(\chi_1),\iota(\chi_{\varrho^2S})\rangle$ of Theorem~\ref{thm:nonsplit}(b) with
Theorem~\ref{thm:W}. Let $K=\F_8(S)$, $b=S^4+S^3+\varrho$ and the notation be as in \S\ref{subsec:nonsplit}; here $b$
is not a square, so $E'[2]$ is a non-split extension of $\Z/2$ by $\mu_2$. Take $\chi=\chi_1$, $\chi'=\chi_{\varrho^2S}$,
$L'=K(y)$ with $y^2+y=\varrho^2S$, and the norm element
\[
\beta_0=A+By,\qquad A=S^2+\varrho^4S+1,\qquad B=\varrho^3(S+1),\qquad \Nm_{L'/K}\beta_0=b,
\]
of \S\ref{subsec:nonsplit}. By Theorem~\ref{thm:W}(a), $(\chi',[\beta_0])$ defines a class
$\xi\in H^1_\fl(K,E'[2])$ lifting $\chi'$; since $b$ is not a square, $[\beta_0]$ is not a class in
$K^\times/K^{\times2}$, only in $\Ks{}^\times/\Ks{}^{\times2}$.

\emph{Places other than $\infty$.} The characters $\chi,\chi'$ are unramified at all finite places; $b$ is a unit
away from $v_b$; and $\beta_0$, $\bar\beta_0$ are integral over $\F_8[S]$ with product $b$, so they are units above every
finite place $v\neq v_b$. By the argument of \S\ref{subsec:A6} the term at such a place vanishes. At $v_b$ the class
$\chi_{v_b}$ is trivial (proof of Theorem~\ref{thm:nonsplit}(a)), so the term $\inv_{v_b}(\chi_{v_b}\cup g_{v_b})$
vanishes too. Hence
\[
\langle\iota(\chi_1),\iota(\chi_{\varrho^2S})\rangle=\inv_\infty\bigl(\chi_{1,\infty}\cup g_\infty\bigr).
\]

\emph{The place $\infty$.} Let $t=1/S$. Take $R'_\infty=R_{\varrho^2}$, the point with $X:=x(R'_\infty)=\varrho^2S$ of
\S\ref{subsec:nonsplit}(a), so $\delta_\varphi(R'_\infty)=\chi'_\infty$, and use Lemma~\ref{lem:trace}(b) with
$\alpha'=\varrho^2S$. Since $\varrho^{-4}=\varrho^3$,
\[
\frac b{X^2}+\alpha'=\varrho^3(S^2+S)+\varrho^4t^2+\varrho^2S=\varrho^3S^2+\varrho^5S+\varrho^4t^2,
\]
and as $\wp(\varrho^5S)=\varrho^{10}S^2+\varrho^5S=\varrho^3S^2+\varrho^5S$, the element $h=\varrho^5S+h_0$ with
$h_0\in t\F_8[[t]]$, $h_0^2+h_0=\varrho^4t^2$, satisfies $h^2+h=b/X^2+\alpha'$. Then
\[
A+Bh=S^2+\varrho^4S+1+\varrho^3(S+1)(\varrho^5S+h_0)=\varrho^3S^2+\varrho^2S+O(1),
\]
using $\varrho^8=\varrho$, $1+\varrho=\varrho^3$ and $\varrho^4+\varrho=\varrho^2$. Choosing $P'_\infty$ with
$x(P'_\infty)=X(y+h)$ as in Lemma~\ref{lem:trace}(b), we have $\beta_0\neq x(P'_\infty)$ because $B\neq X$, so
$\gamma_\infty\neq1$, and Lemma~\ref{lem:trace}(b) gives
\[
g_\infty=[\,b\,X\,(A+Bh)\,],\qquad v_\infty(g_\infty)=-4-1-2=-7 .
\]
In the language of Theorem~\ref{thm:W}(b): the global lift $\xi$ and the local Kummer class $\delta_2(R'_\infty)$ have
the same image $\chi'_\infty$ in $H^1(K_\infty,\Z/2)$, and their difference is $j_*(g_\infty)$, a $\mu_2$-class of odd
valuation. The character $\chi_1$ is unramified at $\infty$ and non-trivial on the residue field $\F_8$
($\Tr_{\F_8/\F_2}(1)=1$), so by Fact~\ref{fact:schmid} and Lemma~\ref{lem:cup}
\[
\inv_\infty\bigl(\chi_{1,\infty}\cup g_\infty\bigr)=\Tr_{\F_8/\F_2}\Res_\infty\Bigl(1\cdot\frac{dg_\infty}{g_\infty}\Bigr)
=v_\infty(g_\infty)\cdot\Tr_{\F_8/\F_2}(1)=1\in\F_2 .
\]
So $\langle\iota(\chi_1),\iota(\chi_{\varrho^2S})\rangle=\frac12$, in agreement with Theorem~\ref{thm:nonsplit}(b). In
this example the local terms of Theorems~\ref{thm:BM} and~\ref{thm:W} agree place by place, although
$\chi'_\infty\neq0$ (compare Remark~\ref{rem:BMcompare}). The same computation for $\chi'=\chi_{\varrho S}$ (with
$X=\varrho S$ and the second norm element of \S\ref{subsec:nonsplit}) gives $A+Bh$ of degree $1$, hence
$v_\infty(g_\infty)=-6$ and $\langle\iota(\chi_1),\iota(\chi_{\varrho S})\rangle=0$, again as in
Theorem~\ref{thm:nonsplit}(b).

\section{The radical}\label{sec:radical}

In this section we prove Theorem~\ref{thm:radical}. Throughout, $K$, $E$, $E'$, $\varphi$ and $\hat\varphi$ are as in
\S\ref{sec:main}, $D_v=\delta_\varphi(E'(K_v))$, and $\kappa=\kappa_{K_v}=\delta_{\hat\varphi}$ is the Kummer map of
$\hat\varphi$ (Proposition~\ref{prop:kummerF}). The local results of \S\ref{subsec:localimages} are due to Kramer
\cite[\S2]{Kr}, the global duality of \S\ref{subsec:global} is a restricted-product form of \cite[Lemmas~3.6,
3.7]{Kr}, and the proof in \S\ref{subsec:proofradical} follows his proof of \cite[Thm.~3.1]{Kr}.

\subsection{Local duality and the local Kummer images}\label{subsec:localimages}
For every place $v$ the pairing
\begin{equation}\label{eq:localpairing}
H^1(K_v,\Z/2)\times K_v^\times/K_v^{\times2}\to\F_2,\qquad(\chi,u)\mapsto\inv_v(\chi,u),
\end{equation}
identifies the discrete group $H^1(K_v,\Z/2)=K_v/\wp K_v$ with the Pontryagin dual of the profinite group
$K_v^\times/K_v^{\times2}$, and vice versa. Indeed $\inv_v(\chi,u)=\chi(\mathrm{rec}_v(u))$ for the reciprocity map
$\mathrm{rec}_v$, which induces an isomorphism of $K_v^\times/K_v^{\times2}$ onto the maximal quotient of exponent $2$ of
$G_v^{\mathrm{ab}}$ (local class field theory \cite[XIII, XIV, \S1]{Se}; see also \cite[III, Thm.~6.10]{Mi}). For a
subgroup $L\subseteq H^1(K_v,\Z/2)$ let $L^\perp\subseteq K_v^\times/K_v^{\times2}$ be its annihilator. The unramified
classes are the classes of the elements of $\mathcal O_v$, and $H^1_\nr(K_v,\Z/2)^\perp=U_v$ by Fact~\ref{fact:schmid}.

\begin{lemma}[{cf.\ \cite[Prop.~2.4]{Kr}}]\label{lem:Dv}
\begin{enumerate}[label=\textup{(\roman*)},leftmargin=*]
\item $D_v\simeq E'(K_v)/\varphi E(K_v)$, and this group is finite.
\item If $a_2\in\mathcal O_v$ and $b\in\mathcal O_v^\times$, then $D_v=H^1_\nr(K_v,\Z/2)$.
\end{enumerate}
In particular $D_v=H^1_\nr(K_v,\Z/2)$ for almost all $v$.
\end{lemma}

\begin{proof}
(i) The isomorphism comes from the Kummer sequence of $\varphi$. As $\varphi$ is \'etale, $\varphi(E(K_v))$ is an open
subgroup of the compact group $E'(K_v)$ (implicit function theorem), hence of finite index.

(ii) Let $R\in E'(K_v)$, $R\neq O$. If $v(x(R))\geq0$, then $\delta_\varphi(R)=[x(R)+a_2]$ with
$x(R)+a_2\in\mathcal O_v$ is unramified. If $v(x(R))<0$, then $b/x(R)^2\in\pi_v\mathcal O_v\subseteq\wp(K_v)$, so
$\delta_\varphi(R)=[b/x(R)^2]=0$ (Lemma~\ref{lem:deltaphi}). Hence $D_v\subseteq H^1_\nr(K_v,\Z/2)$. Conversely, $E$
has good reduction at $v$ (its discriminant $b^2$ is a unit), so for unramified $\chi$ the class $\iota(\chi)_v$ is
inflated from $H^1(\Gal(K_v^{\nr}/K_v),E(K_v^{\nr}))=0$ \cite[I, Prop.~3.8]{Mi}; that is, $\chi\in D_v$.
\end{proof}

\begin{lemma}[{cf.\ \cite[Thm.~2.1]{Kr}}]\label{lem:annihilator}
For every place $v$, $D_v^\perp=\kappa(E(K_v))$.
\end{lemma}

\begin{proof}
The inclusion $\kappa(E(K_v))\subseteq D_v^\perp$ is Lemma~\ref{lem:iso}, since every $\chi\in D_v$ is of the form
$\delta_\varphi(\varphi P)$ with $P\in E(\Kvs)$. Both groups have finite index in $H^1_\fl(K_v,\mu_2)$, and we show that
the indices agree. By the Kummer sequence of $\hat\varphi$ and Proposition~\ref{prop:kummerF}, the index of
$\kappa(E(K_v))$ is the order of the kernel of $\hat\varphi_*\colon H^1(K_v,E')\to H^1(K_v,E)$, which by local duality for
abelian varieties \cite[III, Thm.~7.8]{Mi} (valid for $p$-primary parts) is dual to
$\operatorname{coker}(\varphi\colon E(K_v)\to E'(K_v))\simeq D_v$. On the other hand the local pairing
$H^1(K_v,\Z/2)\times H^1_\fl(K_v,\mu_2)\to\frac12\Z/\Z$ is perfect (\S\ref{subsec:localimages}), so
$H^1_\fl(K_v,\mu_2)/D_v^\perp\simeq\Hom(D_v,\frac12\Z/\Z)$ has order $|D_v|$. Hence the two indices are both $|D_v|$,
and the inclusion is an equality.
\end{proof}

This is a special case of \cite[Prop.~B.1]{Ces2}; Kramer obtains it by computing the index of $\kappa(E(K_v))$
directly \cite[Props.~2.2, 2.3]{Kr}.

\subsection{Global duality}\label{subsec:global}
Let
\[
\mathcal P=\textstyle\prod'_v\bigl(H^1(K_v,\Z/2);\,H^1_\nr(K_v,\Z/2)\bigr),\qquad
\mathcal P^\vee=\prod'_v\bigl(K_v^\times/K_v^{\times2};\,U_v\bigr)
\]
be the restricted products. They are locally compact, with compact open subgroups $\prod_vH^1_\nr(K_v,\Z/2)$ and
$\prod_vU_v$, and since the local pairings are perfect and $H^1_\nr(K_v,\Z/2)^\perp=U_v$, the pairing
$\langle(\chi_v),(u_v)\rangle=\sum_v\inv_v(\chi_v,u_v)$, a finite sum, identifies each of them with the Pontryagin dual
of the other. For a subgroup $Y$ of $\mathcal P$ or of $\mathcal P^\vee$ we write $Y^\perp$ for its annihilator in the other group; for
closed subgroups $Y_1,Y_2$ one has $(Y_1\cap Y_2)^\perp=\overline{Y_1^\perp+Y_2^\perp}$, as both sides have annihilator
$Y_1\cap Y_2$. Let $\mathcal A\subseteq \mathcal P$ and $\mathcal A'\subseteq \mathcal P^\vee$ be the images of $H^1(K,\Z/2)$ and of
$K^\times/K^{\times2}$; both maps are injective \cite[Lemma~3.6]{Kr}.

\begin{lemma}[{cf.\ \cite[Lemmas~3.6, 3.7]{Kr}}]\label{lem:globalduality}
\begin{enumerate}[label=\textup{(\roman*)},leftmargin=*]
\item $\mathcal A$ and $\mathcal A'$ are discrete, and $\mathcal A^\perp=\mathcal A'$.
\item Let $\mathcal W=\prod_vL_v\subseteq \mathcal P$ for finite subgroups $L_v\subseteq H^1(K_v,\Z/2)$ with
$L_v=H^1_\nr(K_v,\Z/2)$ for almost all $v$. Then $\mathcal W$ is compact and open, $\mathcal W^\perp=\prod_vL_v^\perp$
is compact and open, and $(\mathcal A\cap\mathcal W)^\perp=\mathcal A'+\mathcal W^\perp$.
\end{enumerate}
\end{lemma}

\begin{proof}
(i) The group $\mathcal A\cap\prod_vH^1_\nr(K_v,\Z/2)$ consists of the everywhere unramified characters and is dual to
$\Pic(C)/2\Pic(C)$; the group $\mathcal A'\cap\prod_vU_v$ consists of the classes of the $u\in K^\times$ with
$\dv u\in2\Div(C)$ and is isomorphic to $\Pic(C)[2]$ by $u\mapsto[\frac12\dv u]$ (constants are squares). Both are
finite \cite[Chap.~XIII]{We}, so $\mathcal A$ and $\mathcal A'$ are discrete. By reciprocity (Fact~\ref{fact:ABHN}),
$\mathcal A'\subseteq\mathcal A^\perp$. Conversely, let $g=(g_v)\in J_K$ with $\sum_v\inv_v(\chi_v,g_v)=0$ for all $\chi\in H^1(K,\Z/2)$, and let
$C_K=J_K/K^\times$. The group $C_K$ is the product of $\Z$ and a compact group, so $C_K/C_K^2$ is a compact Hausdorff
group of exponent $2$, whose continuous characters separate points. By the existence theorem of global class field
theory \cite[XIII, \S9]{We}, every continuous character of order $2$ of $C_K$ is the composite of the reciprocity map
with a character $\chi$ of $G_K$, and its value at the class of $g$ is $\sum_v\inv_v(\chi_v,g_v)=0$ by local--global
compatibility. So $g\in K^\times J_K^2$, i.e.\ $\mathcal A^\perp\subseteq\mathcal A'$. For the components at a
finite set of places containing generators of the class group this is \cite[Lemma~3.7]{Kr}, and it is a special case
of the Poitou--Tate duality of \cite[Prop.~4.10(c)]{Ces}.

(ii) If $L_v=H^1_\nr(K_v,\Z/2)$ for $v\notin S$, $S$ finite, then $\mathcal W^\perp=\prod_{v\in S}L_v^\perp\times
\prod_{v\notin S}U_v$, and $L_v^\perp$ is closed of finite index $|L_v|$, hence open. So $\mathcal W$ and
$\mathcal W^\perp$ are compact and open, and $\mathcal A$ is closed by (i). Therefore
$(\mathcal A\cap\mathcal W)^\perp=\overline{\mathcal A'+\mathcal W^\perp}=\mathcal A'+\mathcal W^\perp$, the last group
being open, hence closed.
\end{proof}

\begin{lemma}[{cf.\ \cite[Prop.~3.8]{Kr}}]\label{lem:ratio}
For a family $\mathcal L=(L_v)$ as in Lemma~\ref{lem:globalduality}\textup{(ii)} put
$\Sel_{\mathcal L}=\{\chi\in H^1(K,\Z/2):\chi_v\in L_v\ \forall v\}$ and
$\Sel_{\mathcal L^\perp}=\{u\in K^\times/K^{\times2}:u_v\in L_v^\perp\ \forall v\}$. These groups are finite, and
\[
\frac{|\Sel_{\mathcal L}|}{|\Sel_{\mathcal L^\perp}|}=2\prod_v\frac{|L_v|}{2} .
\]
\end{lemma}

\begin{proof}
Write $\mathcal W=\prod_vL_v$. Then $\Sel_{\mathcal L}\simeq\mathcal A\cap\mathcal W$ and
$\Sel_{\mathcal L^\perp}\simeq\mathcal A'\cap\mathcal W^\perp$ are finite, being discrete and compact. Write
$r(\mathcal L)$ for the left-hand side of the formula.

Let $\mathcal L'=(L'_v)$ be a second family with $L_v\subseteq L'_v$ for all $v$, and $\mathcal W'=\prod_vL'_v$. The
pairing induces a perfect pairing between the finite groups $\mathcal W'/\mathcal W$ and
$\mathcal W^\perp/\mathcal W'^\perp$. The image $X$ of $\mathcal A\cap\mathcal W'$ in $\mathcal W'/\mathcal W$ has order
$|\Sel_{\mathcal L'}|/|\Sel_{\mathcal L}|$, and the image $Y$ of $\mathcal A'\cap\mathcal W^\perp$ in
$\mathcal W^\perp/\mathcal W'^\perp$ has order $|\Sel_{\mathcal L^\perp}|/|\Sel_{\mathcal L'^\perp}|$. The preimage of $X$
in $\mathcal W'$ is $(\mathcal A+\mathcal W)\cap\mathcal W'$, and $\mathcal A+\mathcal W$ is open, hence closed, with
annihilator $\mathcal A'\cap\mathcal W^\perp$; so the annihilator of the preimage is
$\overline{(\mathcal A'\cap\mathcal W^\perp)+\mathcal W'^\perp}=(\mathcal A'\cap\mathcal W^\perp)+\mathcal W'^\perp$, whose
image in $\mathcal W^\perp/\mathcal W'^\perp$ is $Y$. Thus $Y$ is the annihilator of $X$, so $|X|\,|Y|=[\mathcal W':\mathcal W]$,
i.e.\ $r(\mathcal L')=r(\mathcal L)\prod_v|L'_v/L_v|$.

Let $\mathcal U=(H^1_\nr(K_v,\Z/2))_v$. Comparing $\mathcal L$ and $\mathcal U$ with the family
$(L_v\cap H^1_\nr(K_v,\Z/2))_v$ gives $r(\mathcal L)=r(\mathcal U)\prod_v|L_v|/2$. Finally $r(\mathcal U)=2$: by the
proof of Lemma~\ref{lem:globalduality}(i), $|\Sel_{\mathcal U}|=|\Hom(\Pic(C),\Z/2)|=2\,|\Pic^0(C)/2|$, because
$\Pic(C)\simeq\Z\oplus\Pic^0(C)$ with $\Pic^0(C)$ finite, and $|\Sel_{\mathcal U^\perp}|=|\Pic(C)[2]|=|\Pic^0(C)[2]|$.
\end{proof}

\subsection{Proof of Theorem~\ref{thm:radical}}\label{subsec:proofradical}
Let $\delta_2\colon E'(F)\to H^1_\fl(F,E'[2])$ be the Kummer map of $[2]_{E'}$,
$\iota_2\colon H^1_\fl(K,E'[2])\to H^1(K,E')$ the natural map, and $j\colon\mu_2=\ker\hat\varphi\hookrightarrow E'[2]$.
Since $\varphi\circ\hat\varphi=[2]_{E'}$, the triples $(\hat\varphi,\hat\varphi,\mathrm{id})$ and
$(j,\mathrm{id},\varphi)$ are morphisms of short exact sequences of fppf sheaves
\begin{gather*}
\bigl(0\to E'[2]\to E'\xrightarrow{2}E'\to0\bigr)\longrightarrow\bigl(0\to\ker\varphi\to E\xrightarrow{\varphi}E'\to0\bigr),\\
\bigl(0\to\mu_2\to E'\xrightarrow{\hat\varphi}E\to0\bigr)\longrightarrow\bigl(0\to E'[2]\to E'\xrightarrow{2}E'\to0\bigr),
\end{gather*}
and they give
\begin{equation}\label{eq:compat}
\hat\varphi_*\circ\delta_2=\delta_\varphi,\qquad \iota\circ\hat\varphi_*=\hat\varphi_*\circ\iota_2,\qquad
j_*\circ\kappa=\delta_2\circ\varphi .
\end{equation}
In particular $\hat\varphi_*\Sel^2(E'/K)\subseteq\Sel^\varphi(E/K)$ and
$\delta_\varphi(E'(K))=\hat\varphi_*\delta_2(E'(K))\subseteq\hat\varphi_*\Sel^2(E'/K)$. The long exact sequence of
$0\to\mu_2\xrightarrow{j}E'[2]\xrightarrow{\hat\varphi}\Z/2\to0$ (Proposition~\ref{prop:B}(i)) shows, for $F=K$ and
$F=K_v$, that the kernel of $j_*$ on $H^1_\fl(F,\mu_2)$ is generated by $[b]=\kappa(T)$, and that the kernel of
$\hat\varphi_*$ on $H^1_\fl(F,E'[2])$ is $j_*H^1_\fl(F,\mu_2)$.

\begin{lemma}[{cf.\ \cite[\S1]{Kr}}]\label{lem:exact}
The sequence
\[
0\to\frac{E(K)[\varphi]}{\hat\varphi\bigl(E'(K)[2]\bigr)}\xrightarrow{\ \kappa\ }\Sel^{\hat\varphi}(E'/K)
\xrightarrow{\ j_*\ }\Sel^2(E'/K)\xrightarrow{\ \hat\varphi_*\ }\Sel^\varphi(E/K)\xrightarrow{\ \iota\ }
\frac{\Sha(E/K)[\varphi]}{\hat\varphi\bigl(\Sha(E'/K)[2]\bigr)}\to0
\]
is exact, and its first term has dimension $\varepsilon_b$.
\end{lemma}

\begin{proof}
By \eqref{eq:compat} the maps are defined, and $\hat\varphi_*\circ j_*=0$. \emph{The first term.} Here
$E(K)[\varphi]=\{O,T\}$. If $b=b_0^2$ with $b_0\in K$, then $T'=(0,b_0)\in E'(K)[2]$ and $\hat\varphi(T')=T$, so the first
term is $0$; otherwise $E'(K)[2]=0$ (Remark~\ref{rem:model}) and the first term has order $2$. In both cases $\kappa$
maps it isomorphically onto $\langle[b]\rangle=\ker j_*\cap\Sel^{\hat\varphi}(E'/K)$. \emph{Exactness at
$\Sel^2(E'/K)$.} Let $y\in\Sel^2(E'/K)$ with $\hat\varphi_*y=0$, so $y=j_*u$ with $u\in K^\times/K^{\times2}$. For each
$v$, $y_v=\delta_2(R_v)$ with $R_v\in E'(K_v)$, and $\delta_\varphi(R_v)=\hat\varphi_*y_v=0$, so $R_v=\varphi(P_v)$ with
$P_v\in E(K_v)$. Then $j_*(u_v)=\delta_2(\varphi P_v)=j_*\kappa(P_v)$, so $u_v\in\kappa(P_v)\langle[b]\rangle\subseteq
\kappa(E(K_v))$. Thus $u\in\Sel^{\hat\varphi}(E'/K)$. \emph{Exactness at $\Sel^\varphi(E/K)$, and surjectivity.} The map
$\iota\colon\Sel^\varphi(E/K)\to\Sha(E/K)[\varphi]$ is surjective with kernel
$\delta_\varphi(E'(K))\subseteq\hat\varphi_*\Sel^2(E'/K)$, and by \eqref{eq:compat} it maps $\hat\varphi_*\Sel^2(E'/K)$ onto
$\hat\varphi_*\iota_2\Sel^2(E'/K)=\hat\varphi\bigl(\Sha(E'/K)[2]\bigr)$. Hence the kernel of the last map is
$\hat\varphi_*\Sel^2(E'/K)$.
\end{proof}

We now prove Theorem~\ref{thm:radical}.

\begin{proof}
(a) Suppose $\chi'=\hat\varphi_*y$ with $y\in\Sel^2(E'/K)$. Then $\iota(\chi')=\hat\varphi_*\iota_2(y)$ by
\eqref{eq:compat}, with $\iota_2(y)\in\Sha(E'/K)$, and Proposition~\ref{prop:F} for $\rho=\varphi$, $\rho^\vee=\hat\varphi$
gives $\langle\iota(\chi),\iota(\chi')\rangle=\langle\varphi_*\iota(\chi),\iota_2(y)\rangle_{E'}=0$ for every $\chi$,
because $\varphi_*\circ\iota=0$.

Conversely, let $\chi'\in\Sel^\varphi(E/K)$ be orthogonal to $\Sel^\varphi(E/K)$. Fix $\beta_0$ as in
Lemma~\ref{lem:norm}, for every place $v$ a point $P'_v$ as in Theorem~\ref{thm:A}(a), and put $R'_v=\varphi(P'_v)$.
Let $g_v\in K_v^\times/K_v^{\times2}$ be the classes of Theorem~\ref{thm:A}(a); they depend only on $\chi'$ and on these
choices, and by Theorem~\ref{thm:A}(b)
\begin{equation}\label{eq:gperp}
\sum_v\inv_v(\chi_v,g_v)=\langle\iota(\chi),\iota(\chi')\rangle=0\qquad\text{for all }\chi\in\Sel^\varphi(E/K).
\end{equation}

\emph{The family $g=(g_v)$ lies in $\mathcal P^\vee$.} Let $\Sigma'$ be the finite set of places at which $\chi'$ is ramified,
$a_2\notin\mathcal O_v$, $b\notin\mathcal O_v^\times$, or $\beta_0$ or $b/\beta_0$ is not a unit, and let $v\notin\Sigma'$.
If $\chi'_v\neq0$, then $g_v=[\Nm w_0]^{-1}$ is the class of a norm from the unramified quadratic extension
$L'_v/K_v$ (\S\ref{subsec:A5}), so $g_v\in U_v$. If $\chi'_v=0$, then $P'_v\in E(K_v)$ and
$g_v=[\beta_0]\kappa(P'_v)^{-1}$ with $\beta_0\in\mathcal O_v^\times$ and $\kappa(P'_v)\in D_v^\perp=U_v$
(Lemmas~\ref{lem:annihilator} and~\ref{lem:Dv}(ii)); so again $g_v\in U_v$.

\emph{Global duality.} Let $\mathcal W=\prod_vD_v$; by Lemma~\ref{lem:Dv}, Lemma~\ref{lem:globalduality}(ii) applies to
it, and the image of $\Sel^\varphi(E/K)$ in $\mathcal P$ is $\mathcal A\cap\mathcal W$. By \eqref{eq:gperp},
$g\in(\mathcal A\cap\mathcal W)^\perp=\mathcal A'+\mathcal W^\perp$. So there are $u\in K^\times/K^{\times2}$ and
$d_v\in D_v^\perp$ with $g_v=u_vd_v$ for all $v$, and by Lemma~\ref{lem:annihilator} $d_v=\kappa(P_v)$ with
$P_v\in E(K_v)$.

\emph{The $2$-Selmer class.} Let $\xi=\Phi^{-1}(\chi',[\beta_0])\in H^1_\fl(K,E'[2])$ be the lift of $\chi'$ of
Proposition~\ref{prop:B}(ii)(d), and put $\xi'=\xi-j_*(u)$, where $u$ is viewed in $H^1_\fl(K,\mu_2)$. Since
$\hat\varphi\circ j=0$, $\hat\varphi_*\xi'=\chi'$. Writing the groups $H^1_\fl(K_v,\mu_2)$ additively, Corollary~\ref{cor:Aprime}
and \eqref{eq:compat} give, for every $v$,
\[
\xi'_v=\delta_2(R'_v)+j_*(g_v)-j_*(u_v)=\delta_2(R'_v)+j_*\kappa(P_v)=\delta_2\bigl(R'_v+\varphi(P_v)\bigr).
\]
Hence $\xi'\in\Sel^2(E'/K)$, and $\chi'=\hat\varphi_*\xi'\in\hat\varphi_*\Sel^2(E'/K)$.

(b) The map $\iota$ induces an isomorphism $\Sel^\varphi(E/K)/\hat\varphi_*\Sel^2(E'/K)\simeq
\Sha(E/K)[\varphi]/\hat\varphi\bigl(\Sha(E'/K)[2]\bigr)$ (Lemma~\ref{lem:exact}), compatible with the pairings, so the
first and second statements of (b) follow from (a); and if $\hat\varphi(y)\in\Sha(E/K)[\varphi]$ with
$y\in\Sha(E'/K)$, then $2y=\varphi\hat\varphi(y)=0$, so $\Sha(E/K)[\varphi]\cap\hat\varphi\,\Sha(E'/K)=
\hat\varphi\bigl(\Sha(E'/K)[2]\bigr)$. The induced form is alternating by Corollary~\ref{cor:bilin}, and it lives on a
finite-dimensional space, since $\Sel^\varphi(E/K)$ is finite (Lemma~\ref{lem:ratio}). A non-degenerate alternating form
over $\F_2$ exists only in even dimension.

(c) By Lemma~\ref{lem:exact},
\[
\dim\Sel^2(E'/K)=\dim\Sel^{\hat\varphi}(E'/K)-\varepsilon_b+\dim\Sel^\varphi(E/K)
-\dim\bigl(\Sha(E/K)[\varphi]/\hat\varphi\bigl(\Sha(E'/K)[2]\bigr)\bigr),
\]
and the last dimension is even by (b).
\end{proof}

\section{The parity of \texorpdfstring{$2$}{2}-Selmer ranks}\label{sec:parity}
We keep the notation of \S\ref{sec:radical}.

\subsection{The Selmer ratio}

\begin{proposition}\label{prop:ratio}
$\dim\Sel^\varphi(E/K)-\dim\Sel^{\hat\varphi}(E'/K)=1+\sum_v\bigl(\dim D_v-1\bigr)$, where almost all terms of the sum
vanish.
\end{proposition}

\begin{proof}
By Lemma~\ref{lem:Dv}, Lemma~\ref{lem:ratio} applies to $\mathcal L=(D_v)_v$. By definition
$\Sel_{\mathcal L}=\Sel^\varphi(E/K)$, and $\Sel_{\mathcal L^\perp}=\Sel^{\hat\varphi}(E'/K)$ by
Lemma~\ref{lem:annihilator} and Proposition~\ref{prop:kummerF}.
\end{proof}

This is the analogue, for the Verschiebung in characteristic $2$, of Cassels' formula for the quotient of the orders of the
Selmer groups of an isogeny and of its dual; compare the counting argument in the proof of \cite[Prop.~3.8]{Kr}.

\subsection{The size of the local images}

\begin{lemma}[{cf.\ \cite[Props.~2.2, 2.3]{Kr}}]\label{lem:Dsize}
Let $v$ be a place, let $\omega_v,\omega'_v$ be N\'eron differentials of $E,E'$ at $v$, and $\varphi^*\omega'_v=\vartheta_v\omega_v$.
Then $\vartheta_v\neq0$ and
\[
|D_v|=2\cdot\frac{c_v(E')}{c_v(E)}\cdot q_v^{\,v(\vartheta_v)} .
\]
Moreover $\sum_v\deg(v)\,v(\vartheta_v)=\deg\omega_E-\deg\omega_{E'}$.
\end{lemma}

\begin{proof}
For a homomorphism $\rho$ of abelian groups with finite kernel and cokernel let
$\mathfrak z(\rho)=|\operatorname{coker}\rho|/|\ker\rho|$ \cite[\S2]{TW}. By the snake lemma $\mathfrak z$ is
multiplicative along morphisms of short exact sequences, and $\mathfrak z(\rho)=|B|/|A|$ if $\rho\colon A\to B$ with
$A$, $B$ finite. As $E(K_v)[\varphi]=\{O,T\}$, Lemma~\ref{lem:Dv}(i) gives $|D_v|=2\,\mathfrak z(\varphi_v)$ for
$\varphi_v\colon E(K_v)\to E'(K_v)$.

Let $\mathcal E$, $\mathcal E'$ be the N\'eron models over $\mathcal O_v$, $E^0(K_v)\subseteq E(K_v)$ the subgroup of
points reducing into the identity component $\mathcal E^0_s$ of the special fibre, and $E^1(K_v)$ the kernel of
reduction; similarly for $E'$. By the N\'eron mapping property $\varphi$ extends to $\mathcal E\to\mathcal E'$, so
$\varphi_v$ respects these filtrations, and $\mathfrak z(\varphi_v)$ is the product of the values of $\mathfrak z$ on
the graded pieces.
\begin{itemize}
\item $E(K_v)/E^0(K_v)$ and $E'(K_v)/E'^0(K_v)$ have orders $c_v(E)$ and $c_v(E')$; this gives the factor $c_v(E')/c_v(E)$.
\item $E^0(K_v)/E^1(K_v)\simeq\mathcal E^0_s(k_v)$ by Hensel's lemma, and similarly for $E'$. The homomorphisms
between $\mathcal E^0_s$ and $\mathcal E'^0_s$ induced by $\varphi$ and $\hat\varphi$ compose, in either order, to
multiplication by $2$, which is non-zero on an elliptic curve and on a one-dimensional torus. There are no non-zero
homomorphisms between connected one-dimensional groups of different kinds (elliptic curve, torus, unipotent group).
Hence $\mathcal E^0_s$ and $\mathcal E'^0_s$ are either isogenous elliptic curves, or isogenous one-dimensional tori
(so both split or both non-split), or both unipotent, hence isomorphic to $\mathbb G_a$ since $k_v$ is perfect. In each
case they have the same number of $k_v$-points (for elliptic curves because isogenous curves over a finite field have
isomorphic Tate modules), and this graded piece contributes the factor $1$.
\item $E^1(K_v)=\hat E(\mathfrak m_v)$ for the formal group $\hat E$ of a minimal Weierstrass model, and the
successive quotients of the filtration by the $\hat E(\mathfrak m_v^N)$ are isomorphic to $k_v$
\cite[VII.2.2, IV.3.2]{Si}. On formal groups $\varphi$ is a power series $f(z)=\vartheta_vz+\sum_{i\geq2}c_iz^i$ with
$c_i\in\mathcal O_v$ \cite[IV.4.3]{Si}; the linear coefficient is $\vartheta_v$ because the parameters $z,z'$ come from
minimal models, whose invariant differentials are the N\'eron differentials and have the form $(1+\cdots)dz$ and
$(1+\cdots)dz'$, so that comparing linear terms in $\varphi^*\omega'_v=f'(z)(1+\cdots)dz=\vartheta_v\omega_v$ gives
$f'(0)=\vartheta_v$; and $\vartheta_v\neq0$ because $\varphi$ is separable \cite[II.4.2(c)]{Si}. Let
$N>v(\vartheta_v)$. For $z\in\mathfrak m_v^N$, $v(c_iz^i)\geq2v(z)>v(\vartheta_v)+v(z)$, so $f$ maps $\hat E(\mathfrak m_v^N)$
injectively into $\hat E'(\mathfrak m_v^{N+v(\vartheta_v)})$; it is onto by Hensel's lemma (for
$w\in\mathfrak m_v^{N+v(\vartheta_v)}$, $z\mapsto(w-\sum_{i\geq2}c_iz^i)/\vartheta_v$ is a contraction of $\mathfrak m_v^N$).
So $f\colon\hat E(\mathfrak m_v^N)\to\hat E'(\mathfrak m_v^N)$ has $\mathfrak z(f)=q_v^{v(\vartheta_v)}$, while
$E^1(K_v)/\hat E(\mathfrak m_v^N)$ and $E'^1(K_v)/\hat E'(\mathfrak m_v^N)$ both have order $q_v^{N-1}$.
\end{itemize}
This proves the formula for $|D_v|$; Kramer expresses the same index through the indices of the kernels of reduction
of an integral, not necessarily minimal, model \cite[Props.~2.2, 2.3 and the remark following them]{Kr}. For the second
statement let $\eta=dx/x$ and $\eta'=dX/X$ be the invariant
differentials of the models \eqref{eq:E} of $E$ and $E'$, where $X$ is the coordinate on $E'$. Since
$X\circ\varphi=x+b/x$ (Lemma~\ref{lem:isog}) and $d(x+b/x)=(1+b/x^2)\,dx$ in characteristic $2$, we have
$\varphi^*\eta'=\eta$. Write $\omega_v=e_v\eta$ and $\omega'_v=e'_v\eta'$ with $e_v,e'_v\in K_v^\times$; then
$\vartheta_v=e'_v/e_v$. The differential $\eta$ is a rational section of $\omega_E$ whose order at $v$ is $-v(e_v)$, so
$\deg\omega_E=-\sum_v\deg(v)\,v(e_v)$, and similarly for $E'$. Hence
$\sum_v\deg(v)\,v(\vartheta_v)=\deg\omega_E-\deg\omega_{E'}$.
\end{proof}

\subsection{Proof of Theorem~\ref{thm:parity}}
\begin{proof}
By Remark~\ref{rem:model} we may write $A\colon y^2+xy=x^3+a_2x^2+b$. Then $A=E'$, and by Remark~\ref{rem:FV} we may
identify $A^{(2)}$ with $E$ and $V$ with $\varphi$, so that $D_v(A)=D_v$, $c_v(A^{(2)})=c_v(E)$ and
$\omega_{A^{(2)}}\simeq\omega_E$. Moreover $\dim A(K)[2]=1-\varepsilon_b$ (Remark~\ref{rem:model}). The exact sequence
$0\to A(K)/2A(K)\to\Sel^2(A/K)\to\Sha(A/K)[2]\to0$ and the Lang--N\'eron theorem ($A(K)$ is finitely generated) give
$s_2(A/K)=\rk A(K)+\dim\Sha(A/K)[2]$. By Theorem~\ref{thm:radical}(c) and Proposition~\ref{prop:ratio},
\[
s_2(A/K)\equiv\dim\Sel^\varphi(E/K)-\dim\Sel^{\hat\varphi}(E'/K)+\varepsilon_b-(1-\varepsilon_b)
\equiv\sum_v\bigl(\dim D_v-1\bigr)\pmod2,
\]
and $\dim D_v=1$ for almost all $v$ by Lemma~\ref{lem:Dv}. By Lemma~\ref{lem:Dsize},
$\dim D_v-1=\ord_2\bigl(c_v(A)/c_v(A^{(2)})\bigr)+f_v\,v(\vartheta_v)$ with $f_v=n\deg(v)$, and summing over $v$ gives the
second congruence.
\end{proof}

\begin{remark}[the Hodge term]\label{rem:hodge}
For a Weierstrass model of $A$ with invariant differential $\omega$, the differential $a_1\omega$ does not depend on the
model ($a_1'=a_1/u$ and $\omega'=u\omega$ in characteristic $2$); it is the Hasse invariant, a non-zero global section
of $\omega_A$ whose order at $v$ is the valuation of $a_1$ on a minimal model at $v$. So
$\deg\omega_A=\sum_v\deg(v)\,v(a_{1,v})$ for every ordinary $A$, where $a_{1,v}$ is $a_1$ of a minimal model at $v$
(compare Propositions~\ref{prop:exnonsplit} and~\ref{prop:exsplit}). Suppose now that $A$ is semistable.
Then this order is $0$ at places of good ordinary and of multiplicative reduction, and at a place of good
supersingular reduction it is the integer $k=k(v)\geq1$ of Lemma~\ref{lem:good}; so
$\deg\omega_A=\sum_v\deg(v)\,k(v)$, the sum over the supersingular places. Moreover the Frobenius twist of a minimal model of $A$ at $v$ is a minimal model of $A^{(2)}$ (proofs of
Lemmas~\ref{lem:good} and~\ref{lem:mult}), so $\omega_{A^{(2)}}$ is the pull-back of $\omega_A$ by the absolute
Frobenius of the curve and $\deg\omega_{A^{(2)}}=2\deg\omega_A$. Hence
\[
n\bigl(\deg\omega_{A^{(2)}}-\deg\omega_A\bigr)=n\deg\omega_A=\sum_{v\ \mathrm{supersingular}}f_vk(v),
\]
which is the sum of the local terms $\dim D_v-1$ of Lemma~\ref{lem:good}, while the Tamagawa term of
Theorem~\ref{thm:parity} is the sum of those of Lemma~\ref{lem:mult}. In \S\ref{subsec:exparity}, $\deg\omega_A=1$
comes from the supersingular place $\infty$. In characteristic $p$ the same argument gives
$\deg\omega_{A^{(p)}}-\deg\omega_A=(p-1)\deg\omega_A$, which is even for odd $p$. Semistability is needed: at a place of
additive reduction the Frobenius twist of a minimal model need not be minimal (in Proposition~\ref{prop:exnonsplit},
$v_\infty(\Delta_{\min})$ is $8$ for $E'$ and $4$, not $16$, for $E\simeq E'^{(2)}$).
\end{remark}

\begin{corollary}[constant field extensions]\label{cor:tower}
Let $A/K$ be an ordinary elliptic curve, $\F_q$ the constant field of $K$, $q=2^n$, and $K_k=K\F_{q^k}$. There is a
constant $C$, depending on $A$ and $K$ but not on $k$, such that
\[
\dim\Sha(A^{(2)}/K_k)[\varphi]\geq nk\bigl(\deg\omega_{A^{(2)}}-\deg\omega_A\bigr)-C\qquad(k\geq1),
\]
where $\varphi\colon A^{(2)}\to A$ is the Verschiebung. If $A$ is semistable, the coefficient of $k$ is $n\deg\omega_A$.
\end{corollary}

\begin{proof}
Put $E=A^{(2)}$, $E'=A$, and apply Proposition~\ref{prop:ratio} and Lemma~\ref{lem:Dsize} over $K_k$, whose constant
field has $2^{nk}$ elements. The sum over the places of $K_k$ of the terms $f_wv_w(\vartheta_w)$ of
Lemma~\ref{lem:Dsize} is $nk(\deg\omega_E-\deg\omega_{E'})$ (the degrees of $\omega$ do not change under constant
field extension), so
\[
\dim\Sel^\varphi(E/K_k)\geq\dim\Sel^\varphi(E/K_k)-\dim\Sel^{\hat\varphi}(E'/K_k)
=1+\ord_2\prod_w\frac{c_w(E')}{c_w(E)}+nk\bigl(\deg\omega_E-\deg\omega_{E'}\bigr),
\]
the product being over the places $w$ of $K_k$. A place $v$ of $K$ of degree $d$ has at most $d$ places above it in
$K_k$, and the Tamagawa numbers are bounded by the orders of the geometric component groups; so the second term is
bounded independently of $k$. Next, $\dim\Sha(E/K_k)[\varphi]=\dim\Sel^\varphi(E/K_k)-\dim E'(K_k)/\varphi E(K_k)$, and
$\dim E'(K_k)/\varphi E(K_k)\leq\dim E'(K_k)/2E'(K_k)\leq\rk E'(K\bar\F_q)+2$, which is finite by the theorem of Lang and
N\'eron \cite{LN} (the points of the $K\bar\F_q/\bar\F_q$-trace over $\bar\F_q$ are torsion). The last assertion is Remark~\ref{rem:hodge}.
\end{proof}

The bound is vacuous unless $\deg\omega_{A^{(2)}}>\deg\omega_A$. If $A$ is constant, both bundles are trivial and
nothing is asserted (and for isotrivial $A$ the difference may vanish); if $A$ is semistable and not isotrivial, then
$12\deg\omega_A=\deg\Delta_{\min}(A)>0$, as $A$ has a place of multiplicative reduction. For the curves $E'$ of
Theorem~\ref{thm:example} one has $\deg\omega_{E'}=1$ and $\deg\omega_E=2$ (proof of
Proposition~\ref{prop:exsplit}), so $|\Sha(E/K)[\varphi]|$ grows at least like $q$, up to a bounded factor, in the
family of Theorem~\ref{thm:example}; Theorem~\ref{thm:example}(c) exhibits a subgroup of order $q/16$ on which the
pairing is non-degenerate.

\subsection{Root numbers at semistable places}\label{subsec:rootnumbers}
In this subsection $A\colon y^2+xy=x^3+a_2x^2+b$ is an ordinary elliptic curve over $K$, $E'=A$, $E\simeq A^{(2)}$, and
$v$ is a place. We use the following values of local root numbers: $w_v(A)=1$ if $A$ has good reduction at $v$,
$w_v(A)=-1$ if $A$ has split multiplicative reduction, and $w_v(A)=1$ if $A$ has non-split multiplicative reduction.
They follow from Deligne's theory of local constants, valid over local fields of any characteristic: an unramified
representation has root number $1$, and for the special representation $\chi\otimes\mathrm{sp}(2)$ with $\chi$
unramified one has $w=-\chi(\mathrm{Frob})$ (see \cite{Ro}, and \cite[\S19]{Ro} for the multiplicative case) (Rohrlich's definitions are used in characteristic
$p$ in \cite[\S1]{TW}). By the product formula for
the constant of the functional equation \cite[Thm.~3.2.1.1]{La}, the global root number $w(A/K)=\prod_vw_v(A)$ is the
sign in the functional equation of $L(A/K,s)$, so that $(-1)^{\ord_{s=1}L(A/K,s)}=w(A/K)$; see \cite[\S1]{TW}. We write $\tilde a$ for the residue class of $a\in\mathcal O_v$.

\begin{remark}\label{rem:uniquemodel}
Over a field $F$ of characteristic $2$, two models $y^2+xy=x^3+a_2x^2+b$ and $y^2+xy=x^3+a'_2x^2+b'$ are isomorphic if
and only if $b'=b$ and $a'_2-a_2\in\wp(F)$. Indeed, a change of coordinates $x=u^2x'+r$, $y=u^3y'+u^2sx'+t$ between them
must have $u=1$ (compare $a_1$), then $r=0$ (compare $a_3$) and $t=0$ (compare $a_4$), and then $a'_2=a_2+s^2+s$ and
$b'=b$ \cite[III.1]{Si}. Hence $[a_2,b)_v$ depends only on $A$ and $v$, and may be computed from any model of this shape
over $K_v$.
\end{remark}

\begin{lemma}[good reduction; for $k=0$ cf.\ {\cite[Prop.~2.4]{Kr}}]\label{lem:good}
Suppose that $A$ has good reduction at $v$, and let $k\geq0$ be the valuation of the coefficient $a_1$ of a minimal
Weierstrass model of $A$ over $\mathcal O_v$; thus $k=0$ if and only if the reduction is ordinary. Then
\[
\dim D_v-1=f_vk\qquad\text{and}\qquad[a_2,b)_v\equiv f_vk\pmod2 .
\]
In particular $w_v(A)=1=(-1)^{\dim D_v-1+[a_2,b)_v}$.
\end{lemma}

\begin{proof}
Let $y^2+a_1xy+a_3y=x^3+a^*_2x^2+a_4x+a_6$ be a minimal model over $\mathcal O_v$. In characteristic $2$ its
discriminant and $j$-invariant are
\[
\Delta=a_1^4b_8+a_3^4+a_1^3a_3^3,\qquad b_8=a_1^2a_6+a_1a_3a_4+a^*_2a_3^2+a_4^2,\qquad j=a_1^{12}/\Delta,
\]
and $\Delta\in\mathcal O_v^\times$. As $A$ is ordinary, $j\neq0$, so $a_1\neq0$; the reduction is supersingular if and
only if its $j$-invariant $\tilde a_1^{12}/\tilde\Delta$ vanishes \cite[Appendix~A]{Si}, i.e.\ if and only if $k>0$. The change of
coordinates $x=a_1^2x'+r$, $y=a_1^3y'+(a_4+r^2)/a_1$ with $r=a_3/a_1$ transforms the model into
\[
y^2+xy=x^3+a'_2x^2+\Delta/a_1^{12},\qquad a'_2=\frac{a^*_2}{a_1^2}+\frac{a_3}{a_1^3}
\]
\cite[III.1]{Si}. By Remark~\ref{rem:uniquemodel}, $b=\Delta/a_1^{12}$ and $[a_2,b)_v=[a'_2,\Delta)_v$, as $a_1^{12}$ is a
square.

\emph{The symbol.} If $k=0$, then $a'_2\in\mathcal O_v$ and $\Delta\in\mathcal O_v^\times$, so $[a'_2,\Delta)_v=0$ by
Fact~\ref{fact:schmid}. Let $k\geq1$. Then $a_3\in\mathcal O_v^\times$, because $\Delta\equiv a_3^4\pmod{\pi_v}$. Put
$s=a_1^3/a_3$, so that $v(s)=3k$, and write $\Delta=a_3^4(1+s)(1+z)$; then $v(z)\geq4k$. If $c\in K_v$ and $y\in\pi_v\mathcal O_v$
satisfy $v(c)+v(y)\geq1$, then $c\,dy/(1+y)$ is a regular differential, so $[c,1+y)_v=0$ by Fact~\ref{fact:schmid}. As
$v(a'_2)\geq-3k$ and $v(a^*_2/a_1^2)\geq-2k$, this gives $[a'_2,1+z)_v=0$ and $[a^*_2/a_1^2,1+s)_v=0$; also
$[a'_2,a_3^4)_v=0$. Since $a_3/a_1^3=1/s$,
\[
\begin{aligned}
[a'_2,\Delta)_v&=[1/s,1+s)_v=\Tr_{k_v/\F_2}\Res\Bigl(\frac{ds}{s(1+s)}\Bigr)\\
&=\Tr_{k_v/\F_2}\Res\Bigl(\frac{ds}{s}-\frac{ds}{1+s}\Bigr)=\Tr_{k_v/\F_2}\bigl(v(s)\bigr)\equiv3kf_v\equiv kf_v .
\end{aligned}
\]

\emph{The local image.} Both Tamagawa numbers are $1$. The change of coordinates above, with $u=a_1$, multiplies the
invariant differential by $u$ \cite[III.1, Table~3.1]{Si}, so the N\'eron differential of $A=E'$ is
$\omega'_v=a_1^{-1}\eta'$, where $\eta'=dX/X$ (by Remark~\ref{rem:uniquemodel} the model obtained is isomorphic to
\eqref{eq:E} by a change of coordinates with $u=1$, which preserves $\eta'$). The Frobenius twist of the minimal model,
with coefficients $a_i^2$ and unit discriminant $\Delta^2$, is a minimal model of $A^{(2)}\simeq E$, and the same
computation with $a_i^2$ in place of $a_i$ gives $\omega_v=a_1^{-2}\eta$ with $\eta=dx/x$. As $\varphi^*\eta'=\eta$
(proof of Lemma~\ref{lem:Dsize}), $\vartheta_v=a_1$, and Lemma~\ref{lem:Dsize} gives $\dim D_v-1=f_v\,v(a_1)=f_vk$.
\end{proof}

\begin{lemma}[multiplicative reduction; for $D_v$ cf.\ {\cite[Prop.~2.5]{Kr}}]\label{lem:mult}
Suppose that $A$ has multiplicative reduction at $v$. Then $A$ has a model $y^2+xy=x^3+a'_2x^2+b$ over $K_v$ with
$a'_2\in\mathcal O_v$; put $m=v(b)$. Then $m>0$, the reduction is split if and only if
$\Tr_{k_v/\F_2}(\tilde a'_2)=0$, and:
\begin{enumerate}[label=\textup{(\roman*)},leftmargin=*]
\item if the reduction is split, then $D_v=0$ and $[a_2,b)_v=0$;
\item if the reduction is non-split, then $\dim D_v-1$ is $0$ for even $m$ and $-1$ for odd $m$, and
$[a_2,b)_v\equiv m\pmod2$.
\end{enumerate}
In both cases $w_v(A)=(-1)^{\dim D_v-1+[a_2,b)_v}$.
\end{lemma}

\begin{proof}
A minimal model of $A$ at $v$ has $c_4=a_1^4$ (characteristic $2$) of valuation $0$, so $a_1\in\mathcal O_v^\times$, and
the change of coordinates of the proof of Lemma~\ref{lem:good} is defined over $\mathcal O_v$. It gives an integral model
$y^2+xy=x^3+a'_2x^2+b'$, and $b'=b$ by Remark~\ref{rem:uniquemodel}; thus $m=v(\Delta)>0$. An integral model with
$a_1=1$ is minimal: for any other integral model $a_1/u$ is integral (characteristic $2$), so $v(u)\leq0$ and the
valuation of the discriminant does not decrease. The same applies to the model $y^2+xy=x^3+a'_2x^2+b^2$ of $E$
(Remark~\ref{rem:uniquemodel}). Their discriminants $b$ and $b^2$ have valuations $m$ and $2m$, so $A$ and $E$ have
reduction of types $\mathrm I_m$ and $\mathrm I_{2m}$. The node of $y^2+xy=x^3+\tilde a'_2x^2$ is $(0,0)$, with tangents
$y^2+xy+\tilde a'_2x^2=0$, which are defined over $k_v$ if and only if $\tilde a'_2\in\wp(k_v)$, i.e.\
$\Tr_{k_v/\F_2}(\tilde a'_2)=0$; the same criterion applies to $E$. The N\'eron differentials are $dX/X$ and $dx/x$, so
$\vartheta_v=1$. By Tate's algorithm \cite[IV.9]{Si2}, the Tamagawa number of a curve of type $\mathrm I_m$ is $m$ if the
reduction is split, and $2$ or $1$ according as $m$ is even or odd if it is non-split. Lemma~\ref{lem:Dsize} gives
$|D_v|=2c_v(A)/c_v(E)$, which is $2m/2m=1$ in the split case and $c_v(A)$ in the non-split case. If the reduction is
split, then $a'_2\in\wp(K_v)$ by Hensel's lemma, so $[a_2,b)_v=[a'_2,b)_v=0$. If it is non-split, Fact~\ref{fact:schmid}
gives $[a'_2,b)_v=\Tr_{k_v/\F_2}(\tilde a'_2)\,v(b)\equiv m$. The root number is $-1$, resp.\ $+1$, and in each case the
exponent $\dim D_v-1+[a_2,b)_v$ has the right parity.
\end{proof}

\subsection{Proof of Theorem~\ref{thm:semistable}}
\begin{proof}
At a place of semistable reduction the reduction is good or multiplicative, and the local identity is
Lemma~\ref{lem:good} or Lemma~\ref{lem:mult}. By reciprocity (Fact~\ref{fact:ABHN}), $\sum_v[a_2,b)_v=0$, so by
Theorem~\ref{thm:parity}
\[
(-1)^{s_2(A/K)}=\prod_v(-1)^{\dim D_v-1}=\prod_v(-1)^{\dim D_v-1+[a_2,b)_v} .
\]
The factor at a place of semistable reduction is $w_v(A)$, and since $w(A/K)=\prod_vw_v(A)$ and $w_v(A)^2=1$, this is
the asserted formula.
\end{proof}

\begin{proofof}{Corollary~\ref{cor:delta}}
Write $\Sha(A/K)[2^\infty]\simeq(\Q_2/\Z_2)^c\oplus M$ with $M$ finite; then
$\dim\Sha(A/K)[2]=c+\dim M[2]=c+\delta(A)$. So $s_2(A/K)=\rk A(K)+c+\delta(A)$, and $\rk A(K)+c$ is the corank of the $2^\infty$-Selmer group. By \cite[Thm.~1.1]{TY}
with $\ell=2$ this corank has the parity of $\ord_{s=1}L(A/K,s)$, and $(-1)^{\ord_{s=1}L(A/K,s)}=w(A/K)$
(\S\ref{subsec:rootnumbers}). Hence $(-1)^{s_2(A/K)}=w(A/K)(-1)^{\delta(A)}$, and the formula follows from
Theorem~\ref{thm:semistable}.
\end{proofof}

\subsection{An example}\label{subsec:exparity}
Let $K=\F_2(S)$, $t=1/S$, and
\[
A\colon y^2+xy=x^3+S^3x^2+S^{12}+S^9,\qquad b=S^{12}+S^9=S^9(S+1)(S^2+S+1).
\]
At the finite places not dividing $b$ the curve has good ordinary reduction. At the places $S$, $S+1$ and $S^2+S+1$
(the last of degree $2$, with residue field $\F_4$) it has multiplicative reduction with $m=9$, $1$, $1$; the residue
of $a_2=S^3$ is $0$, $1$, $1$, and $\Tr_{\F_2/\F_2}(1)=1$, $\Tr_{\F_4/\F_2}(1)=0$, so the reduction is split, non-split
and split respectively (Lemma~\ref{lem:mult}). At $\infty$, the minimal model $y^2+txy+y=x^3$, with $a_1=t$, $a_3=1$ and
$\Delta=1+t^3$, is transformed by the change of coordinates of the proof of Lemma~\ref{lem:good} into $y^2+xy=x^3+t^{-3}x^2+(1+t^3)t^{-12}$,
which is the given model; so $A$ has good supersingular reduction at $\infty$, with $k=1$. By Lemmas~\ref{lem:good}
and~\ref{lem:mult}:
\begin{center}\small
\begin{tabular}{lcccccc}
\toprule
$v$ & reduction & $c_v(A)$ & $c_v(A^{(2)})$ & $\dim D_v-1$ & $[a_2,b)_v$ & $w_v(A)$\\
\midrule
$S$ & split $\mathrm I_9$ & $9$ & $18$ & $-1$ & $0$ & $-1$\\
$S+1$ & non-split $\mathrm I_1$ & $1$ & $2$ & $-1$ & $1$ & $+1$\\
$S^2+S+1$ & split $\mathrm I_1$ & $1$ & $2$ & $-1$ & $0$ & $-1$\\
$\infty$ & good, supersingular & $1$ & $1$ & $1$ & $1$ & $+1$\\
\bottomrule
\end{tabular}
\end{center}
At all other places $\dim D_v=1$ and $[a_2,b)_v=0$. Hence $\sum_v(\dim D_v-1)=-2$, and $s_2(A/K)$ is even. In
Theorem~\ref{thm:parity} the Tamagawa term $\ord_2\prod_vc_v(A)/c_v(A^{(2)})=-3$ is odd, and it is compensated by the
Hodge term: $\vartheta_v$ is a unit at all finite places and $\vartheta_\infty=t$, so
$\deg\omega_{A^{(2)}}-\deg\omega_A=1$ by Lemma~\ref{lem:Dsize}. The global root number is $w(A/K)=(-1)(+1)(-1)(+1)=+1$,
as predicted by Theorem~\ref{thm:semistable}. This agrees with the $L$-function: the conductor of $A$ has degree
$1+1+2=4$, so $L(A/K,T)$ is a polynomial of degree $4-4=0$ by the Grothendieck--Ogg--Shafarevich formula
\cite{Ul2}, i.e.\ $L(A/K,T)=1$. In particular $\rk A(K)=0$, as the rank is at most the order of vanishing of the
$L$-function \cite{Ta66}. In fact $\Sha(A/K)=0$. The point $(S^3,S^6)$ lies on $A$ and has order $3$: in characteristic
$2$ the $3$-division polynomial of $y^2+xy=x^3+a_2x^2+b$ is $x^4+x^3+b$, and $x=S^3$ is a root. The surface attached to $A$
is rational ($\deg\omega_A=1$), with trivial lattice $U\oplus A_8$ (the fibre $\mathrm I_9$ at $S=0$; the other singular
fibres over $\bar\F_2$ are of type $\mathrm I_1$), so $A(\bar\F_2(S))$ is finite of order $\sqrt{\det A_8}=3$ \cite{Shi90}, \cite{OS}.
Thus $A(K)\simeq\Z/3$. As $\rk A(K)=\ord_{T=1/q}L(A/K,T)=0$, the Birch and Swinnerton-Dyer formula holds \cite{KT}, and
with $L(A/K,1/q)=1$, $\prod_vc_v(A)=9$ and $\deg\omega_A=1$ it gives $|\Sha(A/K)|=|A(K)_{\mathrm{tors}}|^2/9=1$ (see
\S\ref{subsec:exrevisited} for the normalisation). Hence $s_2(A/K)=0$.

\subsection{The examples of \S\ref{sec:example} revisited}\label{subsec:exrevisited}
In the examples of \S\ref{sec:example} the groups $\Sel^\varphi(E/K)$, $E(K)$ and $\Sha(E/K)$ can be determined
completely. The $2$-primary part of $\Sha$ is obtained by descent, from Theorem~\ref{thm:radical} and the pairing
matrices of \S\ref{sec:example}; the Birch and Swinnerton-Dyer formula then gives the odd part, and an independent
check of the $2$-primary part. We use the following facts, for a non-isotrivial elliptic curve $A$ over $K=\F_q(S)$. The $L$-function is a
polynomial $L(A/K,T)$ in $T=q^{-s}$ of degree $\deg\mathfrak n_A-4$, where $\mathfrak n_A$ is the conductor
\cite{Ul2}. With $A_m=\sum_{s\in\mathbb P^1(\F_{q^m})}a_s(m)$, where $a_s(m)=q^m+1-\#A_s(\F_{q^m})$ for fibres of good
reduction and $a_s(m)\in\{1,-1,0\}$ for split multiplicative, non-split multiplicative and additive fibres, one has
$L(A/K,T)=\exp\bigl(\sum_mA_mT^m/m\bigr)$; for $A\colon y^2+xy=x^3+a_2x^2+b$ and a fibre of good or multiplicative
reduction, $a_s(m)=-(-1)^{\Tr a_2(s)}\sum_{x\in\F_{q^m}^\times}(-1)^{\Tr(x+\sqrt{b(s)}/x)}$, a Kloosterman sum. By Tate,
$\rk A(K)\leq r:=\ord_{T=1/q}L(A/K,T)$ \cite{Ta66}; if equality holds, then $\Sha(A/K)$ is finite and
\[
L^*(A)=\frac{\operatorname{Reg}(A)\,|\Sha(A/K)|\,\prod_vc_v(A)\,q^{1-\deg\omega_A}}{|A(K)_{\mathrm{tors}}|^2}
\]
\cite{KT,Ul2}, where $L^*(A)$ is the value at $T=1/q$ of $L(A/K,T)(1-qT)^{-r}$ and $\operatorname{Reg}(A)$ is the determinant of the
height pairing of Shioda \cite{Shi90}, the canonical height pairing divided by $\log q$. For a rational elliptic
surface ($\deg\omega_A=1$) and sections $P,Q$, $\langle P,P\rangle=2+2(P\cdot O)-\sum_v\operatorname{contr}_v(P)$ and
$\langle P,Q\rangle=1+(P\cdot O)+(Q\cdot O)-(P\cdot Q)-\sum_v\operatorname{contr}_v(P,Q)$, where
$\operatorname{contr}_v=0$ if $P$ or $Q$ meets the identity component, $\operatorname{contr}_v=\frac12$ for fibres
of types $\mathrm I_2$, $\mathrm{III}$, and $\operatorname{contr}_v(P)=1$ for $\mathrm I_0^*$; and the free part of
the Mordell--Weil lattice over $\bar\F_q(S)$ has determinant $|A(\bar\F_q(S))_{\mathrm{tors}}|^2/\det\mathcal T$, $\mathcal T$ the
root lattice of the reducible fibres \cite{Shi90}. Moreover, for a rational elliptic surface $X$ all of
$H^2(X_{\bar\F_q},\Q_\ell(1))$ is spanned by divisor classes, the N\'eron--Severi group modulo the classes of the zero
section and of the fibre components is $A(\bar\F_q(S))$ \cite[Thm.~1.3]{Shi90}, and consequently
\begin{equation}\label{eq:Lrational}
L(A/K,T)=\det\bigl(1-qT\,\mathrm{Fr}\mid A(\bar\F_q(S))\otimes\Q\bigr),
\end{equation}
where $\mathrm{Fr}$ is the $q$-Frobenius, an automorphism of finite order of the Mordell--Weil group over $\bar\F_q(S)$
\cite{Ul2}; in particular $\ord_{T=1/q}L(A/K,T)=\rk A(K)$. Reduction types and Tamagawa numbers are computed by Tate's
algorithm \cite[IV.9]{Si2}.

\begin{proposition}\label{prop:exnonsplit}
Let $K=\F_8(S)$, $E$, $E'$ and $b=S^4+S^3+\varrho$ be as in Theorem~\ref{thm:nonsplit}.
\begin{enumerate}[label=\textup{(\alph*)},leftmargin=*]
\item $\Sel^\varphi(E/K)=\langle\chi_1,\chi_{\varrho S},\chi_{\varrho^2S}\rangle$ and $\Sel^{\hat\varphi}(E'/K)=\langle[b]\rangle$.
\item $P'=(\varrho^5+\varrho S,\ S^2+\varrho^6S)\in E'(K)$ has infinite order, $\delta_\varphi(P')=\chi_1+\chi_{\varrho S}$,
$E'(K)=\Z P'$ and $E(K)=\Z\hat\varphi(P')\oplus\{O,T\}$.
\item $\Sha(E'/K)[2^\infty]=0$, and $\Sha(E/K)[2^\infty]=\Sha(E/K)[\varphi]\simeq(\Z/2)^2$ is spanned by $\iota(\chi_1)$
and $\iota(\chi_{\varrho^2S})$.
\item $L(E/K,T)=1-8^4T^4$, $\Sha(E'/K)=0$ and $\Sha(E/K)=\Sha(E/K)[\varphi]$.
\end{enumerate}
Parts \textup{(a)}--\textup{(c)} are proved by descent, without $L$-functions.
\end{proposition}

\begin{proof}
\emph{Local data.} At $v_b$ both curves have split multiplicative reduction ($a_2=0$), of types $\mathrm I_1$ and
$\mathrm I_2$, so $D_{v_b}=0$ (Lemma~\ref{lem:mult}). At $\infty$, with $t=1/S$ and $e_0=1+t+\varrho t^4$, the change of
coordinates with $u=t^{-1}$ followed by $y\mapsto y+t$ gives the model $y^2+txy=x^3+t^2x+t^3+\varrho t^6$ of $E'$, with
$v(\Delta)=8$; in Tate's algorithm the polynomial $U^3+a_{2,1}U^2+a_{4,2}U+a_{6,3}$ reduces to $U^3+U+1$, with the
three roots $\varrho,\varrho^2,\varrho^4\in\F_8$, so $E'$ has type $\mathrm I_0^*$ and $c_\infty(E')=4$. For $E$ the change
of coordinates with $u=t^{-1}$ followed by $y\mapsto y+t^{-1}e_0$ gives $y^2+txy=x^3+e_0x$, with $v(\Delta)=4$; after
$x\mapsto x+1$ one has $v(a_6)=1$, so $E$ has type $\mathrm{II}$ and $c_\infty(E)=1$. Both changes of coordinates have
$u=t^{-1}$, so $\vartheta_\infty=1$, and Lemma~\ref{lem:Dsize} gives $|D_\infty|=2\cdot4/1=8$. Ogg's formula \cite{Sai} gives the
conductor exponent $8+1-5=4$ at $\infty$, so $\deg\mathfrak n=4+4$ and $\deg L(E/K,T)=4$.

(a) A class $\chi\in\Sel^\varphi(E/K)$ is unramified at all finite places ($D_v=H^1_\nr$ at the good places, $D_{v_b}=0$),
so it is represented by a polynomial in $S$, and after replacing $cS^{2j}$ by $\sqrt cS^j$ by a unique polynomial
$c_0+\sum_{j\ \mathrm{odd}}c_jS^j$ with $c_0\in\{0,1\}$. The same expression is the normal form of $\chi_\infty$ in
$K_\infty/\wp K_\infty$, so restriction to $K_\infty$ is an isomorphism from the group of characters unramified at all
finite places onto $H^1(K_\infty,\Z/2)$, and $\Sel^\varphi(E/K)$ is identified with the set of $d\in D_\infty$ whose
global representative $\chi_d$ is trivial at $v_b$. This condition is automatic: by reciprocity
$\sum_v\inv_v(\chi_d,b)=0$, the term at $\infty$ vanishes because $[b]=\kappa(T)\in D_\infty^\perp$
(Lemma~\ref{lem:annihilator}), the terms at the other finite places $v\neq v_b$ vanish, and the term at $v_b$ is the
obstruction to $\chi_{d,v_b}=0$. Hence $\dim\Sel^\varphi(E/K)=\dim D_\infty=3$, and Theorem~\ref{thm:nonsplit}(a) gives
the first assertion. A class $u\in\Sel^{\hat\varphi}(E'/K)$ has even valuation at the good finite places, where
$\kappa(E(K_v))=D_v^\perp=U_v$, so $u\in\langle[b]\rangle$; and $[b]=\kappa(T)$. (As a check,
$3-1=1+(0-1)+(3-1)$, as in Proposition~\ref{prop:ratio}.)

(b) The point $P'$ satisfies the equation of $E'$, and $\delta_\varphi(P')=[\varrho^5+\varrho S]=\chi_1+\chi_{\varrho S}$
as $\Tr(\varrho^5)=1$. The surface of $E'$ is rational ($\deg\Delta_{\min}=4\cdot1+8=12$). As $P'$ is an integral
section ($\deg x\leq2$, $\deg y\leq3$), $(P'\cdot O)=0$; at $\infty$ it passes through the singular point of the
reduction, i.e.\ it meets a non-identity component of the fibre $\mathrm I_0^*$, and the fibres $\mathrm I_1$ do not
contribute. So $\langle P',P'\rangle=2-1=1$. All self-heights $\langle Q,Q\rangle$ on this surface are integers, the only reducible fibre being of type
$\mathrm I_0^*$, where $\operatorname{contr}(Q)\in\{0,1\}$; so $P'$ is not a multiple of
another point modulo torsion, and there is no non-zero torsion, since a torsion section $Q$ would have
$0=\langle Q,Q\rangle\geq2-1$. Next, $E(K)/\hat\varphi E'(K)\simeq\kappa(E(K))\subseteq\Sel^{\hat\varphi}(E'/K)=\langle\kappa(T)\rangle$,
so $E(K)=\hat\varphi E'(K)+\{O,T\}$, and $E(K)_{\mathrm{tors}}=\{O,T\}$ because $\varphi$ maps it into
$E'(K)_{\mathrm{tors}}=0$. It remains to show that $\rk E'(K)\leq1$. The group $\delta_\varphi(E'(K))$ is the kernel of
$\iota$ on $\Sel^\varphi(E/K)$, so it lies in the radical of the form $B$ of Corollary~\ref{cor:bilin}, which by (a) and
Theorem~\ref{thm:nonsplit}(b) is $\langle\chi_1+\chi_{\varrho S}\rangle$; hence $|E'(K)/\varphi E(K)|\leq2$. In the exact
sequence
\[
E(K)/\hat\varphi E'(K)\xrightarrow{\ \varphi\ }E'(K)/2E'(K)\longrightarrow E'(K)/\varphi E(K)\longrightarrow0
\]
the first group is generated by the class of $T$, and $\varphi(T)=O$. So $|E'(K)/2E'(K)|\leq2$, and as $E'(K)$ is
torsion-free, $\rk E'(K)\leq1$. Thus $E'(K)=\Z P'$ and $\delta_\varphi(E'(K))=\langle\chi_1+\chi_{\varrho S}\rangle$.

(c) By (a) and (b), $\Sha(E/K)[\varphi]\simeq\Sel^\varphi(E/K)/\delta_\varphi(E'(K))$ has order $4$ and is spanned by
$\iota(\chi_1)$ and $\iota(\chi_{\varrho^2S})$, and $\Sha(E'/K)[\hat\varphi]\simeq\Sel^{\hat\varphi}(E'/K)/\kappa(E(K))=0$.
By Theorem~\ref{thm:nonsplit}(b) the Cassels--Tate pairing is non-degenerate on $\Sha(E/K)[\varphi]$, so
$\hat\varphi(\Sha(E'/K)[2])=0$ by Theorem~\ref{thm:radical}(b); as $\Sha(E'/K)[\hat\varphi]=0$, this gives
$\Sha(E'/K)[2]=0$, hence $\Sha(E'/K)[2^\infty]=0$. If $x\in\Sha(E/K)[2^\infty]$, then
$\varphi(x)\in\Sha(E'/K)[2^\infty]=0$, so $x\in\Sha(E/K)[\varphi]$.

(d) The surface of $E'$ is rational and $\deg L=4$, so by \eqref{eq:Lrational}
$L(E'/K,T)=\prod_{i=1}^4(1-8\zeta_iT)$ with roots of unity $\zeta_i$, of which exactly one is $1$ because
$\rk E'(K)=1$. As $L$ has integer coefficients, $\prod_i(x-\zeta_i)=(x-1)g(x)$ where $g$ is a product of cyclotomic
polynomials $\Phi_n$, $n\geq2$, of degree $3$: $g\in\{\Phi_2^3,\ \Phi_2\Phi_3,\ \Phi_2\Phi_4,\ \Phi_2\Phi_6\}$, and
$\sum_i\zeta_i=-2,\,-1,\,0,\,1$ respectively. On the other hand $-8\sum_i\zeta_i=A_1=\sum_{s\in\F_8}a_s(1)$ (the fibre
at $\infty$ is additive, and $b$ has no root in $\F_8$), where $a_s(1)=-\mathrm{Kl}(\sqrt{b(s)})$ with
$\mathrm{Kl}(c)=\sum_{x\in\F_8^\times}(-1)^{\Tr(x+c/x)}$. The sum $\mathrm{Kl}(c)$ depends only on the orbit of $c$ under
$c\mapsto c^2$, and with the traces listed in \S\ref{subsec:nonsplit} one finds
\[
\mathrm{Kl}(1)=-5,\qquad\mathrm{Kl}(\varrho)=-1,\qquad\mathrm{Kl}(\varrho^3)=3
\]
(for instance, for $c=\varrho$ and $x=\varrho^i$, $0\leq i\leq6$, the trace of $\varrho^i+\varrho^{1-i}$ is $1,1,1,0,0,0,1$;
as a check, $\sum_{c\neq0}\mathrm{Kl}(c)=1$ and $\sum_{c\neq0}\mathrm{Kl}(c)^2=55$). The values of $b$ at
$s=0,1,\varrho,\varrho^2,\varrho^3,\varrho^4,\varrho^5,\varrho^6$ are
$\varrho,\varrho,\varrho^5,\varrho^6,1,1,\varrho^6,\varrho^5$, so $A_1=2\cdot1+4\cdot(-3)+2\cdot5=0$. Hence
$g=\Phi_2\Phi_4$ and $L(E'/K,T)=(1-8T)(1+8T)(1+64T^2)=1-8^4T^4$; isogenous curves have the same $L$-function. Now
$r=1=\rk E'(K)$, and $L^*=(1+8T+64T^2+512T^3)|_{T=1/8}=4$. For $E'$: $\operatorname{Reg}=1$, $\prod c_v=1\cdot4$,
$\deg\omega_{E'}=1$ and trivial torsion, so $4=4\,|\Sha(E'/K)|$ and $\Sha(E'/K)=0$. For $E$:
$\operatorname{Reg}(E)=\langle\hat\varphi P',\hat\varphi P'\rangle=2$, $\prod c_v=2\cdot1$ (types $\mathrm I_2$ and
$\mathrm{II}$), $\deg\omega_E=1$ ($\deg\Delta_{\min}=2\cdot4+4$) and $|E(K)_{\mathrm{tors}}|=2$, so
$4=2\,|\Sha(E/K)|\cdot2/4$ and $|\Sha(E/K)|=4$, in accordance with (c).
\end{proof}

In particular the radical $\langle\chi_1+\chi_{\varrho S}\rangle$ of the matrix of Theorem~\ref{thm:nonsplit}(b) is
$\delta_\varphi(E'(K))=\hat\varphi_*\Sel^2(E'/K)$, as predicted by Theorem~\ref{thm:radical}, and $s_2(E'/K)=1$, in
accordance with Theorem~\ref{thm:parity}: $\sum_v(\dim D_v-1)=-1+2$.

\begin{proposition}\label{prop:exsplit}
Let $k=2$, $K=\F_{256}(S)$, and let $E$, $E'$, $\Lambda$, $Z$ be as in Theorem~\ref{thm:example}.
\begin{enumerate}[label=\textup{(\alph*)},leftmargin=*]
\item $\Sel^\varphi(E/K)=\{\chi_\lambda:\lambda\in Z^\perp\}$ and $\delta_\varphi(E'(K))=\{\chi_\lambda:\lambda\in Z\}$.
\item $L(E/K,T)=(1-256T)^3$ and $\rk E(K)=3$.
\item $\Sha(E'/K)=0$, and $\Sha(E/K)=\Sha(E/K)[\varphi]\simeq Z^\perp/Z\simeq(\Z/2)^4$.
\end{enumerate}
\end{proposition}

\begin{proof}
Write $\zeta=\zeta_3$ and $\Lambda_\infty=t^4\Lambda(1/t)=1+\zeta t^2+\zeta^2t^3$.

\emph{Local data.} At the four roots of $\Lambda$ the curves have split multiplicative reduction of types $\mathrm I_2$
($E'$) and $\mathrm I_4$ ($E$), so $D_v=0$ there. The substitution $y\mapsto y+\Lambda$ gives the global model
$y^2+xy=x^3+\Lambda x$ of $E'$, with $\deg a_i\leq i$, so the surface of $E'$ is rational; at $\infty$ it becomes
$y^2+txy=x^3+\Lambda_\infty x$, with $v(\Delta)=4$, and after $x\mapsto x+1$ one has $\pi^2\mid a_6$, $\pi^3\nmid b_8$: type
$\mathrm{III}$, $c_\infty(E')=2$. For $E$ one finds the minimal model $y^2+t^2xy=x^3+\Lambda_\infty^2x$ at $\infty$, with
$v(\Delta)=8$; Tate's algorithm leads to the polynomial $U^3+U=U(U+1)^2$ and then to $Y^2+Y+\zeta$, which has two
roots in $\F_{256}$ since $\Tr_{\F_{256}/\F_2}(\zeta)=0$: type $\mathrm I_1^*$. For $\mathrm I_n^*$ with $n$ odd the
geometric component group is $\Z/4$, and the two far components correspond to the roots of this quadratic
\cite[IV, \S9]{Si2}; as both are rational, $c_\infty(E)=4$. The Néron differentials
at $\infty$ are $t^{-1}\,dx/x$ for $E'$ and $t^{-2}\,dx/x$ for $E$, so $\vartheta_\infty=t$, and Lemma~\ref{lem:Dsize} gives
$|D_\infty|=2\cdot\frac24\cdot256=2^8$. Ogg's formula \cite{Sai} gives the conductor exponent $4+1-2=3$ at $\infty$, and
$\deg\omega_E=2$.

(a) For $\lambda,d\in\F_q$ with $\Tr(\lambda)=0$, put $c_\lambda=1+\lambda$ and $X=S^2+c_\lambda S+d$. Then
$\Lambda=X(S^2+c_\lambda S+c_\lambda^2+d+\zeta)+r$ with $\deg r\leq1$, and $r/X\in t\F_q[[t]]\subseteq\wp(K_\infty)$, so
$X+\Lambda^2/X^2\equiv X+\Lambda/X\equiv c_\lambda^2+\zeta\pmod{\wp K_\infty}$ (cf.\ \eqref{eq:division}), and
$\Tr(c_\lambda^2+\zeta)=\Tr(\lambda)=0$. So there is $R\in E'(K_\infty)$
with $x(R)=X$ and $\delta_\varphi(R)=[X]=[\lambda S+d]$. These classes form a group of order $2^7\cdot2=|D_\infty|$, hence
all of $D_\infty$. As in the proof of Proposition~\ref{prop:exnonsplit}, $\Sel^\varphi(E/K)$ is identified with the set of
$d\in D_\infty$ whose representative is trivial at the places $S=a$, $a\in Z$; for $\chi_{\lambda S+d}$ this means
$\Tr(\lambda a+d)=0$ for $a\in Z$, i.e.\ $d\in\wp(\F_q)$ and $\lambda\in Z^\perp$ (note $1\in Z$, so $\lambda\in Z^\perp$
implies $\Tr(\lambda)=0$). The second assertion is then Theorem~\ref{thm:example}(c), as $\delta_\varphi(E'(K))$ is the
kernel of $\iota$ on $\Sel^\varphi(E/K)$.

(b) The isogenous curves $E$ and $E'$ have the same $L$-function, and we compute it for $E'$, for which
$\sqrt{b(s)}=\Lambda(s)$. Since $\Lambda$ is additive, $\sum_{s\in\F_{q^m}}(-1)^{\Tr(\Lambda(s)/x)}$ vanishes unless $\Lambda^*(1/x)=0$, where
$\Lambda^*(c)=c^{1/4}+(\zeta c)^{1/2}+\zeta^2c$ is the adjoint of $\Lambda$ for the trace form. The non-zero roots of
$\Lambda^*$ are the roots of $c^3+c+\zeta$, which lie in $\F_{16}$, and $\Tr_{\F_{q^m}/\F_2}(1/c)=2m\Tr_{\F_{16}/\F_2}(1/c)=0$.
Hence $A_m=-3q^m$ for all $m$, and $L(E/K,T)=(1-qT)^3$. Let $Z\smallsetminus\{0\}=\{c,c',c''\}$; these are the roots of
$S^3+\zeta S+\zeta^2$, so $c''=c+c'$, $c^2+cc'+c'^2=\zeta$ and $cc'c''=\zeta^2$. Choose $w_c\in\F_q$ with
$w_c^2+w_c=c^2+\zeta$ (possible as $\Tr(c^2+\zeta)=0$). Since $\Lambda=(S^2+cS)(S^2+cS+c^2+\zeta)$, the point
$R_c=(S^2+cS,\ w_c(S^2+cS))$ lies on $y^2+xy=x^3+\Lambda x$; it is an integral section, so $(R_c\cdot O)=0$. At $S=a$,
$a\in Z$, the singular point of the fibre is $(0,0)$, and $R_c$ passes through it exactly for $a\in\{0,c\}$; at $\infty$
one has $t^2x(R_c)=1+ct$, and $R_c$ passes through the singular point $(1,0)$. So $R_c$ meets the non-identity
components at $S=0$, $S=c$ and $\infty$, and $\langle R_c,R_c\rangle=2-\frac32=\frac12$. For $c\neq c'$ the sections $R_c$
and $R_{c'}$ have the same $x$-coordinate only at $S=0$ and $\infty$. Both singular points of the surface there are of
type $A_1$. Such a point is resolved by one blow-up, whose exceptional curve, the non-identity component of the
fibre, is the projectivised tangent cone, i.e.\ the conic defined by the quadratic part of the equation, and a smooth
curve through the point meets it in the point given by its tangent direction:
$y^2+xy+\zeta^2Sx=0$ at $S=0$, and $x'^2+y^2+ty+\zeta t^2=0$ at $\infty$ (with $x'=t^2x+1$ and $t^3y$ in place of $y$).
The tangent directions of $R_c$, hence its intersection points with these conics, are $(x:y:S)=(c:w_cc:1)$ and
$(x':y:t)=(c:w_c:1)$, which are different for different $c$. Hence $(R_c\cdot R_{c'})=0$ and $\langle R_c,R_{c'}\rangle=1-0-(\frac12+\frac12)=0$. Thus the $R_c$ are linearly independent, $\rk E'(K)\geq3=\ord_{T=1/q}L$, and equality holds. For (b) only this
independence is needed; the value $\frac18$ of the determinant of the Gram matrix is used in (c).

(c) The lattice spanned by the $R_c$ has determinant $\frac18$. The free part of the Mordell--Weil lattice over
$\bar\F_q(S)$ has determinant $|\mathrm{tors}|^2/2^5$, and the torsion contains $T'=(0,\Lambda)$; hence the torsion is
$\{O,T'\}$ and $E'(K)/\{O,T'\}$ is spanned by the $R_c$, with $\operatorname{Reg}(E')=\frac18$. With $L^*=1$, $\prod c_v(E')=2^4\cdot2$
and $\deg\omega_{E'}=1$ the formula gives $1=\frac18|\Sha(E'/K)|\cdot32/4$, so $\Sha(E'/K)=0$. For $E$: the point
$Q=(\Lambda,\Lambda^2)\in E(K)$ satisfies $\varphi(Q)=T'$, so $E(K)_{\mathrm{tors}}=\langle Q\rangle\simeq\Z/4$ ($\varphi$ maps
it into $\{O,T'\}$ with kernel $\{O,T\}$). The image $\varphi(E(K))$ contains $E'(K)_{\mathrm{tors}}$ and has index
$|\delta_\varphi(E'(K))|=4$ in $E'(K)$, so $\operatorname{Reg}(E)=4^2\operatorname{Reg}(E')/2^3=\frac14$. With $\prod c_v(E)=4^4\cdot4$ and
$\deg\omega_E=2$ the formula gives $1=\frac14|\Sha(E/K)|\cdot1024\cdot256^{-1}/16$, so $|\Sha(E/K)|=16$ (the formula
applies to $E$ as well, since $\rk E(K)=3$). By (a),
$\Sha(E/K)[\varphi]\simeq\Sel^\varphi(E/K)/\delta_\varphi(E'(K))\simeq Z^\perp/Z$ has order $16$, so it is all of $\Sha(E/K)$.
\end{proof}

\subsection{Places of additive reduction: proof of Theorem~\ref{thm:local}}\label{subsec:additive}
Let $F$ be a local field of characteristic $2$ with valuation $v$, maximal ideal $\mathfrak m$ and residue field of
order $q=2^f$. For an ordinary elliptic curve $A\colon y^2+xy=x^3+a_2x^2+b$ over $F$ put
\[
\epsilon(A):=w(A/F)\,(-1)^{\dim D(A)-1+[a_2,b)_F}\in\{\pm1\};
\]
it depends only on the isomorphism class of $A$ (Remark~\ref{rem:uniquemodel}). Theorem~\ref{thm:local}
and Remark~\ref{thm:additive} concern the identity $\epsilon(A)=1$. By Lemmas~\ref{lem:good} and~\ref{lem:mult} and the values of the root numbers recalled in
\S\ref{subsec:rootnumbers}, $\epsilon(A)=1$ if $A$ has semistable reduction.

\subsubsection{Potentially multiplicative reduction}
For $d\in F$ and $s\in F^s$ with $s^2+s=d$, the map $(x,y)\mapsto(x,y+sx)$ is an isomorphism from
$y^2+xy=x^3+a_2x^2+b$ to $y^2+xy=x^3+(a_2+d)x^2+b$ over $F(s)$, and for $\sigma$ with $\sigma s=s+1$ it differs from
its conjugate by $(x,y)\mapsto(x,y+x)$, which is $P\mapsto-P$. Hence $A$ is the quadratic twist by $\chi_{a_2}$ of
$A_0\colon y^2+xy=x^3+b$.

\begin{proposition}\label{prop:potmult}
Suppose $v(b)>0$. If $\chi_{a_2}=0$, then $D(A)=0$; if $\chi_{a_2}\neq0$, then $\dim D(A)=1-[a_2,b)_F$. If $\chi_{a_2}$
is ramified, then $A$ has additive, potentially multiplicative reduction, $w(A/F)=1$, and $\epsilon(A)=1$.
\end{proposition}

\begin{proof}
We use Tate's uniformisation over complete non-archimedean fields of any characteristic \cite{Roq} (for $p$-adic
fields see \cite[V, \S\S3, 5]{Si2}). There is a unique $\mathsf q\in F^\times$ with $v(\mathsf q)>0$ such that the Tate curve $E_{\mathsf q}$ has $j$-invariant $1/b$. Its equation $y^2+xy=x^3+a_4(\mathsf q)x+a_6(\mathsf q)$ has $a_1=1$,
$a_2=a_3=0$ and discriminant $\Delta(\mathsf q)=\mathsf q\prod_{n\geq1}(1-\mathsf q^n)^{24}$, so its model of the shape \eqref{eq:E} is
$y^2+xy=x^3+\Delta(\mathsf q)$ (proof of Lemma~\ref{lem:good}); comparing $j$-invariants, $E_{\mathsf q}\simeq A_0$ and $b=\Delta(\mathsf q)\in\mathsf qF^{\times2}$. The map $u\mapsto(X(u,\mathsf q),Y(u,\mathsf q))$, given by power series with coefficients in $\Z[[\mathsf q]]$,
induces an isomorphism $L^\times/\mathsf q^{\Z}\simeq E_{\mathsf q}(L)$ for every finite extension $L$ of $F$,
separable or not \cite{Roq}; being given by universal series, it is compatible with all $F$-embeddings between such
fields, so it is Galois-equivariant and compatible with the inclusion $L\subseteq L^{1/2}$. The power series $a_4,a_6$ have integer coefficients, so
$a_i(\mathsf q)^2=a_i(\mathsf q^2)$ in characteristic $2$: the Frobenius twist of $E_{\mathsf q}$ is $E_{\mathsf q^2}$,
and the Frobenius $E_{\mathsf q}\to E_{\mathsf q^2}$ is $u\mapsto u^2$. The Verschiebung
$V\colon E_{\mathsf q^2}\to E_{\mathsf q}$ is induced by the identity of $L^\times$: every $u\in L^\times$ is the square
of some $u_1$ in the local field $L_1=L^{1/2}$, and on $L_1$-points $V(u)=V(\mathrm{Frob}(u_1))=[2]u_1=u_1^2=u$.

If $\chi_{a_2}=0$, then $A\simeq E_{\mathsf q}$ and $V\colon F^\times/\mathsf q^{2\Z}\to F^\times/\mathsf q^\Z$ is surjective, so
$D(A)\simeq A(F)/V(A^{(2)}(F))=0$. Let $\chi_{a_2}\neq0$, let $L$ be the quadratic extension cut out by $\chi_{a_2}$ and
$\sigma$ its non-trivial automorphism. Then $A(F)=\{P\in E_{\mathsf q}(L):\sigma P=-P\}$, and $A^{(2)}$ is the twist of
$E_{\mathsf q^2}$ by the same character ($a_2^2\equiv a_2$ modulo $\wp(F)$), so
\[
A(F)\simeq\{u\in L^\times/\mathsf q^\Z:\ u\,\sigma(u)\in\mathsf q^\Z\},\qquad
A^{(2)}(F)\simeq\{u\in L^\times/\mathsf q^{2\Z}:\ u\,\sigma(u)\in\mathsf q^{2\Z}\},
\]
and $V$ is induced by the identity. For $u$ in the first group write $\Nm_{L/F}(u)=\mathsf q^m$; the parity of $m$
depends only on the class of $u$, and $u\mapsto m\bmod2$ is a homomorphism $A(F)\to\Z/2$ with kernel $V(A^{(2)}(F))$. It
is surjective if and only if $\mathsf q\in\Nm_{L/F}(L^\times)$ (if $m$ is odd, then $u\mathsf q^{(1-m)/2}$ has norm
$\mathsf q$), i.e.\ if and only if $[a_2,\mathsf q)_F=0$ (Fact~\ref{fact:norm}). As $b\in\mathsf qF^{\times2}$, this
gives $\dim D(A)=1-[a_2,b)_F$.

If $\chi_{a_2}$ is ramified, the representation of the Weil--Deligne group attached to $A$ is the twist by the
ramified character $\chi_{a_2}$ of that of $E_{\mathsf q}$; so $A$ has additive, potentially multiplicative reduction,
and $w(A/F)=1$. Indeed, for a Weil--Deligne representation $\sigma\otimes\mathrm{sp}(2)$ one has
$w(\sigma\otimes\mathrm{sp}(2))=w(\sigma)^2\det(-\mathrm{Frob}\mid\sigma^{I})$ \cite{Ro}; for $\sigma=\chi_{a_2}$ ramified,
$\sigma^I=0$, so $w(A/F)=w(\chi_{a_2})^2=w(\chi_{a_2})w(\chi_{a_2}^{-1})=\chi_{a_2}(-1)=1$, because $-1=1$ in $F$. Both
identities come from Tate's local functional equation and Deligne's definition of local constants, which apply to every
non-archimedean local field \cite{Ro}. Hence
$\epsilon(A)=(-1)^{(1-[a_2,b)_F)-1+[a_2,b)_F}=1$.
\end{proof}

For unramified $\chi_{a_2}$ the proposition recovers the computation of $D_v$ in Lemma~\ref{lem:mult}.

\subsubsection{Local constancy}
\begin{lemma}\label{lem:constancy}
Let $y^2+\alpha_1xy+\alpha_3y=x^3+\alpha_2x^2+\alpha_4x+\alpha_6$ be an integral Weierstrass equation of an ordinary
elliptic curve $A$ over $F$. There is an integer $N\geq1$ such that for all $\alpha'_i\in\mathcal O_F$ with
$\alpha'_i\equiv\alpha_i\pmod{\mathfrak m^N}$ the equation with coefficients $\alpha'_i$ defines an ordinary elliptic
curve $A'$ over $F$ with the same Kodaira type as $A$ and with $\epsilon(A')=\epsilon(A)$.
\end{lemma}

\begin{proof}
By \cite[Prop.~3.2]{CI}, applied with $K'=K=F$ and the identity of $\mathcal O_F/\mathfrak m^N$ (this case is
allowed, cf.\ \cite[Rem.~3.3]{CI}), for $N$ large the
curves $A$ and $A'$ have the same root number \cite[Prop.~3.2(ix)]{CI}, the same Kodaira type (v), the same valuation
of the minimal discriminant (ii), and the same Tamagawa number, the groups of rational points of the component groups
of the N\'eron models being isomorphic over the residue field (iii). The Frobenius twists $A^{(2)}$ and $A'^{(2)}$ are given by the equations
with coefficients $\alpha_i^2$ and $\alpha_i'^2$, which are congruent modulo $\mathfrak m^{2N}$; we take $N$ large
enough for the same conclusions to hold for them. As $A$ is ordinary, $\alpha_1\neq0$, and the model
$y^2+xy=x^3+a_2x^2+b$ of $A$ has $a_2=\alpha_2/\alpha_1^2+\alpha_3/\alpha_1^3$ and $b=\Delta/\alpha_1^{12}$
(proof of Lemma~\ref{lem:good}), and similarly for $A'$. For $N$ large, $\alpha_1'\neq0$, $a_2'-a_2\in\mathfrak m\subseteq\wp(F)$
and $b'/b=1+\varepsilon$ with $v(\varepsilon)\geq\max(1,1-v(a_2))$; then $[a_2',b')_F=[a_2,b)_F+[a_2,1+\varepsilon)_F$ and
$[a_2,1+\varepsilon)_F=\Tr\Res\bigl(a_2\,d\varepsilon/(1+\varepsilon)\bigr)=0$ (Fact~\ref{fact:schmid}). Finally, by
Lemma~\ref{lem:Dsize} and its proof,
\[
\dim D(A)-1=\ord_2\frac{c(A)}{c(A^{(2)})}+f\,v(\vartheta),\qquad
v(\vartheta)=\frac{v(\Delta_{\min}(A^{(2)}))-v(\Delta_{\min}(A))-v(b)}{12},
\]
because the N\'eron differentials of $A^{(2)}$ and $A$ are $e\,dx/x$ and $e'\,dx/x$ in the models \eqref{eq:E} with
$12v(e)=v(b^2)-v(\Delta_{\min}(A^{(2)}))$ and $12v(e')=v(b)-v(\Delta_{\min}(A))$. So $\dim D(A')=\dim D(A)$.
\end{proof}

\subsubsection{Globalisation}
\begin{lemma}\label{lem:globalise}
Let $F=\F_q((t))$, let $\alpha_1,\dots,\alpha_6\in\F_q[[t]]$ be the coefficients of an integral Weierstrass equation of
an ordinary elliptic curve over $F$, and let $N\geq1$. There are $a_1,\dots,a_6\in\F_q[t]$ with
$a_i\equiv\alpha_i\pmod{t^N}$ such that $y^2+a_1xy+a_3y=x^3+a_2x^2+a_4x+a_6$ is an ordinary elliptic curve over $\F_q(t)$
with semistable reduction at every place other than $t=0$.
\end{lemma}

\begin{proof}
In characteristic $2$ a Weierstrass equation has $c_4=a_1^4$ and
\begin{equation}\label{eq:disc2}
\Delta=a_3^4+a_1^3a_3^3+a_1^4\bigl(a_1^2a_6+a_1a_3a_4+a_2a_3^2+a_4^2\bigr)
\end{equation}
\cite[III.1]{Si}. We may assume $N>v(\alpha_1)=:k$ and $N>v(\Delta(\alpha))$. Let $a_1$ be the sum of $t^N$ and of the
truncation of $\alpha_1$ modulo $t^N$; then $a_1=t^kg$ with $g\in\F_q[t]$ monic of degree $N-k$ and $g(0)\neq0$. By
the Chinese remainder theorem there is $a_3\in\F_q[t]$ of degree less than $2N-k$ with $a_3\equiv\alpha_3\pmod{t^N}$
and $a_3\equiv1\pmod g$. Let $a_2,a_4,a_6$ be the truncations of $\alpha_2,\alpha_4,\alpha_6$ modulo $t^N$. Then
$\Delta(a)\equiv\Delta(\alpha)\not\equiv0\pmod{t^N}$, so the equation defines an elliptic curve $\tilde A$ over
$\F_q(t)$, which is ordinary because $a_1\neq0$.

Let $w\neq(t)$ be a finite place, with monic irreducible polynomial $\pi_w$. If $\pi_w\nmid a_1$, then $c_4$ is a unit
at $w$, and the reduction of the (integral) equation is smooth or nodal \cite[VII.5.1]{Si}. If $\pi_w\mid a_1$, then
$\pi_w\mid g$, so $a_3\equiv1$ and $\Delta\equiv a_3^4\equiv1$ modulo $\pi_w$ by \eqref{eq:disc2}: good
(supersingular) reduction. At
$\infty$, with $s=1/t$, the equation with coefficients $s^{iN}a_i(1/s)\in\F_q[s]$ is integral (as $\deg a_i\leq iN$) and
isomorphic to the given one, and its $c_4$ has the value $1$ at $s=0$, $a_1$ being monic of degree $N$; as the
other coefficients vanish at $s=0$, the reduction is the nodal cubic $y^2+xy=x^3$. So $\tilde A$ has multiplicative
reduction at $\infty$; in particular it is not isotrivial.
\end{proof}

\subsubsection{Proof of Theorem~\ref{thm:local}, and the argument of Remark~\ref{thm:additive}}
Theorem~\ref{thm:local} follows from Lemmas~\ref{lem:good} and~\ref{lem:mult} and Proposition~\ref{prop:potmult}: an
ordinary curve with potentially multiplicative reduction has $v(b)>0$, and its reduction is multiplicative if
$\chi_{a_2}$ is unramified. \qed

\smallskip
We now give the argument of Remark~\ref{thm:additive}; it uses Lemmas~\ref{lem:constancy} and~\ref{lem:globalise},
which are unconditional, and \cite[Thm.~1.1.5]{CF}. Let $A$ have additive reduction; we may assume $F=\F_q((t))$. Choose a minimal Weierstrass equation of $A$, let
$N$ be as in Lemma~\ref{lem:constancy}, and let $\tilde A$ over $K=\F_q(t)$ be the curve of Lemma~\ref{lem:globalise}.
At the place $v_0$ given by $t=0$, the curve $\tilde A$ has the Kodaira type of $A$ and
$\epsilon(\tilde A_{K_{v_0}})=\epsilon(A)$ (Lemma~\ref{lem:constancy}); so $v_0$ is the only place of additive reduction
of $\tilde A$, and Corollary~\ref{cor:delta} gives
\[
\epsilon(A)=(-1)^{\delta(\tilde A)} .
\]
Let $X\to\mathbb P^1_{\F_q}$ be the minimal regular model of $\tilde A$; it is a smooth, proper, geometrically
connected surface over $\F_q$, with a section, and $\Br(X)\simeq\Sha(\tilde A/K)$ \cite[\S4]{Gro}, \cite[\S3]{Ta66},
\cite{Ul2}. By \cite[Thm.~1.1.5]{CF} the quotient of $\Br(X)[2^\infty]$ by its divisible subgroup is a finite group
with a non-degenerate alternating form. Such a group is isomorphic to $G\times G$ for some finite abelian group $G$;
hence $\delta(\tilde A)$ is even, and $\epsilon(A)=1$.

\begin{remark}
The argument shows that the identity of Remark~\ref{thm:additive} for all $F$ is equivalent to the evenness of
$\delta(\tilde A)$ for all ordinary $\tilde A$ over the fields $\F_q(t)$ with at most one place of additive reduction.
In Propositions~\ref{prop:exnonsplit} and~\ref{prop:exsplit} the curve $E'$ has additive, potentially good reduction only
at $\infty$, of types $\mathrm I_0^*$ and $\mathrm{III}$, and the identity holds there: in the first case
$w(E'/K)=-1$ ($L$ has a simple zero at $T=1/8$) and $w_{v_b}=-1$, so $w_\infty=+1=(-1)^{\dim D_\infty-1}$; in the second
$w(E'/K)=-1$ and the four places of type $\mathrm I_2$ contribute $(-1)^4$, so $w_\infty=-1=(-1)^{8-1}$; in both cases
$[a_2,b)_\infty=0$ as $a_2=0$.
\end{remark}

\begin{remark}[numerical verification]\label{rem:numerics}
Independently of \cite{CF}, we tested the identity of Remark~\ref{thm:additive} on the curves
$A\colon y^2+xy=x^3+a_2x^2+b$ over $K=\F_2(t)$ with
$a_2\in\{0,\,1,\,t,\,t+1,\,t^3,\,t^3+t,\,t^3+t+1,\,t^5,\,t^5+t^3+t\}$ and $b\in\F_2[t]$ non-zero of degree at most $7$.
These curves are semistable at the finite places, and all but two of them have additive, potentially good reduction at
$\infty$. As $\deg L(A/K,T)=\deg\mathfrak n-4$ for non-constant $A$ over $\F_2(t)$ \cite[Lect.~1, \S12]{Ul2}, the condition
$\deg\mathfrak n\leq13$ bounds the number of point counts needed; for the $1017$ curves with $\deg L(A/K,T)\leq9$ we computed $L(A/K,T)$ by counting points,
read off $w(A/K)$ from the functional equation, obtained $w_\infty(A)$ by dividing by the root numbers at the finite
places (\S\ref{subsec:rootnumbers}), computed $\dim D_\infty$ from Lemma~\ref{lem:Dsize} with Tate's algorithm, and
$[a_2,b)_\infty$ as a residue. The identity $w_\infty(A)=(-1)^{\dim D_\infty-1+[a_2,b)_\infty}$ holds in all $1017$
cases; the Kodaira types $\mathrm{II}$, $\mathrm{III}$, $\mathrm{IV}^*$, $\mathrm I_n^*$ ($0\leq n\leq20$, $n\neq5$),
$\mathrm{III}^*$ and $\mathrm{II}^*$ occur, with both signs of $w_\infty$ and both values of the symbol. For all
$2295$ curves of the family the value of $\dim D_\infty$ given by Lemma~\ref{lem:Dsize} agrees with a direct search for
local points. The computations use Python~3 with the NumPy and \texttt{galois} packages; the programs are provided as
supplementary material.
\end{remark}

\subsubsection{Curves with finite Tate--Shafarevich group}
When $\Sha$ is finite, the parity of the rank can be read off directly from the Birch and Swinnerton-Dyer formula for
$A$ and $A^{(2)}$, by the isogeny-invariance argument of Cassels as used by Dokchitser and Dokchitser \cite{DD}. The rank parity and
the identity $\prod_v\epsilon_v=1$ below need neither Theorem~\ref{thm:radical} nor \cite{TY,CF}; only the last
assertion uses Theorem~\ref{thm:parity}.

\begin{proposition}\label{prop:finSha}
Let $A/K$ be an ordinary elliptic curve with $\Sha(A/K)$ finite. Then
$\rk A(K)\equiv\sum_v(\dim D_v(A)-1)\pmod2$, the local identities $\epsilon_v(A):=w_v(A)(-1)^{\dim D_v(A)-1+[a_2,b)_v}$
satisfy $\prod_v\epsilon_v(A)=1$, and $\dim\Sha(A/K)[2]$ is even.
\end{proposition}

\begin{proof}
Put $E'=A$, $E=A^{(2)}$, let $r$ be their common rank, $g$ the genus of $K$, and $C(X)=\prod_vc_v(X)$. As $E$ and $E'$
are isogenous, $\Sha(E/K)$ is finite too, and the Birch and Swinnerton-Dyer formula holds for both, with the same
$L$-function \cite{KT}, \cite{Ul2}: the quantity
\[\operatorname{Reg}(X)\,|\Sha(X/K)|\,C(X)\,q^{1-g-\deg\omega_X}/|X(K)_{\mathrm{tors}}|^2\]
is the same for $X=E$ and $X=E'$.
Canonical heights satisfy $\hat h_{E'}(\varphi P)=2\hat h_E(P)$; if $i$ is the index of $\varphi(E(K))$ modulo torsion
in $E'(K)$ modulo torsion, computing the determinant of the height pairing on $\varphi(E(K))$ in two ways gives
$2^r\operatorname{Reg}(E)=i^2\operatorname{Reg}(E')$. The orders of $\Sha(E/K)$ and $\Sha(E'/K)$ are squares, by the
theorem of Liu, Lorenzini and Raynaud on Brauer groups of surfaces \cite{LLR,LLRc,Gei}, since $\Sha$ is the Brauer group
of the minimal regular model of the curve, which has a section \cite[\S4]{Gro}. Taking $2$-adic valuations, and
writing $q=2^n$,
\[
r\equiv\ord_2\frac{C(E')}{C(E)}+n\bigl(\deg\omega_E-\deg\omega_{E'}\bigr)=\sum_v(\dim D_v-1)\pmod2,
\]
the last equality by Lemma~\ref{lem:Dsize}. As $\Sha$ is finite, $r=\ord_{s=1}L(A/K,s)$ \cite{KT}, so
$(-1)^r=w(A/K)=\prod_vw_v(A)$ (\S\ref{subsec:rootnumbers}); with $\sum_v[a_2,b)_v=0$ this gives $\prod_v\epsilon_v(A)=1$.
Finally $\rk A(K)+\dim\Sha(A/K)[2]\equiv\sum_v(\dim D_v-1)$ by Theorem~\ref{thm:parity}, so $\dim\Sha(A/K)[2]$ is even.
\end{proof}

If $\dim\Sha(A/K)[2]\leq1$, then $\delta(A)$ is even without any finiteness hypothesis: either $\Sha(A/K)[2^\infty]$ is
finite, so that $\Sha(A/K)$ is finite \cite{KT} and $\Sha(A/K)[2]=0$ by the proposition, or it is $\Q_2/\Z_2$ and
$\delta(A)=0$. Thus the identity of Remark~\ref{thm:additive} holds for $A$ over $F$ as soon as the curve $\tilde A$ of
\S\ref{subsec:additive} can be chosen with $\Sha(\tilde A/K)$ finite or with $\dim\Sha(\tilde A/K)[2]\leq1$.

\begin{theorem}\label{thm:finSha}
Let $F=\F_q((s))$ and $A\colon y^2+xy=x^3+a_2x^2+b$ an ordinary elliptic curve over $F$. Then
$w(A/F)=(-1)^{\dim D(A)-1+[a_2,b)_F}$ in each of the following cases:
\begin{enumerate}[label=\textup{(\arabic*)},leftmargin=*]
\item $\chi_{a_2}$ is unramified, $v(b)=-2^em$ with $m$ odd and positive and $e\geq0$, and $b\in F^{2^e}$; in
particular whenever $\chi_{a_2}$ is unramified and $v(b)$ is odd and negative;
\item $b\in\F_q^\times$, and $a_2\in F$ is arbitrary.
\end{enumerate}
The reduction is potentially good in both cases. It is additive in case \textup{(1)} when $12\nmid v(b)$
\textup{(}e.g.\ when $e\leq1$\textup{)}, and in case \textup{(2)} when $\chi_{a_2}$ is ramified.
\end{theorem}

\begin{proof}
The reduction is potentially good because $j(A)=1/b$ and $v(b)\leq0$. In case (1), $j(A)\in\mathfrak m_F$, so good
reduction would be supersingular, and then $b=\Delta/a_1^{12}$ for a minimal model (proof of Lemma~\ref{lem:good}), so
$12\mid v(b)$. In case (2), $A$ is the twist by $\chi_{a_2}$ of the curve $y^2+xy=x^3+b$ with good reduction, so its
reduction is additive if $\chi_{a_2}$ is ramified. In both cases we construct a curve $\tilde A$ over
$K=\F_q(t)$, with a place $v_0$ such that $K_{v_0}\simeq F$ and $\tilde A_{K_{v_0}}\simeq A$, with semistable reduction
at all other places, and with $\Sha(\tilde A/K)$ finite. Then $\epsilon_v(\tilde A)=1$ for $v\neq v_0$ by
Theorem~\ref{thm:local}, so Proposition~\ref{prop:finSha} gives $\epsilon_{v_0}(\tilde A)=1$, which is the assertion.
Finiteness of $\Sha$ will be checked after a finite constant extension $K'=\F_{q'}(t)$: the kernel of
$\Sha(\tilde A/K)\to\Sha(\tilde A/K')$ lies in the finite group $H^1(\Gal(K'/K),\tilde A(K'))$.

(1) Write $b=cs^{-2^em}u$ with $c\in\F_q^\times$ and $u$ a $1$-unit. As $b$, $c$ and $s^{-2^em}$ are $2^e$-th powers, so is
$u$, and comparing residues, $u$ is the $2^e$-th power of a $1$-unit. The group of $1$-units is a $\Z_2$-module, so
uniquely divisible by the odd integer $m$; hence $u=u_1^{-2^em}$ for a $1$-unit $u_1$, and replacing the uniformiser $s$
by $su_1$ we may assume $b=cs^{-2^em}$. As $\mathfrak m_F\subseteq\wp(F)$ and $\chi_{a_2}$ is unramified, we may assume
$a_2=a\in\F_q$ (Remark~\ref{rem:uniquemodel}). Let $\tilde A\colon y^2+xy=x^3+ax^2+ct^{2^em}$ over $K$, and
$v_0=\infty$, $t=1/s$. At $t=0$ the reduction is multiplicative ($c_4=1$), at the other finite places it is good.
Over $K'=\F_{q'}(t)$, where $\F_{q'}$ contains $z$ with $z^2+z=a$ and $\gamma$ with $\gamma^{2^em}=c$, the substitution
$y\mapsto y+zx$, followed by the transport of structure along the $\F_{q'}$-automorphism $t\mapsto\gamma^{-1}t$ of $K'$,
transforms $\tilde A$ into $y^2+xy=x^3+t^{2^em}$, the $e$-th Frobenius twist
of Ulmer's curve $E_m\colon y^2+xy=x^3-t^m$, to which it is isogenous. For every $m$ prime to the characteristic and every finite constant field, Ulmer proves the Birch and
Swinnerton-Dyer conjecture for $E_m$, via the Tate conjecture for its elliptic surface, which is dominated by a Fermat
surface \cite[Prop.~6.4]{Ul02}; as $m$ is odd this applies, and $\rk E_m(K')=\ord_{s=1}L(E_m/K',s)$; by
\cite{KT} (see also \cite{Mi75}) $\Sha(E_m/K')$ is then finite, including its $2$-primary part. Hence so is
$\Sha(\tilde A/K')$.

(2) Let $v_0$ be the place $t=0$ and identify $F$ with $K_{v_0}=\F_q((t))$. As $t\F_q[[t]]\subseteq\wp(F)$, there is
$\tilde a\in\F_q[t^{-1}]$ with $a_2-\tilde a\in\wp(F)$, and $\tilde A\colon y^2+xy=x^3+\tilde ax^2+b$ satisfies
$\tilde A_{K_{v_0}}\simeq A$ (Remark~\ref{rem:uniquemodel}). At every other place, including $\infty$, the model is
integral with unit discriminant $b$, so the reduction is good. Over $K'=K(z)$ with $z^2+z=\tilde a$, the substitution
$y\mapsto y+zx$ identifies $\tilde A$ with the constant curve $A_b\colon y^2+xy=x^3+b$. For a constant curve the rank equals the
order of vanishing of the $L$-function, by Tate's theorem on homomorphisms of abelian varieties over finite fields
\cite{TaEnd}, \cite{Mi68}; so $\Sha(A_b/K')$ is finite, including its $2$-primary part, by \cite{KT}. As $K'/K$ is finite Galois, $\Sha(\tilde A/K)$ is finite.
\end{proof}

For example, $y^2+xy=x^3+t$ over $\F_2(t)$ falls under case (1) at $\infty$, where it has type $\mathrm{II}^*$; its
$L$-function is $1$, and $w_\infty=-1=(-1)^{\dim D_\infty-1}$ with $\dim D_\infty=0$. Many of the curves of
Remark~\ref{rem:numerics} (those with constant $a_2$ and $\deg b$ odd, and those with $b=1$) fall under
Theorem~\ref{thm:finSha}.

\section{Points of prime order}\label{sec:ell}
In this section $K$ is a global field of any characteristic $p\geq0$, $\ell$ is a prime ($\ell=p$ is allowed), $E/K$ is
an elliptic curve with a point $T\in E(K)$ of order $\ell$, $\varphi\colon E\to E'=E/\langle T\rangle$, and we identify
$\langle T\rangle$ with $\Z/\ell$ by $T\mapsto1$. The isogeny $\varphi$ is \'etale. For $R\in E'(K)$ and
$Q\in E(\Ks)$ with $\varphi(Q)=R$, the Kummer class $\delta_\varphi(R)\in H^1(K,\Z/\ell)$ is $\sigma\mapsto\sigma Q-Q$,
and $\Sel^\varphi(E/K)\subseteq H^1(K,\Z/\ell)$; $\iota\colon H^1(K,\Z/\ell)\to H^1(K,E)$ is induced by
$\Z/\ell\simeq\langle T\rangle\subseteq E$. We show that the proof of Theorem~\ref{thm:BM} extends to this situation,
the role of the norm equation being played by one norm equation in the cyclic extension cut out by a Selmer class.

\emph{Conventions.} For $\ell=2$ no sign matters; for odd $\ell$ the formulas below depend on the following
choices, which we fix. (i) The class of a divisor $\sum n_i(Q_i)$ of degree $0$ in $\Pic^0(X_{\chi,\s})=E(\Ks)$ is
$\sum n_iQ_i$, i.e.\ $E$ is identified with its dual by the inverse of the sum map, $P\mapsto[(P)-(O)]$. The
polarisation $\phi_{\mathcal L}\colon P\mapsto[t_P^*\mathcal L\otimes\mathcal L^{-1}]$ attached to $\mathcal L=\mathcal O((O))$
sends $P$ to $[(-P)-(O)]$; with that identification every pairing in this
section changes sign. (ii) Invariants $\inv_v\colon\Br K_v\to\Q/\Z$ are normalised as in \cite[XIII, \S3]{Se}. (iii)
The symbol $(\psi,u)$ is defined by the cocycle below; it equals $\delta\psi\cup u$ exactly, where
$\delta\colon H^1(K,\Z/\ell)\to H^2(K,\Z)$ is the connecting map of $0\to\Z\xrightarrow{\ell}\Z\to\Z/\ell\to0$, because
$\delta\psi(\sigma,\tau)=(\tilde\psi(\sigma)+\tilde\psi(\tau)-\tilde\psi(\sigma\tau))/\ell$. (iv) Schmid's formula is used in
the normalisation of \cite[XIV, \S5]{Se}; other normalisations of the Artin--Schreier character differ from it by a sign
for odd $p$, which changes only the overall sign.

For $\psi\in H^1(K,\Z/\ell)$ let $\tilde\psi(\sigma)\in\{0,\dots,\ell-1\}$ be the representative of $\psi(\sigma)$, and
for $u\in K^\times$ let $(\psi,u)\in\Br K$ be the class of the $2$-cocycle
$(\sigma,\tau)\mapsto u^{\lfloor(\tilde\psi(\sigma)+\tilde\psi(\tau))/\ell\rfloor}$; this is the cyclic algebra
$\delta\psi\cup u$ \cite[XIV, \S1]{Se}, it is bi-additive, it vanishes if and only if $u$ is a norm from the extension
cut out by $\psi$ when $\psi\neq0$ \cite[XIV, \S1, Prop.~2]{Se}, and $\sum_v\inv_v(\psi_v,u)=0$
(Fact~\ref{fact:ABHN}). For $\ell=2$ it is \eqref{eq:cycliccocycle}. If $\ell=p$ and $\psi$ is the Artin--Schreier
character of $\alpha\in K_v$, then $\inv_v(\psi,u)=\frac1p\Tr\Res(\alpha\,du/u)$ \cite{Sch}, \cite[XIV, \S5,
Prop.~15]{Se}.

For $\chi\in H^1(K,\Z/\ell)$ let $X_\chi$ be the torsor of $\iota(\chi)$: $X_{\chi,\s}=E_\s$, with the twisted action
$\sigma^*Q=\sigma Q+\chi(\sigma)T$ on points, and the induced actions on divisors and on functions. As $T$ is rational,
$\sigma^*$ commutes with translations by multiples of $T$. For a divisor $D$ and $j\in\Z$ we write $D+jT$ for the
translate of $D$ by $jT$; divisor classes are taken as in convention (i).

\begin{lemma}\label{lem:cycliclift}
Let $\chi,\chi'\in H^1(K,\Z/\ell)$, $\psi=\chi-\chi'$ \textup{(}for $\ell=2$ this is the character $\chi+\chi'$ of
\S\ref{subsec:proofBM}\textup{)}, and let $\Delta\in\Div^0(X_{\chi,\s})$ have class $-T$ and satisfy
$\sigma^*\Delta=\Delta+\tilde\psi(\sigma)T$ for all $\sigma\in G_K$. Then $N(\Delta)=\sum_{j=0}^{\ell-1}(\Delta+jT)$ is the
divisor of a function $f\in K(X_\chi)^\times$. If $\iota(\chi),\iota(\chi')\in\Sha(E/K)$, then for all
$P_v\in X_\chi(K_v)$ outside the support of $N(\Delta)$,
\[
\langle\iota(\chi),\iota(\chi')\rangle=\sum_v\inv_v\bigl(\psi_v,f(P_v)\bigr).
\]
\end{lemma}

\begin{proof}
The divisor $N(\Delta)$ is invariant under all $\sigma^*$ and has class $\ell\cdot(-T)=O$, so it is the divisor of a
function $f_1\in\Ks(X_\chi)^\times$ with $\sigma^*f_1/f_1\in\Ks{}^\times$; by Hilbert's Theorem~90, $f_1$ can be
multiplied by a constant so as to become $G_K$-invariant. Put
\[
a''_\sigma=\sum_{j=0}^{\tilde\psi(\sigma)-1}(\Delta+jT)+\bigl(\sigma^*(O)-(O)\bigr).
\]
Since $\sigma^*(O)=(\chi(\sigma)T)$, the class of $a''_\sigma$ is $-\tilde\psi(\sigma)T+\chi(\sigma)T=\chi'(\sigma)T$, so
$a''$ lifts the cocycle $\sigma\mapsto\chi'(\sigma)T$ representing $\iota(\chi')$. The second summand is a coboundary.
For the first, $\sigma^*(\Delta+jT)=\Delta+(\tilde\psi(\sigma)+j)T$, so $\sigma^*a''_\tau+a''_\sigma-a''_{\sigma\tau}$ is
the sum of the translates $\Delta+jT$ for $0\leq j<\tilde\psi(\sigma)+\tilde\psi(\tau)$ minus the sum for
$0\leq j<\tilde\psi(\sigma\tau)$. As $\Delta+(j+\ell)T=\Delta+jT$ and
$\tilde\psi(\sigma\tau)=\tilde\psi(\sigma)+\tilde\psi(\tau)-\ell\lfloor(\tilde\psi(\sigma)+\tilde\psi(\tau))/\ell\rfloor$,
\[
(da'')(\sigma,\tau)=\Bigl\lfloor\frac{\tilde\psi(\sigma)+\tilde\psi(\tau)}\ell\Bigr\rfloor N(\Delta)=\dv\bigl(f^{\lfloor(\tilde\psi(\sigma)+\tilde\psi(\tau))/\ell\rfloor}\bigr),
\]
and $f'(\sigma,\tau)=f^{\lfloor(\tilde\psi(\sigma)+\tilde\psi(\tau))/\ell\rfloor}$ is a $2$-cocycle with values in
$\Ks(X_\chi)^\times$, as $f$ is invariant. By Definition~\ref{def:PS} and Lemma~\ref{lem:ev}(b) the local terms are
the classes $(\psi_v,f(P_v))$.
\end{proof}

\begin{theorem}\label{thm:ell}
Let $\chi'\in\Sel^\varphi(E/K)$, let $X'=X_{\chi'}$, and let $\pi'\colon X'\to E'$ be the $K$-morphism given by $\varphi$
on points. Let $\mathfrak b'\in\Div^0(X'_\s)$ be $K$-rational, of class $-T$, and let $g\in K(E')^\times$ with
$\dv g=\pi'_*\mathfrak b'$. Then for all $\chi\in\Sel^\varphi(E/K)$ and all $R_v\in E'(K_v)$ with
$\delta_\varphi(R_v)=\chi_v$, outside the support of $\dv g$,
\[
\langle\iota(\chi),\iota(\chi')\rangle=\sum_v\inv_v\bigl(\chi_v+\chi'_v,\ g(R_v)\bigr).
\]
\end{theorem}

\begin{proof}
On $E_\s$ the actions of $X_{-\chi}$ and $X'$ differ by $\sigma^*_{X_{-\chi}}Q=\sigma^*_{X'}Q-(\chi+\chi')(\sigma)T$. As
$\mathfrak b'$ is invariant under the action of $X'$, the divisor $\Delta=\mathfrak b'$ on $X_{-\chi}$ satisfies the
hypothesis of Lemma~\ref{lem:cycliclift} for the pair $(-\chi,\chi')$, with $\psi=-\chi-\chi'$. The map
$\pi\colon X_{-\chi}\to E'$ given by $\varphi$ is a $K$-morphism, and as $\varphi$ is \'etale with kernel
$\langle T\rangle$, $\pi^*\pi'_*\mathfrak b'=N(\mathfrak b')$; so one may take $f=\pi^*g$. The points of $X_{-\chi}(K_v)$
are the $Q\in E(\Kvs)$ with $\sigma Q-Q=\chi(\sigma)T$ for $\sigma\in G_v$, i.e.\ the points above the $R_v\in E'(K_v)$ with
$\delta_\varphi(R_v)=\chi_v$. Lemma~\ref{lem:cycliclift} gives
$\langle\iota(-\chi),\iota(\chi')\rangle=\sum_v\inv_v(-\chi_v-\chi'_v,g(R_v))$, and the theorem follows by bilinearity
(Corollary~\ref{cor:bilinear}).
\end{proof}

It remains to produce $\mathfrak b'$. Choose $u\in K(E)^\times$ with $\dv u=2(T)-(O)-(2T)$ (for $\ell=2$ this is
$2(T)-2(O)$), and put
\[
c_u=\prod_{j=0}^{\ell-1}u\circ\tau_{jT}\in K^\times,
\]
where $\tau_{jT}$ is the translation by $jT$; $c_u$ is constant, since its divisor is $0$, and $K$-rational. For
$\ell=2$ and $E\colon y^2+xy=x^3+a_2x^2+b^2$ one may take $u=x$, and then $c_u=x\cdot(b/x)=b$ by (E2).

\begin{proposition}\label{prop:bprime}
Let $\chi'\in\Sel^\varphi(E/K)$, $\chi'\neq0$, let $L'$ be the cyclic extension of degree $\ell$ cut out by $\chi'$, and
$\rho\in G_K$ with $\chi'(\rho)=1$.
\begin{enumerate}[label=\textup{(\alph*)},leftmargin=*]
\item $c_u$ is a norm from $L'$.
\item Let $\lambda\in L'$ with $\Nm_{L'/K}\lambda=c_u^{-1}$, $W=\lambda u\in L'(X')$, and
$F=\sum_{k=0}^{\ell-1}\prod_{m=0}^{k-1}\rho^{m*}W$. Then $F\neq0$, $\rho^*F=F/W$, and $\mathfrak b'=\dv F+(O)-(T)$ is a
$K$-rational divisor of degree $0$ and class $-T$ on $X'$.
\item $\mathfrak b'+\pi'^*(O)$ is effective of degree $\ell$; hence $g$ in Theorem~\textup{\ref{thm:ell}} lies in the
Riemann--Roch space $L(\ell(O))$ of $E'$.
\end{enumerate}
\end{proposition}

\begin{proof}
For $\sigma\in G_K$ and $F\in\Ks(X')$ we have $\sigma^*F=\sigma(F)\circ\tau_{-\tilde\chi'(\sigma)T}$, where $\sigma(F)$ is
obtained by applying $\sigma$ to the coefficients. A translate $u\circ\tau_{-mT}$ has divisor
$2((m+1)T)-(mT)-((m+2)T)$; writing $e_m=((m+1)T)-(mT)$, this is $e_m-e_{m+1}$.

(a) Let $D_0=(O)-(T)$, a divisor of class $-T$ on $X'$, and for $\sigma\in G_K$ with $k=\tilde\chi'(\sigma)$ put
$U_\sigma=\prod_{m=0}^{k-1}u\circ\tau_{-mT}$. Then $\dv U_\sigma=e_0-e_k=\sigma^*D_0-D_0$. As $u$ has coefficients in
$K$, $\sigma^*U_\tau=\prod_{m=k}^{k+\tilde\chi'(\tau)-1}u\circ\tau_{-mT}$, so
\[
\frac{\sigma^*U_\tau\cdot U_\sigma}{U_{\sigma\tau}}=c_u^{\lfloor(\tilde\chi'(\sigma)+\tilde\chi'(\tau))/\ell\rfloor},
\]
a cocycle of constants representing $(\chi',c_u)$. If $(\chi',c_u)=0$, i.e.\ $c_u^{\lfloor\cdot\rfloor}=d\mu$ with $\mu$ a
$1$-cochain of constants, then $\sigma\mapsto U_\sigma/\mu_\sigma$ is a $1$-cocycle with values in $\Ks(X')^\times$, hence
a coboundary $\sigma^*h/h$ by Hilbert's Theorem~90, and $D_0-\dv h$ is $K$-rational of class $-T$. Conversely, the
class $(\chi',c_u)$ is the image of the class $-T\in\Pic(X'_\s)^{G_K}$ under the map
$\Pic(X'_\s)^{G_K}\to\Br K$ of the Hochschild--Serre spectral sequence, so its image in $\Br X'$ vanishes. For every $v$,
$X'(K_v)\neq\emptyset$ because $\chi'$ is a Selmer class, so $\Br K_v\to\Br X'_{K_v}$ is injective
(Lemma~\ref{lem:ev}(a)) and $(\chi'_v,c_u)=0$. Hence $(\chi',c_u)=0$ (Fact~\ref{fact:ABHN}), and $c_u$ is a norm from
$L'$.

(b) We have $\rho^{m*}W=\rho^m(\lambda)\,u\circ\tau_{-mT}$, so $\prod_{m=0}^{\ell-1}\rho^{m*}W=\Nm_{L'/K}(\lambda)c_u=1$. The
$k$-th summand $\prod_{m<k}\rho^{m*}W$ of $F$ has divisor $e_0-e_k=(T)+(kT)-(O)-((k+1)T)$; the summand with
$k=\ell-1$ is the only one with a double pole at $(O)$, so $F\neq0$. Then
\[
\rho^*F=\sum_{k=0}^{\ell-1}\prod_{m=1}^{k}\rho^{m*}W=\frac1W\sum_{k=1}^{\ell}\prod_{m=0}^{k-1}\rho^{m*}W=\frac FW,
\]
since the term $k=\ell$ is $1$, the term $k=0$ of $F$. So $\rho^*\dv F=\dv F-\dv W=\dv F-e_0+e_1$, and
$\rho^*((O)-(T))=(T)-(2T)=e_0-e_1+(O)-(T)$. Hence $\rho^*\mathfrak b'=\mathfrak b'$. If $\sigma\in G_{L'}$, then $\chi'(\sigma)=0$ and
$\sigma^*F=\sigma(F)=F$, because $F$ has coefficients in $L'$; so $\sigma^*$ fixes $\mathfrak b'$. As $G_K$ is generated by $\rho$ and $G_{L'}$, $\mathfrak b'$ is
$K$-rational; its class is that of $(O)-(T)$, i.e.\ $-T$.

(c) By the computation of the divisors of the summands, $F\in L\bigl(2(O)+\sum_{j=2}^{\ell-1}(jT)\bigr)$ (for $\ell=2$:
$F\in L(2(O))$), so $\dv F+(O)-(T)\geq-\sum_{j=0}^{\ell-1}(jT)=-\pi'^*(O)$. Thus $\pi'_*\mathfrak b'=\pi'_*Z-\ell(O)$ with
$Z=\mathfrak b'+\pi'^*(O)$ effective of degree $\ell$, and $g\in L(\ell(O))$.
\end{proof}

\begin{remark}\label{rem:ellcases}
(1) For $\ell=2$ and $E$ as in \eqref{eq:E}, with $u=x$ and $c_u=b$, the divisor $\mathcal Z-(O)-(T)$ of the proof of
Theorem~\ref{thm:BM} (\S\ref{subsec:proofBM}) is $K$-rational on $X'$ and of class $T=-T$, and $g=x_{E'}+t_0$; so
Theorem~\ref{thm:ell} contains Theorem~\ref{thm:BM}. The norm equation of Proposition~\ref{prop:bprime} is that of
Lemma~\ref{lem:norm}: $\lambda=1/\beta_0$, and then $F=1+x/\beta_0$ vanishes on $\mathcal Z$.

(2) For $\ell=p=\operatorname{char}K$ and $E$ ordinary, $\varphi$ is the Verschiebung, and the formula of
Theorem~\ref{thm:ell} is explicit given a solution of the norm equation, the function $g\in L(p(O))$, and the Kummer map
$\delta_\varphi$; the local terms are residues, by Schmid's formula. The explicit form of $g$ for $\ell=3$, and examples
in characteristic $3$, will be treated elsewhere.

(3) Over number fields, Theorem~\ref{thm:ell} is a reformulation of known methods: for $\ell=2$ it is a form of the
formulas of Fisher \cite{Fi}, and for $\ell=3$ it is close to the method of van Beek and Fisher \cite{BF}, which also
uses one norm equation.
\end{remark}

\section{Further remarks}\label{sec:remarks}
\begin{enumerate}[leftmargin=*,label=\arabic*.,ref=\arabic*]
\item\label{it:flat} \emph{Comparison with flat duality.} On $\Sha(E/K)[\varphi]$ the homogeneous-space pairing coincides
with the pairing induced by Milne's flat duality over every sufficiently small open $U$ (Proposition~\ref{prop:milne}). This is not used elsewhere: the
non-divisibility statements of Corollary~\ref{cor:bilin} and Theorems~\ref{thm:example}(d) and~\ref{thm:nonsplit}(c)
follow from bilinearity and functoriality, and the non-degeneracy used for parity is Theorem~\ref{thm:radical}. We do
not know whether the two pairings agree on all of $\Sha(A/K)[2^\infty]$ (Remark~\ref{rem:milnescope}).
\item \emph{General isogenies.} Let $\rho\colon A\to B$ be an isogeny of elliptic curves over $K$ with finite flat
kernel. The method of the proof of Proposition~\ref{prop:F} reduces the analogue of Theorem~\ref{thm:W} for
$(\rho,\hat\rho)$ to a local compatibility: for $R\in B(K_v)$ and $\beta\in H^1_\fl(K_v,B[\hat\rho])$, the Brauer class
measuring the difference between a divisorial lift, in the sense of Definition~\ref{def:PS}, of the image of $\beta$ in
$H^1(K_v,B)$ and its translate by $R$ should be $\pm\,\delta_\rho(R)\cup\beta$, the cup product being taken with respect to
the Weil pairing $A[\rho]\times B[\hat\rho]\to\mu_{\deg\rho}$. For $\deg\rho$ prime to the characteristic this is an
instance of the compatibility between divisorial and cohomological descriptions of the local Tate pairing (cf.\
\cite[I, Rem.~3.5]{Mi}, \cite[Prop.~34]{PS}); for infinitesimal kernels we have not settled it. For
$(\varphi,\hat\varphi)$, Theorem~\ref{thm:A} and Corollary~\ref{cor:Aprime} prove the formula directly.
\item \emph{Brauer groups of surfaces.} The group $\Sha(E/K)$ is closely related to the Brauer group of the associated
elliptic surface, on whose $\ell$-primary part, $\ell\neq p$, Artin and Tate defined an analogous pairing; Feng
\cite{Fe} proved that this pairing is alternating, as conjectured by Tate. That the order of the Brauer group of a
surface over a finite field is a square when finite, including the $p$-part and $p=2$, is a theorem of Liu, Lorenzini
and Raynaud \cite{LLR}, with the corrections of \cite{LLRc,Gei}. It would be interesting to interpret the algebra
$(\chi+\chi',f)$ of Remark~\ref{rem:BMint}, and Theorem~\ref{thm:radical}, on the surface.
\end{enumerate}

\section*{Acknowledgements}
The author thanks the anonymous referees for their careful reading, and in particular for the observation behind
Proposition~\ref{prop:finSha} and the cases of Theorem~\ref{thm:finSha}, for the purity argument in
Lemma~\ref{lem:BMgood}, and for suggesting the argument of Appendix~\ref{app:milne}.

The programs used for the numerical verification described in Remark~\ref{rem:numerics}, which is not used in any
proof, are available from the author on request. The author used AI for literature search and to help
write these programs; the author reviewed the results and takes full responsibility for the content of this article.

\appendix
\section{Two proofs on flat cohomology}\label{app:flat}
Throughout, $F$ is a field of characteristic $2$ containing $a_2$ and $b$, and $E,E',\varphi,\hat\varphi$ are as in
\S\ref{sec:prelim}.

\subsection{Proof of Proposition~\ref{prop:kummerF}}
\begin{proof}
\emph{Gradings} (see \cite[Exp.~I, \S4.7]{SGA3}). An action of $\mu_2=\Spec F[\zeta]/(\zeta^2-1)$ on an affine $F$-scheme $\Spec A$ is a coaction
$\rho\colon A\to A\otimes F[\zeta]/(\zeta^2-1)$; writing $\rho(a)=a_0\otimes1+a_1\otimes\zeta$, coassociativity and the
counit axiom say that $a\mapsto a_i$ are complementary idempotent projections, and $\rho$ being an algebra
homomorphism says that $A=A_0\oplus A_1$ is a $\Z/2$-grading; conversely a grading defines an action, and
equivariant morphisms correspond to graded homomorphisms. In characteristic $2$ the map $D(a):=a_1$ is a derivation
with $D^2=D$ (indeed $(ab)_1=a_0b_1+a_1b_0=aD(b)+D(a)b$), and $A_i=\{a:D(a)=ia\}$. Writing $\zeta=1+z$,
$\rho(a)=a\otimes1+D(a)\otimes z$.

\emph{The derivation of the translation action.} The group $\ker\hat\varphi\simeq\mu_2$ acts on $E'$ by translation.
This action is the identity on the underlying topological space, so it is given by a coaction on
$\mathcal O(U)$ for every open $U$ and on the function field $F(E')$, i.e.\ by a derivation $D$ of $F(E')$ with
$D^2=D$ preserving $\mathcal O(U)$ for all $U$. Since $E'$ is commutative, the action commutes with translations by
points of $E'(\Fs)$, so (after extending scalars to $\Fs$) $D$ commutes with $\tau_Q^*$ for $Q\in E'(\Fs)$. The
derivations of $F(E')/F$ form a one-dimensional $F(E')$-vector space spanned by $D_\omega(f):=df/\omega$, where
$\omega=dx/x$ is the invariant differential of $E'$ ($\omega=dx/(2y+a_1x+a_3)$, \cite[III.5.1]{Si}). As $D_\omega$ is
translation invariant, $D=cD_\omega$ with $c$ a function invariant under all translations, hence constant:
$c\in F$. Now $D_\omega(x)=x$, and differentiating the equation of $E'$ gives $x\,dy=(x^2+y)\,dx$, i.e.\
$D_\omega(y)=y+x^2$. Since $D_\omega^2$ is a derivation (characteristic $2$) and $D_\omega^2(x)=x=D_\omega(x)$, and
derivations of $F(E')$ are determined by their value on the separating element $x$, we get $D_\omega^2=D_\omega$.
Hence $c^2D_\omega=D^2=D=cD_\omega$, so $c\in\{0,1\}$. If $c=0$ the action would be trivial; but for the universal
point $\zeta_0\in(\ker\hat\varphi)(F[z]/(z^2))$ one has $z(\zeta_0+O)=z(\zeta_0)=z\neq0$, so $\zeta_0+O\neq O$.
Therefore $D=D_\omega$:
\[
D(x)=x,\qquad D(y)=y+x^2 .
\]

\emph{The fibres.} The torsor $\delta_{\hat\varphi}(P)$ is the scheme-theoretic fibre $\hat\varphi^{-1}(P)$ with the
translation action of $\ker\hat\varphi$. For $P\neq O$ this fibre lies in the affine open $U=E'\smallsetminus\{O\}$,
and it is stable under the action, so its ideal is a graded ideal of $\mathcal O(U)=F[x,y]/(f_{E'})$ and its coordinate
ring $A$ is a graded quotient, of dimension $\deg\hat\varphi=2$, on which $D$ acts as above. Two $\mu_2$-torsors whose
graded coordinate rings are isomorphic as graded algebras are isomorphic, and the Kummer torsor of $u\in F^\times$ is
$\Spec F[e]/(e^2-u)$ with $e$ of weight one. Let $P=(X_P,Y_P)$. The fibre is defined by $x^2=X_P$ and $y^2+a_2x^2=Y_P$
together with $f_{E'}=0$; substituting $y^2=Y_P+a_2x^2$ into $f_{E'}=0$ gives $xy=x^3+b+Y_P=xX_P+Y_P+b$ in $A$.

If $X_P\neq0$, then $x$ is a unit in $A$ and $y=X_P+(Y_P+b)/x$, so $F[x]/(x^2-X_P)\to A$ is surjective, hence an
isomorphism (both have dimension $2$). Since $D(x)=x$, $x$ has weight one, so $\hat\varphi^{-1}(P)$ is the Kummer torsor
of $X_P$: $\delta_{\hat\varphi}(P)=[X_P]$.

If $P=T$, i.e.\ $X_P=0$, $Y_P=b$, the relations become $x^2=0$, $y^2=b$, $xy=0$; as $y$ is a unit, $x=0$ in $A$ and
$A=F[y]/(y^2-b)$. Now $D(y)=y+x^2=y$ in $A$, so $y$ has weight one and $\delta_{\hat\varphi}(T)=[b]$.
\end{proof}

\subsection{Cup products and the proof of Lemma~\ref{lem:cup}}
For abelian sheaves (or bounded complexes) $A,B,C$ on a site and a pairing
$A\otimes^{\mathbb L}B\to C$, the cup product is $H^p(A)\otimes H^q(B)\to H^{p+q}(A\otimes^{\mathbb L}B)\to H^{p+q}(C)$,
the first map being induced by the canonical map $R\Gamma(A)\otimes^{\mathbb L}R\Gamma(B)\to R\Gamma(A\otimes^{\mathbb L}B)$.

\smallskip\noindent\emph{\'Etale--fppf comparison.} Let $\epsilon$ be the morphism from the big fppf site to the big \'etale site of
$\Spec F$; a commutative group scheme $G$ over $F$ is a sheaf on both, and $\epsilon^*G=G$. The comparison maps are induced by the unit $X\to R\epsilon_*\epsilon^*X$ of the adjunction. Since $\epsilon^*$ is
exact and compatible with tensor products, $R\epsilon_*$ is lax monoidal and the unit is compatible with the products
$R\Gamma(A)\otimes^{\mathbb L}R\Gamma(B)\to R\Gamma(A\otimes^{\mathbb L}B)$ on both sites; hence the comparison maps
$H^n_{\mathrm{\acute et}}(F,X)\to H^n_\fl(F,\epsilon^*X)$ for bounded complexes $X$ of \'etale sheaves commute with cup
products and with the connecting maps of distinguished triangles. For a smooth commutative group scheme $G$ one has
$R\epsilon_*G=G$ \cite[Thm.~11.7]{Gro}, so the comparison maps are isomorphisms for bounded complexes of smooth
commutative group schemes (by d\'evissage), in particular for $\Z$, $\Z/2$, $\mathbb G_m$; and \'etale cohomology over
$\Spec F$ (on the big or the small site) is Galois cohomology of $G_F$, with the usual cochain cup products.

\begin{proof}
\emph{Models.} Let $P=[\Z\xrightarrow{2}\Z]$ in degrees $-1,0$ and $Q=[\mathbb G_m\xrightarrow{2}\mathbb G_m]$ in
degrees $0,1$. The projection $P\to\Z/2$ (in degree $0$) and the inclusion $\mu_2\to Q$ (in degree $0$) are
quasi-isomorphisms of fppf sheaves, the latter because squaring is fppf surjective with kernel $\mu_2$. Define
$\mu\colon P\otimes Q\to Q$ by
\[
\mu(n\otimes g)=g^n\ \text{ on }P^0\otimes Q^0\text{ and }P^0\otimes Q^1,\qquad
\mu(m\otimes g')=g'^{\,m}\ \text{ on }P^{-1}\otimes Q^1,
\]
and $\mu=0$ on $P^{-1}\otimes Q^0$.
It is a map of complexes: for $m\otimes g\in P^{-1}\otimes Q^0$, $d(m\otimes g)=2m\otimes g-m\otimes g^2$ is sent to
$g^{2m}(g^2)^{-m}=1$; for $n\otimes g\in P^0\otimes Q^0$, $\mu(d(n\otimes g))=\mu(n\otimes g^2)=g^{2n}=d(\mu(n\otimes g))$;
for $m\otimes g'\in P^{-1}\otimes Q^1$, $\mu(d(m\otimes g'))=\mu(2m\otimes g')=g'^{2m}=d(\mu(m\otimes g'))$. On $H^0$ it
induces $(n,\zeta)\mapsto\zeta^n$, the Weil pairing. Since $P$ consists of flat sheaves, $P\otimes Q$ represents
$\Z/2\otimes^{\mathbb L}\mu_2$; its $H^1$ is $\mathbb G_m/\mathbb G_m^2=0$, so $P\otimes Q$ lies in $D^{\leq0}$, and a
morphism $P\otimes Q\to\mu_2[0]$ in the derived category is determined by its effect on $H^0$. Hence $\mu$ represents
the pairing $\Z/2\otimes^{\mathbb L}\mu_2\to\mu_2$. All terms of $P$ and $Q$ are smooth, so by the comparison above the
fppf cup product is computed by Galois hypercochains: $\mathrm{Tot}\,C^\bullet(G_F,P)\otimes
\mathrm{Tot}\,C^\bullet(G_F,Q)\to\mathrm{Tot}\,C^\bullet(G_F,Q)$, cup product of cochains followed by $\mu$.

\emph{Computation.} Let $\tilde\chi\colon G_F\to\{0,1\}\subset\Z$ lift $\chi$ and $Y(\sigma,\tau)=\chi(\sigma)\chi(\tau)$, so
$d\tilde\chi=2Y$ and $dY=0$. In degree $1$, $\mathrm{Tot}\,C^\bullet(G_F,P)$ is $C^1(G_F,P^0)\oplus C^2(G_F,P^{-1})$, and
$(\tilde\chi,\pm Y)$ (the sign depends on the convention for the total differential) is a cocycle mapping to $\chi$.
Under $\mu_2\simeq Q=\mathrm{cone}(\mathbb G_m\xrightarrow2\mathbb G_m)[-1]$, the connecting map of the Kummer sequence is
induced by the inclusion of $\mathbb G_m$ as the degree-$1$ term of $Q$, so $\delta_{\mathrm{Kum}}(u)$ is represented by
the cocycle $(1,u)\in C^1(G_F,Q^0)\oplus C^0(G_F,Q^1)$. Their product has components
$\mu(\pm Y\cup u)=u^{\pm Y}\in C^2(G_F,Q^0)$ and $\mu(\tilde\chi\cup u)=u^{\tilde\chi}\in C^1(G_F,Q^1)$ (the products
with the component $1$ vanish). The map $i\colon\mu_2\to\mathbb G_m$ is the projection $Q\to\mathbb G_m[0]$ onto the
degree-$0$ term, which sends this hypercocycle to the $2$-cocycle $u^{\pm Y}$, i.e.\ to $(\chi,u)^{\pm1}$ by
\eqref{eq:cycliccocycle}.
\end{proof}

\section{Comparison with the flat pairing on \texorpdfstring{$\Sha(E/K)[\varphi]$}{Sha(E/K)[phi]}}\label{app:milne}
We show that on $\Sha(E/K)[\varphi]$ the homogeneous-space pairing coincides with the pairing induced by the flat
duality of \cite[III, \S9]{Mi} over any sufficiently small open subscheme of the curve. The argument is formal, given
Corollary~\ref{cor:Aprime}. Notation is that of \S\ref{sec:W} and \S\ref{sec:radical}. All classes that are paired
below are killed by $2$, so signs, the order of the two arguments (Corollary~\ref{cor:bilin}), and the choice between
an identification $E\simeq E^t$ and its negative (cf.\ convention (i) of \S\ref{sec:ell}) play no role and are ignored.

\subsection{Cohomology with compact support}\label{subsec:Hc}
Let $C$ be the smooth projective curve with function field $K$, $U\subsetneq C$ a non-empty open subscheme, and
$D(U)$ the derived category of sheaves of abelian groups on the flat site of $U$. We use cohomology with compact
support in the form given by Demarche and Harari \cite[\S2]{DH}, who make precise the definition of
\cite[III, \S0 and \S8]{Mi}; as explained in \cite[\S1]{DH}, the cone below has to be formed with complexes, and the
completions $K_v$, not the henselisations, have to be used for flat sheaves such as $\mu_2$. Let
$i\colon Z'=\bigsqcup_{v\notin U}\Spec K_v\to U$. For a sheaf $\mathcal F$ let $I^\bullet(\mathcal F)$ be an injective
resolution; $i^*$ is exact and $i^*I^\bullet(\mathcal F)$ is an injective resolution of $i^*\mathcal F$ \cite[\S2]{DH}.
Put
\[\Gamma_c(U,I^\bullet(\mathcal F))=\mathrm{cone}\bigl(\Gamma(U,I^\bullet(\mathcal F))\to\Gamma(Z',i^*I^\bullet(\mathcal F))\bigr)[-1],
\qquad H^r_c(U,\mathcal F)=H^r\bigl(\Gamma_c(U,I^\bullet(\mathcal F))\bigr).\]
Thus an $n$-cochain is a pair $(\alpha,(\beta_v)_{v\notin U})$ with $\alpha\in\Gamma(U,I^n(\mathcal F))$ and
$\beta_v\in\Gamma(K_v,i_v^*I^{n-1}(\mathcal F))$, and $d(\alpha,(\beta_v))=(d\alpha,(\alpha_v-d\beta_v))$, where
$\alpha_v$ is the restriction of $\alpha$. A homotopy equivalence between two injective resolutions, and a homotopy
between two such, commute with restriction to $Z'$ and so induce the same on cones; hence $H^r_c(U,\mathcal F)$ does
not depend on the resolution.
\begin{enumerate}[leftmargin=*,label=(C\arabic*),ref=(C\arabic*)]
\item\label{C1} There is an exact sequence \cite[Prop.~2.1(1)]{DH}
\[\cdots\to\textstyle\bigoplus_{v\notin U}H^{r-1}(K_v,\mathcal F)\xrightarrow{\partial_c}H^r_c(U,\mathcal F)\to
H^r(U,\mathcal F)\to\bigoplus_{v\notin U}H^r(K_v,\mathcal F)\to\cdots,\]
with maps induced by $(\beta_v)\mapsto(0,(\beta_v))$ and $(\alpha,(\beta_v))\mapsto\alpha$.
\item\label{C2} A short exact sequence $0\to\mathcal M'\to\mathcal M\to\mathcal M''\to0$ gives a long exact sequence
of the groups $H^r_c(U,-)$ \cite[Prop.~2.1(2)]{DH}. We use the following construction of it. By the horseshoe lemma
there are injective resolutions forming an exact sequence $0\to I'^\bullet\to I^\bullet\to I''^\bullet\to0$; it is
split in each degree, so it stays exact after applying $\Gamma(U,-)$ and $\Gamma(Z',i^*-)$, and we get a short exact
sequence of complexes $0\to\Gamma_c(U,I'^\bullet)\to\Gamma_c(U,I^\bullet)\to\Gamma_c(U,I''^\bullet)\to0$, compatible with
the maps of \ref{C1}.
\item\label{C4} $H^3(U,\mathbb G_m)=0$ and $H^2(U,\mathbb G_m)=\Br U$ \cite[II, Prop.~2.1]{Mi} (flat and \'etale
cohomology agree for the smooth group $\mathbb G_m$), so by \ref{C1} $H^3_c(U,\mathbb G_m)$ is the cokernel of
$\Br U\to\bigoplus_{v\notin U}\Br K_v$. As $\Br U$ is the group of classes in $\Br K$ unramified on $U$ and
$U\neq C$, class field theory (Fact~\ref{fact:ABHN} and the exactness of $0\to\Br K\to\bigoplus_v\Br K_v\to\Q/\Z\to0$)
shows that $\sum_{v\notin U}\inv_v$ induces an isomorphism $\operatorname{tr}\colon H^3_c(U,\mathbb G_m)\to\Q/\Z$ with
$\operatorname{tr}\circ\,\partial_c=\sum_{v\notin U}\inv_v$.
\end{enumerate}

\begin{lemma}\label{lem:C3}
Let $e\colon\mathcal F\to\mathcal G[r]$ be a morphism in $D(U)$, i.e.\ an element of
$\operatorname{Ext}^r_U(\mathcal F,\mathcal G)$.
\begin{enumerate}[label=\textup{(\roman*)},leftmargin=*]
\item $e$ induces homomorphisms $e_*\colon H^s_c(U,\mathcal F)\to H^{s+r}_c(U,\mathcal G)$, and $(e'\circ e)_*=e'_*\circ e_*$.
\item $e_*$ is compatible with the maps of \ref{C1}: with $e_*$ on $H^s(U,-)$, and
$e_*\circ\partial_c=\partial_c\circ\bigl((e_v)_*\bigr)_v$, where $e_v=i_v^*e$ and
$(e_v)_*\colon H^{s-1}(K_v,\mathcal F)\to H^{s+r-1}(K_v,\mathcal G)$.
\item If $e\in\operatorname{Ext}^1_U(\mathcal M'',\mathcal M')$ is the class of the extension in \ref{C2}, then $e_*$ is
the connecting map of the long exact sequence of \ref{C2}.
\end{enumerate}
\end{lemma}

\begin{proof}
(i), (ii) The morphism $e$ is represented by a map of complexes
$\epsilon\colon I^\bullet(\mathcal F)\to I^\bullet(\mathcal G)[r]$, unique up to homotopy. The maps $\Gamma(U,\epsilon)$
and $\Gamma(Z',i^*\epsilon)$ commute with restriction, so they define a map of cones; a homotopy $h$ between two choices
gives the homotopies $\Gamma(U,h)$, $\Gamma(Z',i^*h)$, which again commute with restriction and give a homotopy of
cones. This defines $e_*$; functoriality and (ii) are read off the definition, since $i_v^*\epsilon$ represents $e_v$.
(iii) With the resolutions of \ref{C2}, let $s\colon I''^\bullet\to I^\bullet$ be a splitting in each degree. Then
$\epsilon=ds-sd$ has values in $I'^\bullet$ and represents $e$. The splittings $\Gamma(U,s)$ and $\Gamma(Z',i^*s)$ commute
with restriction, so the connecting map of the exact sequence of cones in \ref{C2}, computed by lifting a cocycle by
$s$ and applying $d$, is induced by $\epsilon$.
\end{proof}

For $\chi\in H^1(U,\Z/2)=\Hom_{D(U)}(\Z,\Z/2[1])$ we write $\chi\cup{-}\colon\mu_2\to\mathbb G_m[1]$ for the composite of
$\chi\otimes1\colon\mu_2=\Z\otimes^L\mu_2\to\Z/2[1]\otimes^L\mu_2$ with the Cartier pairing
$\ev\colon\Z/2\otimes^L\mu_2\to\Z/2\otimes\mu_2\to\mathbb G_m$. Over $K_v$ it induces the cup product
$g\mapsto\chi_v\cup g$, $H^1_\fl(K_v,\mu_2)\to\Br K_v$, of Corollary~\ref{cor:Aprime}: the cup product of
Appendix~\ref{app:flat} with a class $x\in H^1(F,\Z/2)=\Hom(\Z,\Z/2[1])$ is the map induced by $x\otimes1$ followed by
the pairing.

\subsection{The flat pairing over \texorpdfstring{$U$}{U}}\label{subsec:milnepairing}
Let $\mathcal E$ be an abelian scheme over $U$ with generic fibre $E$; we identify $\mathcal E$ with its dual by the
canonical principal polarisation. The Poincar\'e biextension defines a morphism
$W=W_{\mathcal E}\colon\mathcal E\otimes^L\mathcal E\to\mathbb G_m[1]$ in $D(U)$ \cite[III, \S9 and Appendix~C]{Mi}; the
adjoint morphism $\mathcal E\to R\mathcal Hom(\mathcal E,\mathbb G_m)[1]$ induces on $\mathcal H^0$ the
Barsotti--Weil isomorphism $\mathcal E\simeq\mathcal Ext^1_U(\mathcal E,\mathbb G_m)$ for abelian schemes, and
$\mathcal Hom_U(\mathcal E,\mathbb G_m)=0$ \cite[III, \S18]{Oo}. For $a\in H^1(U,\mathcal E)=\Hom_{D(U)}(\Z,\mathcal E[1])$
put $e(a)=W\circ(a\otimes1)\colon\mathcal E\to\mathbb G_m[2]$, and for $c\in H^1_c(U,\mathcal E)$
\[\langle a,c\rangle_U=\operatorname{tr}\bigl(e(a)_*c\bigr)\in\Q/\Z,\]
with $e(a)_*$ as in Lemma~\ref{lem:C3}. The class $e(a)\in\operatorname{Ext}^2_U(\mathcal E,\mathbb G_m)$ is the image
of $a$ under the Barsotti--Weil isomorphism and the edge map
$H^1(U,\mathcal Ext^1_U(\mathcal E,\mathbb G_m))\to\operatorname{Ext}^2_U(\mathcal E,\mathbb G_m)$, so this is the pairing
$H^1(U,\mathcal E)\times H^1_c(U,\mathcal E)\to\Q/\Z$ that underlies \cite[III, Thm.~9.4]{Mi}, which states that its
kernels, for torsion classes $c$, are the divisible subgroups; \cite[III, Cor.~9.5]{Mi} is deduced from that theorem.
For the duality theorems in this setting see also \cite{DH,Su}.

The map $H^1(U,\mathcal E)\to H^1(K,E)$ is injective, and $\Sha(E/K)$ is the kernel of
$H^1(U,\mathcal E)\to\bigoplus_{v\notin U}H^1(K_v,E)$ \cite[II, Lemma~5.5]{Mi}.
For $a,a'\in\Sha(E/K)$, regarded in $H^1(U,\mathcal E)$, choose $a'_c\in H^1_c(U,\mathcal E)$ mapping to $a'$, which is
possible by \ref{C1}, and put
\[\langle a,a'\rangle_{\fl,U}=\langle a,a'_c\rangle_U .\]
This does not depend on $a'_c$: two choices differ by $\partial_c(s)$ with $s\in\bigoplus_{v\notin U}E(K_v)$, and by
Lemma~\ref{lem:C3}(ii) $e(a)_*\partial_c(s)=\partial_c\bigl((e(a_v)_*s_v)_v\bigr)=0$, because $a_v=0$. We call
$\langle\ ,\ \rangle_{\fl,U}$ the \emph{flat pairing over $U$}. We do not prove that it is independent of $U$ on all of
$\Sha(E/K)$, nor compare it further with the pairing of \cite[III, Cor.~9.5]{Mi}; the next proposition shows that on
$\Sha(E/K)[\varphi]$ it is the same for all sufficiently small $U$.

\begin{proposition}\label{prop:milne}
There is a non-empty open $U_0\subsetneq C$ such that for every non-empty open $U\subseteq U_0$ and all
$\chi,\chi'\in\Sel^\varphi(E/K)$
\[\langle\iota(\chi),\iota(\chi')\rangle_{\fl,U}=\sum_v\inv_v(\chi_v\cup g_v)=\langle\iota(\chi),\iota(\chi')\rangle,\]
where $g_v\in K_v^\times/K_v^{\times2}$ are the classes of Theorem~\textup{\ref{thm:W}(b)} attached to $\chi'$ and the
last term is the homogeneous-space pairing of Definition~\textup{\ref{def:PS}}. Thus, for all such $U$, the flat
pairing over $U$ and the homogeneous-space pairing coincide on $\Sha(E/K)[\varphi]\times\Sha(E/K)[\varphi]$.
\end{proposition}

The second equality is Corollary~\ref{cor:Aprime}. For the first, fix $\chi,\chi'$ and $\beta_0$; we show that the equality holds for every $U$ satisfying
\ref{U1}--\ref{U3} below. As these conditions hold for all sufficiently small $U$ and $\Sel^\varphi(E/K)$ is finite
(Lemma~\ref{lem:ratio}), there is $U_0$ as in the proposition. Let $\xi\in H^1_\fl(K,E'[2])$
be the lift $\Phi^{-1}(\chi',[\beta_0])$ of $\chi'$ and $\Sigma$ the finite set of \S\ref{subsec:A6}. We require:
\begin{enumerate}[leftmargin=*,label=(U\arabic*),ref=(U\arabic*)]
\item\label{U1} $a_2$ is integral and $b$ is a unit at every $v\in U$. Then the equations \eqref{eq:E} of $E$ and $E'$
define abelian schemes $\mathcal E$, $\mathcal E'$ over $U$ (their discriminants are $b^2$ and $b$), the section
$T=(0,b)$ is disjoint from the zero section, $\varphi$ and $\hat\varphi$ extend to isogenies, and we have exact sequences
of sheaves on the flat site of $U$
\[0\to\Z/2\xrightarrow{\,i\,}\mathcal E\xrightarrow{\varphi}\mathcal E'\to0,\qquad
0\to\mu_2\xrightarrow{\,j'\,}\mathcal E'\xrightarrow{\hat\varphi}\mathcal E\to0,\]
the kernel of $\hat\varphi$ being the Cartier dual of that of $\varphi$; as $\operatorname{Aut}_U(\mu_2)=1$, there is
no ambiguity in the identification with $\mu_2$. Here $j'$ is $j\colon\mu_2\to\mathcal E'[2]$ followed by the inclusion
$\mathcal E'[2]\subset\mathcal E'$.
\item\label{U2} $\xi$ is the restriction of a class in $H^1_\fl(U,\mathcal E'[2])$, again denoted $\xi$. This holds for
small $U$ since $\mathcal E'[2]$ is finite and locally free over $U$, so that
$H^1_\fl(K,E'[2])=\varinjlim_UH^1_\fl(U,\mathcal E'[2])$.
\item\label{U3} $U\cap\Sigma=\emptyset$.
\end{enumerate}
Condition \ref{U3} implies \ref{U1}, by the definition of $\Sigma$, and that $\chi$ is unramified on $U$, so
$\chi\in H^1(U,\Z/2)$. Put $a=i_*\chi\in H^1(U,\mathcal E)$, let
$y\in H^1(U,\mathcal E')$ be the image of $\xi$, and $a'=\hat\varphi_*y\in H^1(U,\mathcal E)$. On the generic fibre, $a$
restricts to $\iota(\chi)$ and, by \eqref{eq:compat} and $\hat\varphi_*\xi=\chi'$, $a'$ restricts to $\iota(\chi')$; so
$a$ and $a'$ are the elements of $H^1(U,\mathcal E)$ corresponding to $\iota(\chi),\iota(\chi')\in\Sha(E/K)$. Moreover,
for $v\notin U$,
\begin{equation}\label{eq:ylocal}
y_v=j'_*(g_v)\quad\text{in }H^1(K_v,E'),
\end{equation}
because $\xi_v-\delta_2(R'_v)=j_*(g_v)$ in $H^1_\fl(K_v,E'[2])$ by Corollary~\ref{cor:Aprime}, and $\delta_2(R'_v)$ maps
to $0$ in $H^1(K_v,E')$.

\begin{lemma}\label{lem:clift}
Let $0\to\mathcal M'\to\mathcal M\to\mathcal M''\to0$ be an exact sequence of sheaves on the flat site of $U$, let
$y\in H^1(U,\mathcal M)$ have image $a'\in H^1(U,\mathcal M'')$, and suppose that for every $v\notin U$ there is
$h_v\in H^1(K_v,\mathcal M')$ with image $y_v$ in $H^1(K_v,\mathcal M)$. Then there is $a'_c\in H^1_c(U,\mathcal M'')$
mapping to $a'$ such that the connecting map $\partial\colon H^1_c(U,\mathcal M'')\to H^2_c(U,\mathcal M')$ of
\ref{C2} satisfies $\partial(a'_c)=\partial_c\bigl((h_v)_{v\notin U}\bigr)$.
\end{lemma}

\begin{proof}
We use the resolutions and the exact sequences of complexes of \ref{C2}. Let $\tilde y\in\Gamma(U,I^1)$ and
$\tilde h_v\in\Gamma(K_v,i_v^*I'^1)$ be cocycles representing $y$ and $h_v$; then $\tilde y_v=\tilde h_v+dm_v$ with
$m_v\in\Gamma(K_v,i_v^*I^0)$. Let $\bar y$,
$\bar m_v$ be the images in the complexes of $\mathcal M''$. As $\tilde h_v$ maps to $0$ there, $\bar y_v=d\bar m_v$, so
$(\bar y,(\bar m_v))$ is a $1$-cocycle of the complex computing $H_c(U,\mathcal M'')$, and its class $a'_c$ maps to $a'$.
To compute $\partial(a'_c)$ we lift this cocycle to the cochain $(\tilde y,(m_v))$ of the complex of $\mathcal M$; its
differential is $(d\tilde y,(\tilde y_v-dm_v))=(0,(\tilde h_v))$, which lies in the complex of $\mathcal M'$ and represents
$\partial_c((h_v))$.
\end{proof}

\begin{lemma}\label{lem:cartier}
Assume \ref{U1}, and let $\partial_{\hat\varphi}\in\Hom_{D(U)}(\mathcal E,\mu_2[1])=\operatorname{Ext}^1_U(\mathcal E,\mu_2)$
be the class of the extension $0\to\mu_2\to\mathcal E'\xrightarrow{\hat\varphi}\mathcal E\to0$. Then
\[W_{\mathcal E}\circ(i\otimes1)=\ev\circ(1\otimes\partial_{\hat\varphi})\colon\ \Z/2\otimes^L\mathcal E\to\mathbb G_m[1].\]
\end{lemma}

\begin{proof}
On the flat site the squaring map of $\mathbb G_m$ is an epimorphism, so
$R\mathcal Hom(\Z/2,\mathbb G_m)=\mathcal Hom(\Z/2,\mathbb G_m)=\mu_2$, and by adjunction
\[\Hom_{D(U)}(\Z/2\otimes^L\mathcal E,\mathbb G_m[1])=\Hom_{D(U)}(\mathcal E,\mu_2[1]),\qquad
\ev\circ(1\otimes\rho)\leftrightarrow\rho .\]
Let $\rho$ correspond to $W_{\mathcal E}\circ(i\otimes1)$; we have to show $\rho=\partial_{\hat\varphi}$. The isogenies
$\varphi$ and $\hat\varphi$ are dual to each other under the canonical polarisations: this holds on the generic fibre
\cite[III, \S6]{Si}, and two homomorphisms of abelian schemes over $U$ which agree on the generic fibre agree. The dual
homomorphism is defined by pulling back the Poincar\'e biextension, so
$W_{\mathcal E}\circ(1\otimes\hat\varphi)=W_{\mathcal E'}\circ(\varphi\otimes1)$ on $\mathcal E\otimes^L\mathcal E'$. Hence
$W_{\mathcal E}\circ(i\otimes\hat\varphi)=W_{\mathcal E'}\circ(\varphi i\otimes1)=0$, i.e.\ $\rho\circ\hat\varphi=0$. By the
exact sequence
\[\Hom_{D(U)}(\mu_2[1],\mu_2[1])\xrightarrow{\ \circ\,\partial_{\hat\varphi}\ }\Hom_{D(U)}(\mathcal E,\mu_2[1])
\xrightarrow{\ \circ\,\hat\varphi\ }\Hom_{D(U)}(\mathcal E',\mu_2[1])\]
of the triangle $\mathcal E'\xrightarrow{\hat\varphi}\mathcal E\xrightarrow{\partial_{\hat\varphi}}\mu_2[1]\to\mathcal E'[1]$,
we get $\rho=\theta\circ\partial_{\hat\varphi}$ with $\theta\in\operatorname{End}_U(\mu_2)$. By Cartier duality
$\operatorname{End}_U(\mu_2)=\operatorname{End}_U(\Z/2)=\Z/2$, as $U$ is connected; so $\theta$ is $0$ or the identity,
and it remains to see that $W_{\mathcal E}\circ(i\otimes1)\neq0$. Its adjoint in the other variable is the composite of
$i\colon\Z/2\to\mathcal E$ with $\mathcal E\to R\mathcal Hom(\mathcal E,\mathbb G_m)[1]$. Since
$\mathcal Hom_U(\mathcal E,\mathbb G_m)=0$, morphisms in $D(U)$ from a sheaf to $R\mathcal Hom(\mathcal E,\mathbb G_m)[1]$
are the homomorphisms to $\mathcal Ext^1_U(\mathcal E,\mathbb G_m)$, and the composite is the Barsotti--Weil isomorphism
composed with the monomorphism $i$, which is not $0$.
\end{proof}

\begin{proofof}{Proposition~\textup{\ref{prop:milne}}}
Let $U$ satisfy \ref{U1}--\ref{U3}. By \eqref{eq:ylocal}, Lemma~\ref{lem:clift} applies to
$0\to\mu_2\to\mathcal E'\xrightarrow{\hat\varphi}\mathcal E\to0$, the class $y$ and $h_v=g_v$: there is
$a'_c\in H^1_c(U,\mathcal E)$ mapping to $a'$ with $\partial(a'_c)=\partial_c((g_v)_{v\notin U})$ in $H^2_c(U,\mu_2)$. By
Lemma~\ref{lem:cartier},
\[e(a)=W_{\mathcal E}\circ(i\otimes1)\circ(\chi\otimes1)=\ev\circ(1\otimes\partial_{\hat\varphi})\circ(\chi\otimes1)
=(\chi\cup{-})\circ\partial_{\hat\varphi}\colon\ \mathcal E\to\mu_2[1]\to\mathbb G_m[2].\]
Therefore, by Lemma~\ref{lem:C3} and \ref{C4},
\[\langle\iota(\chi),\iota(\chi')\rangle_{\fl,U}=\operatorname{tr}\bigl((\chi\cup{-})_*\,\partial(a'_c)\bigr)
=\operatorname{tr}\bigl((\chi\cup{-})_*\,\partial_c((g_v))\bigr)=\operatorname{tr}\,\partial_c\bigl((\chi_v\cup g_v)_v\bigr)
=\sum_{v\notin U}\inv_v(\chi_v\cup g_v).\]
By \ref{U3} and \S\ref{subsec:A6} the terms with $v\in U$ in $\sum_v\inv_v(\chi_v\cup g_v)$ vanish.
\end{proofof}

\begin{remark}\label{rem:milnescope}
(1) Another choice of $h_v$ in the proof differs from $g_v$ by an element of $\kappa(E(K_v))$, which is orthogonal to
$\chi_v$ by Lemma~\ref{lem:iso}, in accordance with the independence of $a'_c$.
(2) Proposition~\ref{prop:milne} depends on Corollary~\ref{cor:Aprime}, hence on the computation of \S\ref{sec:W}.
Theorem~\ref{thm:BM} and Theorem~\ref{thm:radical} do not depend on it.
(3) The comparison is proved on $\Sha(E/K)[\varphi]$ only. We do not know whether the two pairings agree on all of
$\Sha(A/K)[2^\infty]$ for an ordinary elliptic curve $A$ over $K$.
\end{remark}

\end{document}